\documentclass[10pt]{amsart}
\usepackage{pstricks,pst-plot,palatino,latexsym}
\usepackage{psfrag,graphicx}
\usepackage{subfigure}
\usepackage{amsmath,amsthm,epsf,amscd,amssymb,verbatim}

\newcommand{\ds}{\displaystyle}
\newcommand{\R}{{\mathbb{R}}}
\newcommand{\real}{{\mathbb{R}}}
\newcommand{\N}{{\mathbb{N}}}
\newcommand{\Z}{{\mathbb{Z}}}

\newcommand{\D}{{\mathbb{D}}}
\newcommand{\E}{{\mathbb{E}}}
\newcommand{\Q}{{\mathbb{Q}}}

\newcommand{\Id}{{\mbox{{\sc Id}}}}
\newcommand{\uu}{{\bf u}}
\newcommand{\ww}{{\bf w}}
\newcommand{\rmref}[1]{{\rm(\ref{#1})}}
\newcommand{\word}[1]{\left\vert{#1}\right\vert_{\mbox{\tiny word}}}
\newcommand{\fp}{\hspace{1mm}\rule{2mm}{2mm}}
\newcommand{\vsp}{\vspace{.1in}}
\newcommand{\vv}{{\bf v}}

\newcommand{\B}{{\mathcal B}}

\newcommand{\RR}{{\mathcal R}}
\newcommand{\DD}{{\mathcal D}}
\newcommand{\EE}{{\mathcal E}}

\newcommand{\NN}{{\mathcal B}}
\newcommand{\MM}{{\mathcal M}}

\newcommand{\U}{{\mathcal U}_E}

\newcommand{\hh}{{\mbox{\bf H}}}
\newcommand{\x}{{\bf x}}
\newcommand{\y}{{\bf y}}
\renewcommand{\a}{{\bf a}}
\renewcommand{\b}{{\bf b}}

\newcommand{\floor}[1]{{\left\lfloor{#1}\right\rfloor}}
\newcommand{\ceil}[1]{{\left\lceil{#1}\right\rceil}}

\newcommand{\braid}{{\beta}}
\newcommand{\Conf}{{\mathcal D}}

\newcommand{\Collapsed}{{\Sigma^{-}}}
\newcommand{\eg}{{\it e.g.}}
\newcommand{\ie}{{\it i.e.}}
\newcommand{\Inv}{{\mbox{\sc Inv}}}
\newcommand{\Seq}{{\bf X}}
\newcommand{\sign}{{\mbox{\sc sign}}}

\newcommand{\Dual}{{\mathbb D}}	
\newcommand{\Sym}{{\mbox{\sc sym}}}
\newcommand{\rel}{{\mbox{\sc rel}}}
\newcommand{\intnum}{\iota}
\newcommand{\pf}{{\em Proof: }}

\newcommand{\eps}{\epsilon}
\newgray{grey1}{0.75}\newgray{grey2}{0.6}\newgray{grey3}{0.45}

\newtheorem{theorem}{Theorem}
\newtheorem{lemma}[theorem]{Lemma}
\newtheorem{proposition}[theorem]{Proposition}
\newtheorem{definition}[theorem]{Definition}
\newtheorem{corollary}[theorem]{Corollary}
\newtheorem{remark}[theorem]{Remark}

\begin{document}
\title[Morse theory on
spaces of braids]{\textsc Morse theory on
spaces of braids and~Lagrangian~dynamics}
\author{R.W. Ghrist}
\address{Department of Mathematics, University of Illinois, Urbana IL, 
        61801 USA}
\author{J.B. Van den Berg}
\address{Department of Applied Mathematics, University of Nottingham, UK }
\author{R.C. VanderVorst}
\address{Department of Mathematics, Free University Amsterdam, 
  De Boelelaan  1081, Amsterdam
  Netherlands; and CDSNS, Georgia Institute of Technology, 
  Atlanta GA, 30332-0160 USA}
\date{\today}
\thanks{The first author was supported by NSF DMS-9971629 and 
  NSF DMS-0134408. 
The second author was supported by an EPSRC Fellowship.
The third author was supported by NWO Vidi-grant 639.032.202.   }

\begin{abstract}
In the first half of the paper we construct a Morse-type theory on certain
spaces of braid diagrams. We define a topological invariant of closed
positive braids which is correlated with the existence of invariant sets 
of {\it parabolic flows} defined on discretized braid spaces.
Parabolic flows, a type of one-dimensional lattice dynamics, 
evolve singular braid diagrams in such a way as
to decrease their topological complexity; algebraic lengths  
decrease monotonically. This topological invariant is derived from  
a Morse-Conley homotopy index.

In the second half of the paper we apply this technology to 
second order Lagrangians via a discrete formulation of the 
variational problem. This culminates in a very general forcing theorem 
for the existence of infinitely many braid classes of closed orbits. 
\end{abstract}

\maketitle
\begin{sloppypar}

\setcounter{tocdepth}{1}
\setcounter{secnumdepth}{3}

%
%
\section{Prelude}
\label{I}\label{prelude}

It is well-known that under the evolution of any scalar uniformly
parabolic equation of the form 
\begin{equation}
\label{PDE}
u_{t} = f(x,u,u_x,u_{xx}) \quad ; \quad 
\partial_{u_{xx}}f\geq\delta>0 ,	
\end{equation} 
the graphs of two solutions $u_1(x,t)$ and $u_2(x,t)$ evolve in such a way 
that the number of intersections of the graphs does not increase
in time. This principle, known in various circles as ``comparison 
principle'' or ``lap 
number'' techniques, entwines the geometry of the graphs
($u_{xx}$ is a curvature term), the topology of the solutions 
(the intersection number is a local linking number), and the 
local dynamics of the PDE. This is a valuable approach 
for understanding local dynamics for a wide variety of flows 
exhibiting parabolic behavior with both classical \cite{Stu30}
and contemporary \cite{Mat1,Angenent4,BF,FR99} implications.

This paper is an extension of this local technique to a global
technique. One such well-established globalization appears in
the work of Angenent on curve-shortening \cite{Angenent3}: evolving 
closed curves on a surface by curve shortening isolates the 
classes of curves dynamically and implies a monotonicity with 
respect to number of self-intersections. 

In contrast, one could consider the following topological globalization.
Superimposing the graphs of a collection of functions 
$u^\alpha (x)$ gives something which resembles the projection of a 
topological braid onto the plane. Assume that the ``height'' of the
strands above the page is given by the slope $u^\alpha_x(x)$, or, 
equivalently, that all of the crossings in the projection are of the 
same sign (bottom-over-top): see Fig.~\ref{fig_ex}[left]. 
Evolving these functions under a parabolic equation (with, say, boundary 
endpoints fixed) yields a flow on a certain space of braid diagrams 
which has a topological monotonicity: {\em linking can be 
destroyed but not created.} This establishes
a partial ordering on the semigroup of positive braids which 
is respected by parabolic dynamics. The idea of topological braid classes 
with this partial ordering is a globalization of the lap number 
(which, in braid-theoretic terms becomes the length of the braid 
in the braid group under standard generators). 

\subsection{Parabolic flows on spaces of braid diagrams.}\label{I1}
\begin{figure}[hbt]
 \begin{center}
 \psfragscanon
 \psfrag{x}[][]{$x$}
 \psfrag{u}[][]{$u$}
 \psfrag{d}[][]{$u_{x}$}
 \includegraphics[angle=0, height=1.0in,
                          width=5.2in]{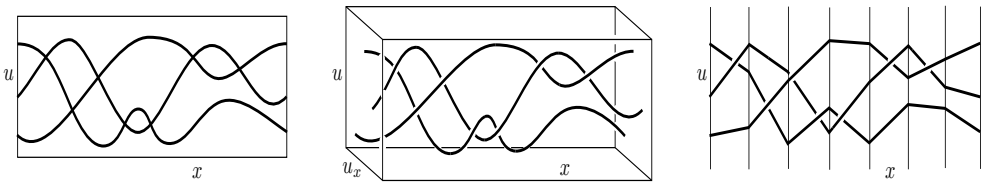}
 \caption{Curves in the $x-u$ plane [left] lift to a braid [center] 
 which is then discretized [right].
 In a discretized isotopy, one slides the anchor points vertically.
 }
 \label{fig_ex}
 \end{center}
 \end{figure}

In this paper, we initiate the study of parabolic flows 
on spaces of braid diagrams. The particular braids in question 
will be (a) {\it positive} -- all crossings are considered to be
of the same sign; (b) {\it closed}\footnote{The theory works equally well
for braids with fixed endpoints.} -- the left and right sides
are identified; and (c) {\it discretized} -- or piecewise linear
with fixed distance between ``anchor points,'' so as to avoid the analytic 
difficulties of working on infinite dimensional spaces of curves. 
See Fig.~\ref{fig_ex} for examples of braid diagrams. 

 







 


The flows we consider evolve the anchor points of the braid diagram 
so that the braid class can change, but only so as to decrease 
complexity: local linking of 
strands may not increase with time. Due to the close similarity with 
parabolic partial differential equations such systems will be referred 
to as {\it parabolic recurrence relations},
and the induced flows as {\it parabolic flows}. These flows are given by
\begin{eqnarray}\label{DPDE}
	\frac{d}{dt}u_i = \RR_i(u_{i-1},u_i,u_{i+1}),
\end{eqnarray}
where the variables $u_i$ represent the vertical positions of the ordered 
anchor points of discrete braid diagrams. The only conditions imposed 
on the dynamics is the monotonicity condition that every $\RR_i$ be 
increasing functions of $u_{i-1}$ and $u_{i+1}$.

While a discretization of a PDE of the form (\ref{PDE}) with 
nearest-neighbor interaction yields a parabolic recurrence relation, the 
class of dynamics we consider is significantly larger in scope
(see, \eg, \cite{Mallet-Paret}). Parabolic recurrence relations are a
sub-class of monotone recurrence relations as studied in \cite{Angenent2}
and \cite{Hir}.

The evolution of braid diagrams yields a situation not unlike that 
considered by Vassiliev in knot theory \cite{Vassiliev}: in our
scenario, the space of all braid diagrams is partitioned by the 
discriminant of singular diagrams into the braid classes. The parabolic 
flows we consider are transverse to these singular varieties (except 
for a set of ``collapsed'' braids) and are co-oriented
in a direction along which the algebraic length of the braid
decreases: this is an algebraic version of curve shortening.

To proceed, two types of noncompactness on spaces of braid diagrams 
must be repaired. Most severe is the problem of braid strands 
completely collapsing onto one another. To resolve 
this type of noncompactness, we assume that the dynamics fixes
some collection of braid strands, a {\em skeleton}, and then 
work on spaces of braid pairs: one free, one fixed. The 
relative theory then leads to forcing results of the type
``Given a stationary braid class, which other braids are 
forced as invariant sets of parabolic flows?'' The second 
type of noncompactness in the dynamics occurs when the braid
strands are free to evolve to arbitrarily large values. In the 
PDE setting, one requires knowledge of boundary conditions at 
infinity to prove theorems about the dynamics. In our braid-theoretic
context, we convert boundary conditions to ``artificial'' braid
strands augmented to the fixed skeleton.

Thus, working on spaces of braid pairs, the dynamics at the discriminant 
permits the construction of a Morse theory in the spirit of Conley to detect 
invariant sets of parabolic flows. 
Conley's extension of the Morse index associates to any sufficiently 
isolated invariant set a space whose homotopy type measures 
not merely the dimension of the unstable manifold (the Morse index) 
but rather the coarse topological features 
of the unstable dynamics associated to this set. We obtain
a well-defined Conley index for braid diagrams from the 
monotonicity properties of parabolic flows. 
To be more precise, relative braid classes (equivalence classes of isotopic 
braid diagrams fixing some skeleton) serve as candidates for isolating
neighborhoods to which the Conley index can be assigned. 
This approach is reminiscent of the ideas of linking of periodic orbits
used by Angenent \cite{Angenent1,Angenent3} and LeCalvez 
\cite{LeCalvez,LeCalvez2}.

Our finite-dimensional approximations to the (infinite-dimensional)
space of smooth topological braids conceivably alter the 
Morse-theoretic properties of the discretized braid classes. 
One would like to know that so long as the discretization is not degenerately
coarse, the homotopy index is independent of both the discretization
and the specific parabolic flow employed. This is true. 
The principal topological result of this work is that the homotopy index is 
indeed an invariant of the topological (relative) braid class:
see Theorems~\ref{topinvariant} and \ref{stabilization}
for details. These theorems seem to evade a simple algebraic-topological
proof. The proof we employ in \S\ref{stable} constructs the appropriate 
homotopy by recasting the problem into singular dynamics 
and applying techniques from singular perturbation theory. 

We thus obtain a topological index which can, like the Morse index,
force the existence of invariant sets. Specifically, a non-vanishing 
homotopy index for a relative braid class indicates that there is an invariant
set in this braid class for any parabolic flow with the appropriate 
skeleton. This is the foundation for the applications to follow in the
remainder of the paper.


The remainder of the paper explores applications of the machinery 
to a broad class of Lagrangian dynamics.

\subsection{Second order Lagrangian dynamics.}\label{prelude2}
Our principal application of the Morse theory on discretized braids is to 
the problem of finding periodic orbits of second order Lagrangian systems:
that is, Lagrangians of the form $L(u,u_x,u_{xx})$ where $L \in C^2(\R^3)$.
An important motivation for studying such systems comes 
from the stationary {\it Swift-Hohenberg model} in physics, which
is described by the fourth order equation 
\begin{eqnarray}\label{SH}
\Biggl( 1+ \frac{d^2}{dx^2}\Biggr)^2 u -\alpha u + u^3 =0,\quad \alpha\in \R.
\end{eqnarray}
This equation is the Euler-Lagrange equation of the second order Lagrangian 
\begin{equation}
	L(u,u_x,u_{xx})=\frac{1}{2}|u_{xx}|^2-|u_x|^2 +
	\frac{1-\alpha}{2}u^2+\frac{1}{4}u^4 .
\end{equation}

We generalize to the broadest possible class of second order Lagrangians. 
One begins with
the conventional convexity assumption, $\partial^2_w L(u,v,w) \ge \delta >0$.
The objective is to find bounded functions $u:~\R \to \R$
which are stationary for the action integral 
$J[u] := \int L(u,u_x,u_{xx})dx$.
Such functions $u$ are bounded solutions of the Euler-Lagrange equations
\begin{equation}
\label{EL}
	\frac{d^2}{dx^2} \frac{\partial L}{\partial u_{xx}}  -
	\frac{d}{dx} \frac{\partial L}{\partial u_x}  +
 	\frac{\partial L}{\partial u} = 0.
\end{equation}
Due to the translation invariance  $x \mapsto x+c$,
the solutions of \rmref{EL} satisfy the energy constraint
\begin{equation}
\label{eq_L}
\Bigl( \frac{\partial L}{\partial u_x} - \frac{d}{dx} \frac{\partial L}
{\partial u_{xx}} \Bigr) u_x + \frac{\partial L}{\partial u_{xx}}u_{xx} -
L(u,u_x,u_{xx}) = E =\mbox{ constant},
\end{equation}
where $E$ is the energy of a solution. To find bounded solutions 
for given values of $E$, we employ the variational principle 
$\delta_{u,T} \int_0^T \bigl(L(u,u_x,u_{xx})+E\bigr)dx=0$,
which forces solutions of \rmref{EL} to have energy $E$.
The Lagrangian problem can be reformulated as a two degree-of-freedom 
Hamiltonian system; in that context, bounded {\it periodic} solutions are 
{\it closed characteristics} of the (corresponding) energy manifold 
$M^3\subset\R^4$. Unlike the case of first-order Lagrangian systems,
the energy hypersurface is {\it not} of contact type in general 
\cite{AVV}, and the recent stunning results in contact homology 
\cite{EGH00} are inapplicable.

The variational principle can be discretized for a certain considerable
class of second order Lagrangians: those for which monotone laps 
between consecutive extrema $\{u_i\}$ are unique and continuous with 
respect to the endpoints. 
We give a precise definition in \S\ref{second}, denoting these as
(second order Lagrangian) {\it twist systems}.
Due to the energy identity \rmref{eq_L} the extrema $\{u_i\}$ are 
restricted to the set $\U= \{ u~|~ L(u,0,0)+E \ge 0\}$, connected components 
of which are called {\it interval components} and denoted by $I_E$.
An energy level is called regular if 
$\frac{\partial L}{\partial u}(u,0,0) \neq 0$ 
for all $u$ satisfying $L(u,0,0)+E=0$.
In order to deal with non-compact interval components $I_E$ certain
asymptotic behavior has to be specified, for example that 
``infinity'' is attracting.
Such Lagrangians are called {\it dissipative}, and are most common 
in models coming from physics, like the Swift-Hohenberg Lagrangian. For 
a precise definition of dissipativity see \S\ref{multi}.
Other asymptotic behaviors may be considered as well, such as 
``infinity'' is repelling, or
more generally that infinity is isolating, implying that
closed characteristics are a priori bounded in $L^\infty$.

Closed characteristics are either {\em simple} or {\em non-simple}  
depending on whether $u(x)$, represented as a closed curve in the 
$(u,u_x)$-plane, is a simple closed curve or not. This distinction is
a sufficient language for the following general forcing theorem:
\begin{theorem}\label{H}
Any dissipative twist system possessing a non-simple closed characteristic 
$u(x)$ at a regular energy value $E$  such that $u(x) \in I_E$, 
must possess an infinite number of (non-isotopic) closed 
characteristics at the same energy level as $u(x)$.
\end{theorem}
This is the optimal type of forcing result: 
there are neither hidden assumptions about nondegeneracy 
of the orbits, nor restrictions to generic behavior. Sharpness 
comes from the fact that there exist systems with finitely many  
simple closed characteristics at each energy level.

The above result raises the following question: when does an energy manifold
contain a non-simple closed characteristic? In general the existence of 
such characteristics depends on the geometry of the energy manifold.
One geometric property that sparks the existence of non-simple
closed characteristics is a singularity or near-singularity of the energy 
manifold. This, coupled with Theorem~\ref{H}, triggers the existence of 
infinitely many closed characteristics. The results that can be proved 
in this context (dissipative twist systems) give a complete classification
with respect to the existence of finitely many versus infinitely many
closed characteristics on singular energy levels. The first result in 
this direction deals with singular energy values for which $I_E=\R$.
\begin{theorem}\label{H3}
Suppose that a dissipative twist system has a singular energy level $E$ 
with $I_{E}=\R$, which contains two or more rest points.
Then the system has infinitely many closed characteristics at energy 
level $E$.\footnote{From the proof of this theorem
in \S\ref{multi} it follows that the statement remains true for 
energy values $E+c$, with  $c>0$ small.}
\end{theorem}
Complementary to the above situation is the case when $I_E$ contains 
exactly one rest point. To have infinitely many closed characteristics, 
the nature of the rest point will come into play. {\it If the rest 
point is a center (four imaginary eigenvalues), then the system has 
infinitely many closed characteristics at each 
energy level sufficiently close to $E$, including $E$.}
If the rest point is not a center, there need not exist infinitely
many closed characteristics as results in \cite{vdB} indicate.

Similar results can be proved for compact interval components (for 
which dissipativity is irrelevant) and semi-infinite 
interval components $I_E\simeq \R^\pm$. 
\begin{theorem}\label{H4}
Suppose that a dissipative twist system has a singular energy level $E$ 
with an interval component $I_E=[a,b]$, or $I_{E}\simeq\R^\pm$, which
contains at least one rest point of saddle-focus/center type.
Then the system has infinitely many closed characteristics at  
energy level $E$. 
\end{theorem}
If an interval component contains no rest points, or only degenerate 
rest points (0 eigenvalues), then there need not exist infinitely many 
closed characteristics, completing our classification.

This classification immediately applies to the Swift-Hohenberg model 
\rmref{SH}, which is a twist system for all parameter values $\alpha \in \R$.
We leave it to the reader to apply the above theorems to the different regimes
of $\alpha$.

\subsection{Additional applications}\label{I3}

The framework of parabolic recurrence relations that we construct
is robust enough to accommodate several other important classes
of dynamics. 

\subsubsection{First-order nonautonomous Lagrangians}
 Finding periodic solutions of first-order Lagrangian systems of
the form $\delta \int L(x,u,u_x)dx=0$, with $L$ being 1-periodic in $x$, can
be rephrased in terms of parabolic recurrence relations of gradient type.
The homotopy index can be used to find periodic solutions $u(x)$ in 
this setting, even though a globally defined Poincar\'e map on $\R^2$ 
need not exist.

\subsubsection{Monotone twist maps}
A monotone twist map (compare \cite{Angenent2,moser}) is a 
(not necessarily area-preserving) map on $\real^2$ of the form 
\[
(u,p_u) \to (u',p_u'),\quad\quad {\frac{\partial u'}{\partial p_u}} > 0.
\]
Periodic orbits $\{ (u_i, p_{u_i}) \}$ are found
by solving a parabolic recurrence relation for the $u$-coordinates 
derived from the twist property.

\subsubsection{Uniformly parabolic PDE's}
The study of the invariant dynamics of Equation \rmref{PDE} can also be 
formulated in terms of parabolic recurrence relations by a spatial 
discretization. The basic theory for braid forcing developed here can 
be adapted to the dynamics of Equation \rmref{PDE}: see 
\cite{GV} for details.

\subsubsection{Lattice dynamics}
The form of a parabolic recurrence relation is precisely that 
arising from a set of coupled oscillators on a [periodic] one-dimensional
lattice with nearest-neighbor attractive coupling.
A similar setup arises in Aubry-LeDaeron-Mather theory of 
the Frenkel-Kontorova model \cite{Aubry-LeD}.
In this setting, a nontrivial homotopy index yields existence of 
invariant states (or stationary, in the exact context) within a 
particular braid class. Related physical systems (e.g., charge 
density waves on a 1-d lattice \cite{Middleton}) are also often 
reducible to parabolic recurrence relations.

\vskip.2cm
\subsection{History and outline}

The history of our approach is the convergence of ideas from knot theory, 
the dynamics of annulus twist maps, and curve shortening. 
We have already mentioned the similarities with Vassiliev's topological 
approach to discriminants in the space of immersed knots. From 
the dynamical systems perspective, the study of parabolic flows and 
gradient flows in relation with embedding data and the Conley index can be 
found in work of Angenent \cite{Angenent1,Angenent2,Angenent3} and Le Calvez 
\cite{LeCalvez,LeCalvez2} on area preserving twist maps. More general studies 
of dynamical properties of parabolic-type flows appear in numerous
works: we have been inspired by the work of Smillie \cite{Smillie},
Mallet-Paret and Smith \cite{MPSmith}, Hirsch \cite{Hir}, and,
most strongly, the work of Angenent on curve shortening \cite{Angenent3}.
Many of our applications to finding closed characteristics of second
order Lagrangian systems share similar goals with the programme of 
Hofer and his collaborators (see, \eg, \cite{EGH00,Hofer2,HofZehn}),
with the novelty that our energy surfaces are all non-compact
and not necessarily of contact type \cite{AVV}.

Clearly there is a parallel between the homotopy index theory presented here 
and Boyland's adaptation of Nielsen-Thurston theory for braid types 
of surface diffeomorphisms \cite{Boy94}. An important difference is 
that we require compactness only at the level of braid diagrams, 
which does not yield compactness on the level of the domains of the 
return maps [if these indeed exist]. Another important observation 
is that the recurrence relations are sometimes not defined on all 
of $\R^2$, which makes it very hard if not impossible to rephrase the 
problem of finding periodic solutions in terms of fixed points 
of 2-dimensional maps.

There are three components of this paper: (a) the precise 
definitions of the spaces involved and flows constructed, 
covered in \S\ref{braidspace}-\S\ref{parabolic}; (b) the 
establishment of existence, invariance, and properties of the 
index for braid diagrams in \S\ref{conley}-\S\ref{Morseth}; and 
(c) applications of the machinery to second order Lagrangian systems 
\S\ref{second}-\S\ref{compu}. Finally, \S\ref{postlude} contains open 
questions and remarks.
\vsp

\noindent {\bf Acknowledgments.} The authors would like express 
special gratitude to Sigurd Angenent and Konstantin Mischaikow 
for numerous enlightening discussions. Special thanks to 
Madjid Allili for his computational work in the earliest stages of this work.
Finally, the hard work of the referee has improved the 
paper in several respects, especially in the definitions of equivalent
relative braid classes.

\tableofcontents

%
%
\section{Spaces of discretized braid diagrams}
\label{braidspace}

\subsection{Definitions}\label{II1}

Recall the definition of a braid (see
\cite{Birman,Hansen} for a comprehensive introduction).
A braid $\beta$ on $n$ strands is a collection
of embeddings $\{\beta^\alpha:[0,1]\to\R^3\}_{\alpha=1}^n$ with disjoint
images such that (a) $\beta^\alpha(0)=(0,\alpha,0)$; (b)
$\beta^\alpha(1)=(1,\tau(\alpha),0)$ for some permutation $\tau$; and
(c) the image of each $\beta^\alpha$ is transverse to all planes 
$\{x={\mbox{const}}\}$.
We will ``read'' braids from left to right with respect to the 
$x$-coordinate. Two such braids are said to be of the same 
{\em topological braid class} if they are homotopic in the space of braids:
one can deform one braid to the other without any intersections 
among the strands. There is a natural group structure on the space of
topological braids with $n$ strands, $B_n$, given by concatenation. 
Using generators $\sigma_i$ which
interchange the $i^{th}$ and $(i+1)^{st}$ strands (with
a positive crossing) yields the presentation for $B_n$:
\begin{equation}
\label{eq_BraidGroup}
	B_n := \left< \sigma_1, \ldots, \sigma_{n-1} \, : \,
	\begin{array}{ccr}
	\sigma_i\sigma_j=\sigma_j\sigma_i & ; & \vert i-j\vert>1 \\
	\sigma_i\sigma_{i+1}\sigma_i=\sigma_{i+1}\sigma_i\sigma_{i+1}
	& ; & i<n-1
	\end{array}
	\right> .
\end{equation}

Braids find their greatest applications in knot theory via
taking their closures. Algebraically, the closed braids on
$n$ strands can be defined as the set of conjugacy 
classes\footnote{Note that we fix the number of strands and do not allow
the Markov move commonly used in knot theory.} 
in $B_n$. Geometrically, one quotients out the range
of the braid embeddings via the equivalence relation
$(0,y,z)\sim(1,y,z)$ and alters the restriction (a) and (b) of
the position of the endpoints to be
$\beta^\alpha (0)\sim\beta^{\tau(\alpha)}(1)$, 
as in Fig. \ref{fig_ex}[center]. Thus, a closed braid
is a collection of disjoint embedded loops in $S^1\times\R^2$
which are everywhere transverse to the $\R^2$-planes.

The specification of a topological braid class (closed or otherwise) may
be accomplished unambiguously by a labeled projection to the
$(x,y)$-plane: a {\it braid diagram}. 
Any braid may be perturbed
slightly so that pairs of strand crossings in the projection are
transversal: in this case, a marking of $(+)$ or $(-)$
serves to indicate whether the crossing is ``bottom over top''
or ``top over bottom'' respectively. Fig.~\ref{fig_ex}[center] 
illustrates a topological braid with all crossings positive.


\subsection{Discretized braids}
In the sequel we will restrict to a class of closed 
braid diagrams which have two special properties: (a) they
are {\it positive} --- that is, all crossings are of $(+)$ 
type; and (b) they are {\it discretized}, or piecewise linear
diagrams with constraints on the positions of anchor points. We 
parameterize such diagrams by the configuration space of anchor 
points.
\begin{definition}\label{PL}
The space of {\em discretized period $d$ braids on $n$ strands},
denoted $\Conf^n_d$, is the space of all pairs $(\uu,\tau)$ where
$\tau\in S_n$ is a permutation on $n$ elements, and $\uu$ is an 
unordered collection of $n$ {\em strands}, $\uu=\{\uu^\alpha\}_{\alpha=1}^n$,
satisfying the following conditions:
\begin{enumerate}
\item[(a)]
	Each strand consists of $d+1$ {\em anchor points}: 
	$\uu^\alpha=(u^\alpha_0,u^\alpha_1,\ldots,u^\alpha_d)\in\R^{d+1}$.
\item[(b)]
	For all $\alpha=1,\ldots,n$, one has 
\[	u^\alpha_d = u^{\tau(\alpha)}_0 .	\]
\item[(c)]
	The following {\em transversality condition} is satisfied:
	for any pair of distinct strands $\alpha$ and $\alpha'$
	such that $u^\alpha_i=u^{\alpha'}_i$ for some $i$, 
\begin{equation}
\label{eq_transverse}
        \bigl(u^\alpha_{i-1}-u^{\alpha'}_{i-1}\bigr)
	\bigl(u^\alpha_{i+1}-u^{\alpha'}_{i+1}\bigr) < 0 .
\end{equation}
\end{enumerate}
The topology on $\Conf^n_d$ is 
the standard topology of $\R^{n+1}$ on the strands and
the discrete topology with respect to the permutation $\tau$, 
modulo permutations which change orderings of strands. 
Specifically, two discretized braids $(\uu,\tau)$ and 
$(\tilde\uu,\tilde\tau)$ are close iff for some permutation 
$\sigma\in S_n$ one has $\uu^{\sigma(\alpha)}$ close to 
$\tilde\uu^\alpha$ (as points in $\R^{n+1}$) for all $\alpha$,
with $\sigma \circ\tilde\tau=\tau\circ\sigma$. 
\end{definition}
\begin{remark} 
\label{convention}
{\rm 
In Equation~(\ref{eq_transverse}), and indeed throughout the paper,
all expressions involving coordinates $u_i$ are considered {\em mod
the permutation $\tau$ at $d$}; thus, for every $j\in \Z$, we 
recursively define 
\begin{equation}
\label{eq_modperm}
	\begin{array}{rcl}
		u^\alpha_{d+j} &:=& u^{\tau(\alpha)}_{j} 
	\end{array} .
\end{equation}
As a point of notation, subscripts always refer to the spatial
discretization and superscripts always denote strands.
For simplicity, we will henceforth suppress the $\tau$ portion of a
discretized braid $\uu$. 
}
\end{remark}

One associates to each configuration $\uu\in\Conf^n_d$ the 
{\it braid diagram} $\braid(\uu)$, given as follows. For 
each strand $\uu^\alpha\in\uu$, consider the piecewise-linear (PL) 
interpolation
\begin{equation}\label{interpolate}
\braid^{\alpha}(s) := u^\alpha_{\floor{d\cdot s}}+(d\cdot s-\floor{d\cdot s})
	(u^\alpha_{\ceil{d\cdot s}}-u^\alpha_{\floor{d\cdot s}}),
\end{equation}
for $s\in[0,1]$.
The braid diagram $\braid(\uu)$ is then defined to be the 
superimposed graphs of all the functions $\braid^{\alpha}$, as
illustrated in Fig.~\ref{fig_ex}[right] for a period six braid on 
four strands (crossings are shown merely for suggestive purposes).
 






This explains the transversality condition of 
Equation~(\ref{eq_transverse}): a failure of this equation to hold 
implies that there is a PL-tangency in the associated
braid diagram. Since all crossings in a discretized braid diagram 
are PL-transverse, the map $\braid(\cdot)$ sends $\uu$ to a 
{\em topological} closed braid diagram once 
a convention for crossings is chosen. Inspired by lifting smooth
curves to a 1-jet extension, we label all crossings of $\braid(\uu)$ as 
positive type. This can be thought of as using the slope of the 
PL-extension of $\uu$ as the ``height'' of the braid strand
(though this analogy breaks down at the sharp corners). With 
this convention, then, the space $\Conf^n_d$ embeds into the space
of all closed positive braid diagrams on $n$ strands. 

\begin{definition}
Two discretized braids $\uu, \uu' \in\Conf^n_d$ are of the 
same {\em discretized braid class}, denoted $[\uu]=[\uu']$, 
if and only if they are in the same path-component of $\Conf^n_d$.
The {\em topological braid class}, $\{\uu\}$, denotes the path component
of $\beta(\uu)$ in the space of positive topological 
braid diagrams.
\end{definition}

The proof of the following lemma is essentially obvious.
\begin{lemma}
\label{braidtype}
If $[\uu]=[\uu']$ in $\Conf^n_d$, then the induced positive braid 
diagrams $\beta$ and $\beta'$ correspond to isotopic closed topological 
braid diagrams. 
\end{lemma}
The converse to this Lemma is not true: two discretizations of 
a topological braid are not necessarily connected in $\Conf^n_d$. 

Since one can write the generators $\sigma_i$ of the braid
group $B_n$ as elements of $\Conf^n_1$, it is clear that all
positive topological braids are representable as discretized braids.
Likewise, the relations for the groups of positive closed braids can
be accomplished by moving within the space of discretized braids; hence,
this setting suffices to capture all the relevant braid theory we 
will use.

\subsection{Singular braids}\label{II1-2}

The appropriate discriminant for completing the space $\Conf^n_d$ 
consists of those ``singular'' braid diagrams admitting tangencies
between strands. 
\begin{definition}
\label{closure}
Denote by $\bar{\Conf}_d^n$ the $nd$-dimensional vector space\footnote{
Strictly speaking $\bar{\Conf}_d^n$ is not a vector space, but
a union of vector spaces. Fixing appropriate
permutations its components are  vector spaces. Consider for instance
$\bar D_1^3$ which is a union of 3 copies of $\R^3$.}  
 of 
all discretized braid diagrams $\uu$ which satisfy properties 
(a) and (b) of Definition \ref{PL}.
Denote by $\Sigma_d^n:=\bar{\Conf}_d^n-\Conf_d^n$ the
set of singular discretized braids.
\end{definition}
We will often suppress the period and strand data and write 
$\Sigma$ for the space of singular discretized braids.
It follows from Definition~\ref{PL} and Equation~(\ref{eq_transverse})
that the set $\Sigma^n_d$ is a semi-algebraic variety in $\bar{\Conf^n_d}$. 
Specifically, for any singular braid $\uu\in \Sigma$ there exists
an integer $i \in \{1,\ldots,d\}$  and indices $\alpha\neq\alpha'$ 
such that $u^\alpha_i = u^{\alpha'}_i$, and 
\begin{equation}
        \bigl(u^\alpha_{i-1}-u^{\alpha'}_{i-1}\bigr)
	\bigl(u^\alpha_{i+1}-u^{\alpha'}_{i+1}\bigr) \geq 0 ,
\end{equation}
where the subscript is always computed mod the permutation $\tau$
at $d$. The number of such distinct occurrences is the codimension of the
singular braid diagram $\uu \in \Sigma$. We decompose $\Sigma$
into the union of strata $\Sigma[m]$ graded by $m$, the codimension of 
the singularity.

Any closed braid (discretized or topological) is partitioned into
components by the permutation $\tau$. 
Geometrically, the components are precisely the connected
components of the closed braid diagram. In our context, a component
of a discretized braid can be specified as $\{u_i^\alpha\}_{i\in\Z}$,
since, by our indexing convention, $i$ ``wraps around'' to the other
side of the braid when $i\not\in\{1,\ldots,d\}$. 

For singular braid diagrams of sufficiently high codimension, entire
components of the braid diagram can coalesce. This can happen in 
essentially two ways: (1) a single component involving multiple 
strands can collapse into a braid with fewer numbers of strands, 
or (2) distinct components can coalesce into a single component.
We define the {\it collapsed singularities}, $\Collapsed$, as follows:
$$
\Collapsed := \{ \uu \in \Sigma~|~ u_i^\alpha = u_i^{\alpha'},~\forall
i \in \Z,~~{\rm for~some}~~\alpha\neq\alpha'\} \subset \Sigma .
$$
Clearly the codimension of singularities in $\Collapsed$ is
at least $d$. 
Since for braid diagrams in $\Collapsed$ the number of strands
reduces, the subspace $\Collapsed$ may be decomposed into a union
of the spaces $\bar\Conf_d^{n'}$ for $n'< n$; i.e., 
$\Collapsed = \cup_{n'<n} \bar \Conf_d^{n'}$. If $n=1$, then 
$\Collapsed = \emptyset$.

\subsection{Relative braid classes}
Evolving certain components of a braid diagram
while fixing the remaining components  motivates working with a class
of ``relative'' braid diagrams. 

Given $\uu \in \bar\Conf_d^n$ and $\vv \in \bar\Conf_d^m$, 
the union $\uu \cup \vv \in \bar\Conf_d^{n+m}$ is 
naturally defined as the unordered union of the strands.  
Given $\vv\in\Conf_d^m$, define 
\[ \Conf_d^n{~\rel~}\vv := \{ \uu\in\Conf^n_d~:~ \uu \cup \vv \in
	\Conf_d^{n+m} \} ,\]
fixing $\vv$ and imposing transversality. 
The path components of $\Conf_d^{n}~\rel~\vv$ 
comprise the 
{\em relative discrete braid classes}, denoted $[\uu~\rel~\vv]$.
%
%
The braid $\vv$ will be called the {\it skeleton} henceforth.
The set of singular braids $\Sigma~{\rel}~\vv$ are those 
braids $\uu$ such that 
$\uu\cup\vv\in\Sigma_d^{n+m}$ 
The collapsed singular braids are denoted by $\Collapsed~{\rel}~\vv$. 
As before, the set $(\Conf_d^n~{\rel}~\vv)\cup(\Sigma~{\rel}~\vv)$
is the closure of $\Conf_d^n~{\rel}~\vv$
in $\R^{nd}$, and is denoted $\bar\Conf_d^n~{\rel}~\vv$.  
We denote by $\{\uu~\rel~\vv\}$ the topological relative braid class:
the set of topological (positive, closed) braids $\uu$ such that 
$\uu\cup\vv$ is a topological (positive, closed) braid diagram. 

Given two relative braid classes $[\uu~{\rel}~ \vv]$ and 
$[\uu'~{\rel}~ \vv']$ in $\Conf_d^n~\rel~\vv$ and $\Conf_d^n~\rel~\vv'$ 
respectively, to what extend are they the same? Consider the set
$$
{\bf D} = \{(\uu,\vv)\in \Conf_d^n\times \Conf_d^m~|~ 
\uu\cup\vv \in \Conf_d^{n+m}\}.
$$
The natural projection $\pi:~(\uu,\vv) \to \vv$ from
${\bf D}$ to $\Conf_d^m$ has as its fiber the braid class 
$[\uu~\rel~\vv]$. The path component of $(\uu,\vv)$ in 
${\bf D}$ will be denoted $\bigl[\uu~\rel~[\vv]\bigr]$. 
This generates the equivalence relation for relative braid 
classes to be used in the remainder of this work: 
$[\uu~\rel~\vv] \sim [\uu'~\rel~\vv']$ if and only if
$\bigl[\uu~\rel~[\vv]\bigr]=\bigl[\uu'~\rel~[\vv']\bigr]$.

Likewise, define $\bigl\{\uu~\rel~\{\vv\}\bigr\}$ 
to be the set of equivalent {\it topological} relative braid
classes. That is, $\{\uu~\rel~\vv\}\sim\{\uu'~\rel~\vv'\}$
if and only if there is a continuous family of 
topological (positive, closed) braid diagram pairs deforming 
$(\uu,\vv)$ to $(\uu',\vv')$.

%
%
%

%
%
\section{Parabolic recurrence relations}
\label{III}
\label{parabolic}
We consider the dynamics of vector fields given by recurrence 
relations on the spaces of discretized braid diagrams. 
These recurrence relations are nearest neighbor interactions --- each 
anchor point on a braid strand influences anchor points to the
immediate left and right on that strand ---  and resemble spatial 
discretizations of parabolic equations. 

\subsection{Axioms and exactness}\label{III1}
Denote by $\Seq$ the sequence space $\Seq:=\R^{\Z}$. 
\begin{definition}
\label{def_parabolic}
A {\em parabolic recurrence relation $\RR$ on $\Seq$} is a sequence of 
real-valued $C^1$ functions $\RR=(\RR_i)_{i\in\Z}$ satisfying
\begin{description}

\item[(A1)]   [{\it monotonicity}]\footnote{Equivalently, one could impose 
$\partial_1 \RR_i\ge 0$ and $\partial_3 \RR_i> 0$ for all $i$.}
$\partial_1 \RR_i>0$ and $\partial_3 \RR_i\ge 0$
for all $i \in \Z$ 

\item[(A2)]  [{\it periodicity}]  For some $d\in \N$, 
$\RR_{i+d} = \RR_i$ for all $i\in \Z$.

\end{description}
\end{definition}

For applications to Lagrangian dynamics a variational structure 
is necessary. At the level of recurrence relations this implies that
$\RR$ is a gradient:
\begin{definition}
A parabolic recurrence relation on $\Seq$ is called {\em exact} if
\begin{description}
\item[(A3)]  [{\it exactness}]
There exists a sequence of $C^2$ {\it generating functions}, 
$(S_i)_{i\in\Z}$, satisfying 
\begin{equation}
	\RR_i(u_{i-1},u_i,u_{i+1}) 
	= \partial_2 S_{i-1}(u_{i-1},u_i) + \partial_1 S_i(u_i,u_{i+1}) ,
\end{equation}
for all $i \in \Z$.
\end{description}
\end{definition}
In discretized Lagrangian problems the action functional naturally 
defines the generating functions $S_i$. 
This agrees with the ``formal'' action in this case: $W(\uu) := \sum_i
S_i(u_i,u_{i+1})$. In this general setting, $\RR = \nabla W$.

\subsection{The induced flow}\label{III2}
In order to define parabolic flows we regard $\RR$ as a vector
field on $\Seq$: consider the differential equations
\begin{eqnarray}
\label{gradflow}
	\frac{d}{dt}u_i = 
	\RR_i(u_{i-1},u_i,u_{i+1}), \quad \uu(t) \in \Seq,~t\in \R .
\end{eqnarray}
Equation \rmref{gradflow} defines a (local) $C^1$ flow 
$\psi^t$ on $\Seq$ under any periodic 
boundary conditions with period $nd$.  
%
To define flows on the finite dimensional spaces $\bar\Conf_d^n$,
one considers the same equations:
\begin{equation}
\label{braidgradflow}
\frac{d}{dt}u^\alpha_i = \RR_i(u^\alpha_{i-1},u^\alpha_i,u^\alpha_{i+1}), 
\quad \uu \in \bar\Conf^n_d.
\end{equation}
where the ends of the braid are identified as per 
Remark~\ref{convention}. Axiom (A2) guarantees that the flow
is well-defined. Indeed, one may consider a  
cover of $\bar\Conf^n_d$ by taking the bi-infinite periodic 
extension of the braids: this yields a subspace of periodic 
sequences in $\Seq^n:=\Seq\times\cdots\times\Seq$ 
invariant under the product flow of (\ref{gradflow}) thanks to Axiom (A2). 
Any flow $\Psi^t$ generated by \rmref{braidgradflow} for some
parabolic recurrence relation $\RR$
is called a {\it parabolic flow on discretized braids}.
In the case of relative classes $\bar\Conf_d^n~\rel~\vv$ a parabolic
flow is the restriction of a parabolic flow on $\bar\Conf_d^{n+m}$ which
fixes the anchor points of the skeleton $\vv$.
We abuse notation and indicate the invariance of the skeleton  
by $\Psi^t(\vv)=\vv$. Indeed, 
for appropriate coverings of the skeletal strands $\vv^\alpha$ it holds 
that $\psi^t(\vv^\alpha)
=\vv^\alpha$.

\subsection{Monotonicity and braid diagrams}\label{III3}

The monotonicity Axiom (A1) in the previous subsection has a 
very clean interpretation in the space of braid diagrams.
Recall from \S\ref{braidspace} 
that any discretized braid $\uu$ has an associated diagram 
$\braid(\uu)$ which can be interpreted as a positive closed
braid. Any such diagram in general position can be expressed in terms of 
the (positive) generators $\{\sigma_j\}_{j=1}^{n-1}$ of the braid group
$B_n$. While this word is not necessarily unique, the length 
of the word is, as one can easily see from the presentation
of $B_n$ and the definition of $\Conf^n_d$. The length of 
a closed braid in the generators $\sigma_j$ is thus precisely
the {\em word metric} $\word{\cdot}$ from geometric group theory. The 
geometric interpretation of $\word{\uu}$ for a braid
$\uu$ is clearly the number of pairwise strand 
crossings in the diagram $\braid(\uu)$. 

\begin{figure}[hbt]
\begin{center}
\psfragscanon
\psfrag{u}[][]{\large $\uu$}
\psfrag{S}[][]{\large $\Sigma[1]$}
\psfrag{0}[][]{\large $i-1$}
\psfrag{1}[][]{\large $i$}
\psfrag{B}[][]{}
\psfrag{2}[][]{\large $i+1$}
\includegraphics[angle=0,width=5in]{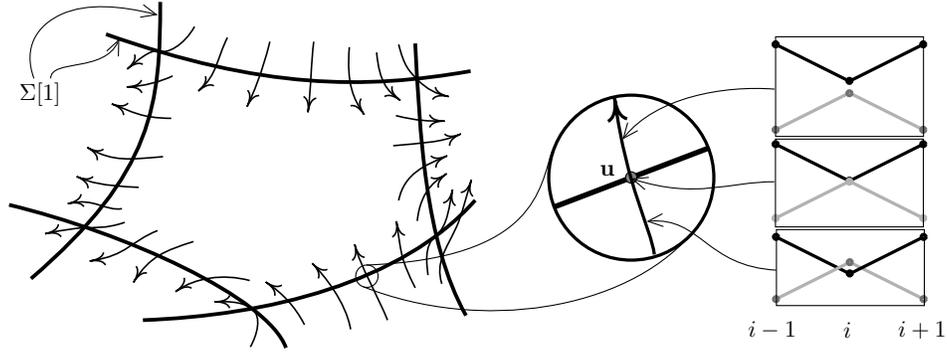}
\caption{A parabolic flow on a discretized braid class is 
transverse to the boundary faces. The local linking of strands
decreases strictly along the flowlines at a singular braid
$\uu$.}
\label{fig_transverse} 
\end{center}
\end{figure}

The primary result of this section is that the word metric acts as 
a discrete Lyapunov function for any parabolic flow on $\bar\Conf_d^n$.
This is really the braid-theoretic version of the lap number
arguments that have been used in several related settings 
\cite{Angenent1,Angenent3, Angenent4,FR99,FusOl,MPSmith,Mat1,Smillie}. 
The result we prove below can be
excavated from these cited works; however, we choose to give a brief
self-contained proof for completeness.
\begin{proposition}
\label{word}
Let $\Psi^t$ be a parabolic flow on $\bar\Conf_d^n$.
\begin{enumerate}
\item[(a)]
         For each point $\uu\in\Sigma-\Collapsed$, the local orbit 
	$\{\Psi^t(\uu) : t\in[-\epsilon,\epsilon]\}$ intersects $\Sigma$ 
	uniquely at $\uu$ for all $\epsilon$ sufficiently small.
\item[(b)]
        For any such $\uu$, the length of the braid diagram
	$\Psi^t(\uu)$ for $t>0$ in the word metric is strictly less 
	than that of the diagram $\Psi^t(\uu)$, $t<0$. 
\end{enumerate}
\end{proposition}
{\it Proof.}
Choose a point $\uu$ in $\Sigma$ representing a singular braid
diagram. We induct on the codimension $m$ of the singularity. In the 
case where $\uu\in\Sigma[1]$ (\ie, $m=1$), 
there exists a unique $i$ and a unique pair
of strands $\alpha\neq\alpha'$ such that $u^\alpha_i=u^{\alpha'}_i$ and
\[
	(u^\alpha_{i-1}-u^{\alpha'}_{i-1})(u^\alpha_{i+1}-u^{\alpha'}_{i+1})>0
.\]
Note that the inequality is strict since $m=1$.
We deduce from \rmref{braidgradflow} that	
\begin{eqnarray*}
	\left.\frac{d}{dt}(u^\alpha_i-u^{\alpha'}_i)\right\vert_{t=0}
	 &=& \RR_i(u^\alpha_{i-1},u^\alpha_i,u^\alpha_{i+1}) -
	     \RR_i(u^{\alpha'}_{i-1},u^{\alpha'}_i,u^{\alpha'}_{i+1}).
\end{eqnarray*}
 From Axiom (A2) one has that 
\[
 	\sign\left(\RR_i(u^\alpha_{i-1},u^\alpha_i,u^\alpha_{i+1}) 
	- \RR_i(u^{\alpha'}_{i-1},u^{\alpha'}_i,u^{\alpha'}_{i+1})\right)
	= \sign(u^\alpha_{i-1}-u^{\alpha'}_{i-1}).
\]
Therefore, as $t\to 0-$, the two strands have two local crossings, 
and as $t\to 0+$, these two strands are locally unlinked 
(see Fig.~\ref{fig_transverse}): the length of the 
braid word in the word metric is thus decreased by two, and
the flow is transverse to $\Sigma[1]$. This proves (a) and (b) 
on $\Sigma[1]$. 

Assume inductively that (a) and (b) are true for every point
in $\Sigma[m]$ for $m<M$. To prove (a) on $\Sigma[M]$, choose 
$\uu\in\Sigma[M]$. There are exactly $M$ distinct pairs of anchor 
points of the braid which coalesce at the braid diagram $\uu$. Since
the vector field $\RR$ is defined by nearest neighbors, 
singularities which are not strandwise consecutive in the braid 
behave independently to first order under the parabolic flow. Thus, it
suffices to assume that for some $i$, $\alpha$, and $\alpha'$ one has 
$\{u^\alpha_{i+j}\}_{j=0}^{M+1}$ and 
$\{u^{\alpha'}_{i+j}\}_{j=0}^{M+1}$ chains of 
consecutive anchor points for the braid diagram $\uu$ such 
that $u^\alpha_{i+j}=u^{\alpha'}_{i+j}$ if and only if $1\leq j\leq M$. 
(Recall that the addition $i+j$ is always done modulo the 
permutation $\tau$ at $d$). Then since
\begin{eqnarray*}
	\left.\frac{d}{dt}(u^\alpha_{i+j}-u^{\alpha'}_{i+j})\right\vert_{t=0}
	 &=& \RR_{i+j}(u^\alpha_{i+j-1},u^\alpha_{i+j},u^\alpha_{i+j+1}) -
	 \RR_{i+j}(u^{\alpha'}_{i+j-1},u^{\alpha'}_{i+j},u^{\alpha'}_{i+j+1}),
\end{eqnarray*}
it follows that for all $j=2,..,(M-1)$, the anchor points $u^\alpha_{i+j}$
and $u^{\alpha'}_{i+j}$ are not separated to first order. 
At the left ``end'' of the singular braid, where $j=0$, 
\[	\RR_i(u^\alpha_{i-1},u^\alpha_i,u^\alpha_{i+1})
	-\RR_i(u^{\alpha'}_{i-1},u^{\alpha'}_i,u^{\alpha'}_{i+1})\neq 0 , \]
so that the vector field $\RR$ is tangent to $\Sigma$ at $\uu$ but
is not tangent to $\Sigma[M]$: the flowline through $\uu$ 
decreases codimension immediately. By the induction hypothesis on (b), the 
flowline through $\uu$ cannot possess intersections with $\Sigma[m]$ for 
$m<M$ which accumulate onto $\uu$ --- the length of the braids are
finite. Thus the flowline intersects $\Sigma$ locally at $\uu$ uniquely. 
This concludes the proof of (a). 

It remains to show that the length of the braid word decreases strictly 
at $\uu$ in $\Sigma[M]$. By (a), the flow $\Psi^t$ is nonsingular in a 
neighborhood of $\uu$; thus, by the Flowbox Theorem, there is 
a tubular neighborhood of local $\Psi^t$-flowlines about $\Psi^t(\uu)$. 
The beginning and ending points of these local flowlines 
all represent nonsingular diagrams with the same word lengths 
as the beginning and endpoints of the path through $\uu$, since the complement
of $\Sigma$ is an open set. Since $\Sigma$ is a codimension-$1$ 
algebraic semi-variety in $\bar\Conf_d^n$, it follows from transversality 
that most of the nearby orbits intersect $\Sigma[1]$, at which braid 
word length strictly decreases. This concludes the proof of (b). 
\fp
\vsp

To put this result in context with the literature, we note that 
the monotonicity in \cite{FusOl,MPSmith} is one-sided: translated 
into our terminology, $\partial_3\RR_i=0$ for all $i$.
One can adapt this proof to generalizations of parabolic recurrence
relations appearing in the work of Le Calvez \cite{LeCalvez,LeCalvez2}:
namely, compositions of twist symplectomorphisms of the annulus reversing
the twist-orientation.

As pointed out above a parabolic flow on $\bar\Conf_d^n~\rel~\vv$ is a 
special case of a parabolic flow on $\bar\Conf_d^{n+m}$
with a fixed skeleton $\vv \in \Conf_d^{m}$, and therefor the analogue of the above
proposition for relative classes thus follows as a special case.

\begin{figure}[h]
\begin{center}
\includegraphics[angle=0,width=4.5in]{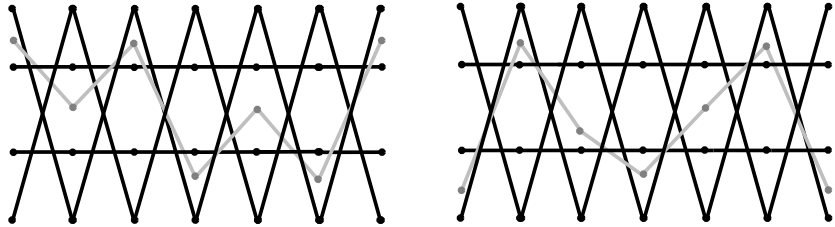}
\caption{Two relative braids with the same linking data 
	but different homotopy indices. The free strands are in grey.}
\label{fig_different} 
\end{center}
\end{figure}

\begin{remark}
{\rm
The information that we derive from relative braid diagrams
is more than what one can obtain from lap numbers alone
(cf. \cite{LeCalvez}). Fig.~\ref{fig_different} gives examples
of two closed discretized relative braids which have the 
same set of pairwise intersection numbers of strands (or lap numbers) 
but which force very different dynamical behaviors. The 
homotopy invariant we define in the next section distinguishes
these braids. The index of the first picture can be computed 
to be trivial, and the index for the second picture is computed 
in \S\ref{XI} to be nontrivial.
}
\end{remark}

%
%
\section{The homotopy index for discretized braids}
\label{IV}
\label{conley}

Technical lemmas concerning existence of certain 
types of parabolic flows are required for showing the existence
and well-definedness of the Conley index on braid classes. We
relegate these results to Appendix~\ref{app_1}.

\subsection{Review of the Conley index}\label{primer}

We include a brief primer of the relevant ideas from Conley's
index theory for flows. For a more comprehensive treatment, 
we refer the interested reader to \cite{Mischaik}. 

In brief, the Conley index is an extension of the Morse index.
Consider the case of a nondegenerate gradient flow: the Morse 
index of a fixed point is then the dimension of the unstable 
manifold to the fixed point. In contrast, the Conley index
is the homotopy type of a certain pointed space (in this 
case, the sphere of dimension equal to the Morse index). 
The Conley index can be defined for sufficiently ``isolated'' 
invariant sets in any flow, not merely for fixed points of gradients. 

Recall the notion of an isolating neighborhood as introduced by Conley 
\cite{Conley}.
Let $X$ be a locally compact metric space. A compact set 
$N \subset X$ is an {\em isolating neighborhood} 
for a flow $\psi^t$ on $X$ if the maximal invariant set $\Inv(N)
:= \{x \in N~|~ {\rm cl}\{\psi^t(x)\}_{t \in \R} \subset N\}$
is contained in the interior of $N$. The invariant set $\Inv(N)$
is then called a {\it compact isolated invariant set} for $\psi^t$.
In \cite{Conley} it is shown that every compact isolated invariant set 
$\Inv(N)$ admits a pair $(N,N^-)$ such that (following the definitions 
given in \cite{Mischaik}) (i) $\Inv(N) = \Inv({\rm cl}(N-N^-))$ with
$N-N^-$ a neighborhood of $\Inv(N)$; (ii) $N^-$ is positively 
invariant in $N$; and (iii) $N^-$ is an {\em exit set} for $N$: 
given $x\in N$ and $t_1>0$ such that 
$\psi^{t_1}(x) \not \in N$, then there exists a $t_0 \in [0,t_1]$ for which
$\{\psi^t(x)~:~t\in [0,t_0]\} \subset N$ and $\psi^{t_0}(x) \in N^-$.
Such a pair is called an {\it index pair} for $\Inv(N)$.
The Conley index, $h(N)$, is then defined as the homotopy type of 
the pointed space $\bigl(N/N^-,[N^-]\bigr)$, abbreviated 
$\bigl[N/N^-\bigr]$. This homotopy class is independent of the 
defining index pair, making the Conley index well-defined.  

A large body of results and applications of the Conley index theory
exists. We recall following \cite{Mischaik} two foundational results. 
\begin{enumerate}
\item[(a)]
       {\bf Stability of isolating neighborhoods:}
Any isolating neighborhood $N$ for a flow $\psi^t$ is an isolating
neighborhood for all flows sufficiently $C^0$-close to $\psi^t$. 
\item[(b)]
        {\bf Continuation of the Conley index:}
Let $\psi^t_\lambda$, $\lambda\in[0,1]$ be a continuous family of 
flows with $N_\lambda$ a family of isolating neighborhoods.
Define the parameterized flow $(t,x,\lambda) \mapsto 
(\psi^t_\lambda (x),\lambda)$ on $X\times [0,1]$,
and $N=\cup_\lambda \bigl( N_\lambda \times \{\lambda\}\bigr)$. 
If $N \subset X\times [0,1]$ is an isolating neighborhood for the
parameterized flow 
then the index $h_\lambda = h(N_\lambda,\psi^t_\lambda)$ 
is invariant under $\lambda$. 
\end{enumerate}
Since the homotopy type of a space is notoriously difficult to 
compute, one often passes to homology or cohomology. One defines 
the Conley homology\footnote{In \cite{ConleyZehn2} \v{C}ech cohomology 
is used. For our purposes ordinary singular (co)homology always suffices.}
of $\Inv(N)$ to be $CH_*(N) := H_*(N,N^-)$, 
where $H_*$ is singular homology. To the
homological Conley index of an index pair $(N,N^-)$
one can also assign the characteristic polynomial
$CP_t(N):= \sum_{k\ge 0} \beta_k t^k$, where
$\beta_k$ is the free rank of $CH_k(N)$. 
Note that, in analogy with Morse homology, if $CH_*(N)\neq 0$, 
then there exists a nontrivial invariant set within the 
interior of $N$. For more detailed description see \S\ref{morseth}.

\subsection{Proper and bounded braid classes}

 From Proposition~\ref{word}, one readily sees that 
complements of $\Sigma$ yield isolating neighborhoods, except for
the presence of the collapsed singular braids $\Collapsed$, which is 
an invariant set in $\Sigma$.  
For the remainder of this paper we restrict our attention to 
those relative braid diagrams whose braid classes prohibit 
collapse. 

\begin{figure}[hbt]
\psset{xunit=1cm,yunit=0.5cm}
\begin{pspicture}(-0.25,-0.75)(10.25,6.75)
 

  \psset{linewidth=1pt}
  \psline{c-c}(0,-0.5)(0,6.5)
  \psline{c-c}(4,-0.5)(4,6.5)
  \psline{c-c}(6,-0.5)(6,6.5)
  \psline{c-c}(10,-0.5)(10,6.5)

  \psset{linewidth=1.5pt}

\psline[linecolor=lightgray](0,1)(1,5)(2,1.3)(3,3.5)(4,1)
\psline[linecolor=lightgray](6,1)(7,3)(8,5)(9,3)(10,1)

\psline(0,0)(4,0)
\psline(0,6)(4,6)
\psline(6,0)(10,0)
\psline(6,6)(10,6)
\psline(0,2)(1,4)(2,2)(3,4)(4,2)
\psline(0,4)(1,2)(2,4)(3,2)(4,4)
\psline(6,2)(7,4)(8,2)(9,4)(10,2)
\psline(6,4)(7,2)(8,4)(9,2)(10,4)
 
  \psset{dotsize={4pt 0},linecolor=black}
\multips(0,0)(0,2){4}{\psdots(0,0)(1,0)(2,0)(3,0)(4,0)}
\multips(6,0)(0,2){4}{\psdots(0,0)(1,0)(2,0)(3,0)(4,0)}

  \psset{linecolor=lightgray}
\psdots(0,1)(1,5)(2,1.3)(3,3.5)(4,1)
\psdots(6,1)(7,3)(8,5)(9,3)(10,1)

\end{pspicture}
\caption{Improper [left] and proper [right] relative braid classes.
	Both are bounded.}
\label{fig_proper} 
\end{figure}
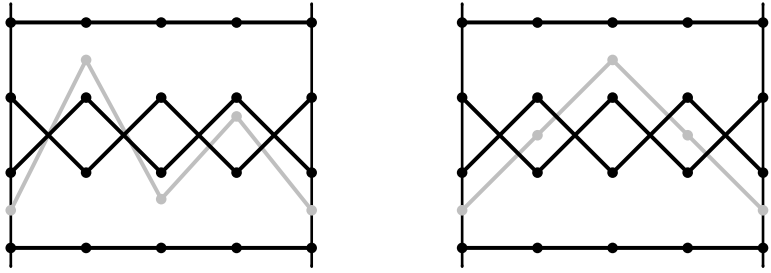

Fix $\vv \in \DD_d^m$, and consider the relative 
braid classes $\{\uu~\rel~\vv\}$ (topological) and $[\uu~\rel~\vv]$
(discretized). 
\begin{definition}\label{proper}
A topological relative braid class $\{\uu~\rel~\vv\}$ is
{\em proper} if it is impossible to find a continuous path of braid diagrams 
$\uu(t)~\rel~\vv$ for $t\in[0,1]$ such that 
$\uu(0)=\uu$, $\uu(t)~\rel~\vv$ defines a braid for all $t\in[0,1)$, 
and $\uu(1)~\rel~\vv$ is a diagram where an entire component of the 
closed braid has collapsed onto itself or onto another component
of $\uu$ or $\vv$.    
A discretized relative braid class $[\uu~\rel~\vv]$ is called {\em proper} if
the associated topological braid class is proper, otherwise, it is 
{\em improper}: see Fig.~\ref{fig_proper}. 
\end{definition}
\begin{definition}\label{bounded}
A topological relative braid class $\{\uu~\rel~\vv\}$ is called 
{\em bounded} if there exists a uniform bound on all 
representatives ${\uu}$ of the equivalence class, i.e. on 
the strands $\beta(\uu)$ (in $C^0([0,1])$).
A discrete relative braid class $[\uu~\rel~\vv]$ is called bounded if 
the set $[\uu~\rel~\vv]$ is bounded.
\end{definition}
Note that if a 
topological class $\{\uu~\rel~\vv\}$ is bounded then the discrete class
$[\uu~\rel~\vv]$ is bounded as well for any period. The converse does
not always hold. Bounded braid classes possess a compactness 
sufficient to implement the Conley index theory without further assumptions. 
It is not hard either to see or to prove that properness and
boundedness are well-defined properties of equivalence classes of braids. 


\subsection{Existence and invariance of the Conley index for braids}\label{IV2}

\begin{theorem}
\label{isolating}
Suppose $[\uu~\rel~\vv]$ is a bounded proper relative braid class and 
$\Psi^t$ is a parabolic flow fixing $\vv$. Then the following 
are true: 
\begin{enumerate}
\item[(a)]
        $N:={\rm cl}[\uu~\rel~\vv] $ is an isolating neighborhood 
	for the flow $\Psi^t$, which thus yields a well-defined 
	Conley index $h(\uu~\rel~\vv):= h(N)$; 
\item[(b)]
        The index $h(\uu~\rel~\vv)$ is independent of
	the choice of parabolic flow $\Psi^t$ so long as 
	$\Psi^t(\vv)=\vv$; 
\item[(c)]
        The index $h(\uu~\rel~\vv)$ is an invariant of 
	$\bigl[\uu~\rel~[\vv]\bigr]$.
\end{enumerate}
\end{theorem}

\begin{definition}
The {\em homotopy index} of a bounded proper discretized braid class
$\bigl[\uu~\rel~[\vv]\bigr]$ in $\DD_d^n~\rel~[\vv]$ is defined to be
$h(\uu~\rel~\vv)$, the Conley index of the braid class 
$[\uu~\rel~\vv]$ with respect to some (hence any) parabolic flow 
fixing any representative $\vv$ of the skeletal braid class $\pi\bigl[\uu~\rel~[\vv]\bigr]
\subset [\vv]$. 
\end{definition}

{\it Proof.} 
Isolation is proved by examining $\Psi^t$ on the boundary $\partial N$.
By Definition \ref{proper} and \ref{bounded} the set $N$ is compact, and 
$\partial N \subset \Sigma\backslash \Collapsed$.
Choose a point $\uu$ on $\partial N$. 
Proposition \ref{word} implies that the parabolic 
flow $\Psi^t$ locally intersects $\partial N$ at $\uu$ alone and 
that furthermore its length in the braid group strictly decreases. 
This implies that under $\Psi^t$, the point $\uu$ exits the set $N$ 
either in forwards or backwards time (if not both). Thus, 
$\uu\not\in\Inv(N)$ and  (a) is proved. 

Denote by $h(\uu~\rel~\vv)$ the index of $\Inv(N)$. 
To demonstrate  (b), consider two parabolic flows $\Psi^t_0$ and 
$\Psi^t_1$ that satisfy all our requirements, and consider the 
isolating neighborhood $N$ valid for both flows. 
Construct a homotopy $\Psi^t_\lambda$, $\lambda \in [0,1]$, by
considering the parabolic recurrence functions 
$\RR^\lambda = (1-\lambda) \RR^0 + \lambda \RR^1$,
where $\RR^0$ and $\RR^1$ give rise to the flows  
$\Psi^t_0$ and $\Psi^t_1$ respectively.
It follows immediately that $ \Psi^t_\lambda(\vv)=\vv$, for 
all $\lambda \in [0,1]$; therefore $N$ is an isolating neighborhood 
for $\Psi^t_\lambda$ with $\lambda \in [0,1]$.
Define $\Inv_\lambda(N)$, $\lambda \in [0,1]$, to be the maximal 
invariant set in $N$ with respect to the flow $\Psi^t_\lambda$.
The continuation property of the Conley index completes the 
proof of (b).


Assume that $[\uu~\rel~\vv]\sim[\uu'~\rel~\vv']$, so that there
is a continuous path $(\uu(\lambda),\vv(\lambda))$, for 
$0\leq\lambda\leq 1$, of 
braid pairs within $\Conf^{n+m}_d$ between the two.
From the proof of Lemma \ref{homotopy} in Appendix A, there exists a 
continuous family of flows $\Psi^t_\lambda$, such that 
$\Psi^t_\lambda( \vv(\lambda)) =  \vv(\lambda)$, for all
$\lambda \in [0,1]$.   
Item (a) ensures that $N_\lambda:={\rm cl}[\uu~\rel~\vv(\lambda)]$
is an isolating neighborhood for all $\lambda \in [0,1]$.
The continuity of $ \vv(\lambda)$ implies that 
the set $N:=\cup_\lambda \bigl(N_\lambda \times \{\lambda\}\bigr)
\subset \bar \Conf_d^n \times [0,1]$ is an
isolating neighborhood for the parameterized flow
$(\Psi^t_\lambda(\uu),\lambda)$ on $\bar \Conf_d^n \times [0,1]$.
%
Therefore via the continuation property of the Conley index, 
$h(\uu~\rel~\vv(\lambda))$ is independent of 
$\lambda \in [0,1]$, which completes the proof of Item (c).
\fp

\subsection{An intrinsic definition}\label{IV4}
 For any bounded proper relative braid class $[\uu ~{\rel}~\vv]$
we can define its index intrinsically, independent of any notions 
of parabolic flows. Denote as before by $N$ the set 
${\rm cl}[\uu ~{\rel}~\vv]$ within $\bar\Conf_d^n$. 
The singular braid diagrams $\Sigma$ partition $\bar\Conf_d^n$
into disjoint cells (the discretized relative braid classes), the closures
of which contain portions of $\Sigma$. For a bounded proper braid 
class, $N$ is compact, and $\partial N$ avoids $\Collapsed$. 

To define the exit set $N^-$, consider any point $\ww$ on 
$\partial N\subset\Sigma$. There exists a small neighborhood
$W$ of $\ww$ in $\bar\Conf_d^n$ for which the subset $W-\Sigma$ 
consists of a finite number of connected components $\{W_j\}$.
Assume that $W_0 = W\cap N$. 
We define $N^-$ to be the set of $\ww$ for which the word
metric is locally maximal on $W_0$, namely, 
\begin{equation}\label{intrinsic}
N^- := {\rm cl}\left\{
	\ww\in\partial N : \word{W_0}\geq\word{W_j} \,\forall 
		j>0 \right\}.
\end{equation}

We deduce that $(N,N^-)$ is an index
pair for any parabolic flow for which $\Psi^t(\vv) = \vv$, and
thus by the independence of $\Psi^t$, the homotopy type
$\bigl[N/N^-\bigr]$ gives the Conley index. The index can be 
computed by choosing a representative $\vv \in [\vv]$ and 
determining $N$ and $N^-$. A rigorous computer assisted approach 
exists for computing the homological index using cube complexes 
and digital homology \cite{AGV}.


\subsection{Three simple examples}\label{examples}

It is not obvious what the homotopy index is measuring 
topologically. Since the space $N$ has one dimension per free 
anchor point, examples quickly become complex. 

\begin{figure}[hbt]
\begin{center}
\psfragscanon
\psfrag{0}[][]{\Large $u_0$}
\psfrag{1}[][]{\Large $u_1$}
\psfrag{2}[][]{\Large $u_2$}
\psfrag{a}[][]{\Large $u_2$}
\psfrag{b}[][]{\Large $u_1$}
\psfrag{S}[][]{\Large $\Sigma$}
\psfrag{U}[][]{\Large $\subset$}
\includegraphics[angle=0,width=5in]{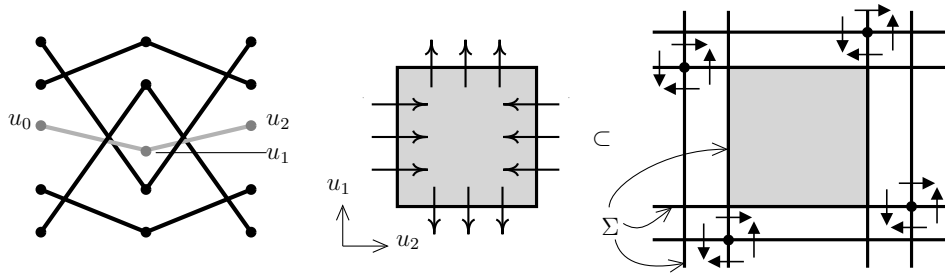}
\caption{The braid of Example 1 [left] and the associated configuration
space with parabolic flow [middle]. On the right is an expanded view
of $\Conf^1_2~\rel~\vv$ where the fixed points of the flow correspond
to the four fixed strands in the skeleton $\vv$. The 
braid classes adjacent to these fixed points are not proper.}
\label{example1} 
\end{center}
\end{figure}

\vspace{0.1in}
\noindent
{\bf Example 1:}
Consider the proper period-2 braid illustrated in 
Fig.~\ref{example1}[left]. (Note that deleting any strand in 
the skeleton yields an improper braid.) 
There is exactly one free strand with two 
anchor points (recall that these are {\it closed} braids and the left 
and right sides are identified). The anchor point in the middle, 
$u_1$, is free to move vertically between the fixed points on the 
skeleton. At the endpoints, one has a singular braid in $\Sigma$ 
which is on the exit set since a slight perturbation sends this 
singular braid to a different braid class with fewer crossings.
The end anchor point, $u_2$ ($=u_0$) can freely move vertically in 
between the two fixed points on the skeleton. The singular boundaries
are in this case {\it not} on the exit set since pushing $u_2$
across the skeleton increases the number of crossings. 

Since the points $u_1$ and $u_2$ can be moved independently,
the configuration space $N$ in this case is the product of
two compact intervals. The exit set $N^-$ consists of those 
points on $\partial N$ for which $u_1$ is a boundary point.
Thus, the homotopy index of this relative braid is 
$[N/N^-]\simeq S^1$. 

\vspace{0.1in}

\begin{figure}[hbt]
\begin{center}
\psfragscanon
\psfrag{0}[][]{$u_0$}
\psfrag{1}[][]{$u_1$}
\psfrag{2}[][]{$u_2$}
\psfrag{3}[][]{$u_3$}
\includegraphics[angle=0,width=5in]{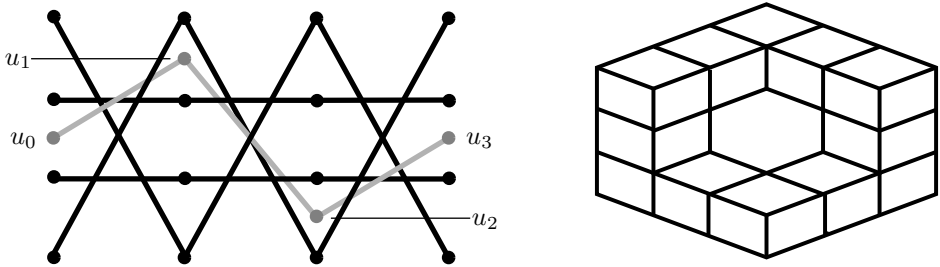}
\caption{The braid of Example 2 and the configuration
space $N$.}
\label{example2} 
\end{center}
\end{figure}

\noindent
{\bf Example 2:}
Consider the proper relative braid presented in 
Fig.~\ref{example2}[left]. Since there is one free strand of
period three, the configuration space $N$ is determined by 
the vector of positions $(u_0,u_1,u_2)$ of the anchor points. 
This example differs greatly from the previous example. For 
instance, the point $u_0$ (as represented in the figure) may
pass through the nearest strand of the skeleton above and below
without changing the braid class. The points $u_1$ and $u_2$ 
may not pass through any strands of the skeleton without 
changing the braid class {\it unless} $u_0$ has already passed
through. In this case, either $u_1$ or $u_2$ (depending on 
whether the upper or lower strand is crossed) becomes free.

To simplify the analysis, consider $(u_0,u_1,u_2)$ as all 
of $\R^3$ (allowing for the moment singular braids and 
other braid classes as well). The position of the skeleton
induces a cubical partition of $\R^3$ by planes, the 
equations being $u_i=v^\alpha_i$ for the various strands
$v^\alpha$ of the skeleton $\vv$. The braid class $N$ is
thus some collection of cubes in $\R^3$. In Fig.~\ref{example2}[right], 
we illustrate this cube complex associated to $N$, claiming 
that it is homeomorphic to $D^2\times S^1$. In this 
case, the exit set $N^-$ happens to be the entire boundary 
$\partial N$
and the quotient space is homotopic to the wedge-sum $S^2\vee S^3$.

\begin{figure}[hbt]
\begin{center}
\includegraphics[angle=0,width=5in]{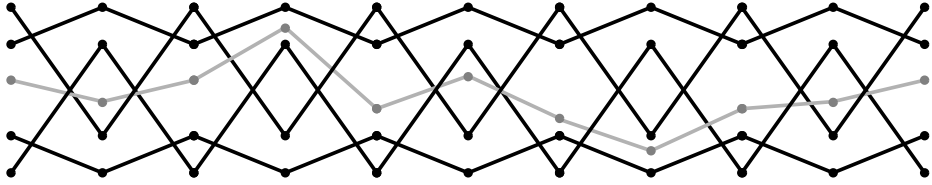}
\caption{The lifted skeleton of Example 1 with one free strand.}
\label{example3} 
\end{center}
\end{figure}

\vspace{0.1in}
\noindent
{\bf Example 3:}
To introduce the spirit behind the forcing theorems of the latter 
half of the paper, we reconsider the period two braid of Example 1. 
Take an $n$-fold cover of the skeleton as illustrated in 
Fig.~\ref{example3}. By weaving a single free strand in and
out of the strands as shown, it is possible to generate numerous
examples with nontrivial index. A moment's meditation suffices to 
show that the configuration space $N$ for this lifted braid 
is a product of $2n$ intervals, the exit set being completely determined 
by the number of times the free strand is ``threaded'' through the
inner loops of the skeletal braid as shown. 

For an $n$-fold cover with one free strand we can select a family of
$3^n$ possible braid classes describes as follows:
the even anchor points of the free strand are always in the middle, 
while for the odd anchor points there are three possible choices.
Two of these braid classes are
not proper. All of the remaining $3^n-2$ braid classes are bounded 
and have homotopy indices equal to a sphere $S^k$ for some $0\leq k\leq n$.  
Several of these strands may be superimposed while maintaining a
nontrivial homotopy index for the net braid: we leave it to the 
reader to consider this interesting situation. 

Stronger results follow from projecting these covers 
back down to the period two setting of Example 1. If the free
strand in the cover is chosen not to be isotopic to a periodic
braid, then it can be shown via a simple argument that 
some projection of the free strand down to the period two 
case has nontrivial homotopy index. 
Thus, the simple period two skeleton of Example 1 is the seed for 
an infinite number of braid classes with nontrivial homotopy indices. 
Using the techniques of \cite{KKV}, one can use this fact to show that
any parabolic recurrence relation ($\RR=0$) 
admitting this skeleton is forced to have positive topological entropy:
cf. the related results from the Nielsen-Thurston theory of disc 
homeomorphisms \cite{Boy94}.

%
%
\section{Stabilization and invariance}
\label{stable}

%
\subsection{Free braid classes and the extension operator}
Via the results of the previous section, the homotopy index
is an invariant of the {\em discretized} braid class: keeping 
the period fixed and moving within a connected component of 
the space of relative discretized braids leaves the index
invariant. The {\em topological} braid class, as defined in 
\S\ref{braidspace}, does not have an implicit notion of period. 
The effect of refining the discretization of a topological closed 
braid is not obvious: not only does the dimension of the 
index pair change, the homotopy types of the isolating 
neighborhood and the exit set may change as well upon 
changing the discretization. It is thus perhaps remarkable
that any changes are correlated under the quotient operation: the 
homotopy index is an invariant of the {\em topological} 
closed braid class. 

\begin{figure}[hbt]
\begin{center}
\includegraphics[angle=0,width=4.5in]{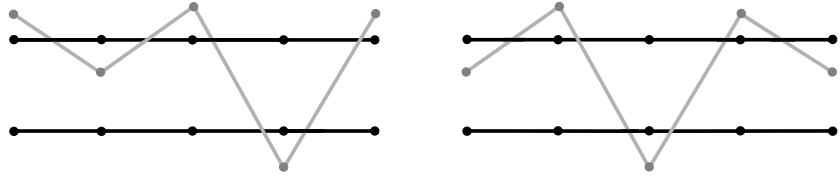}
\caption{An example of two non-free discretized braids which are
	of the same topological braid class but define disjoint 
	discretized braid classes in $\Conf_4^1~\rel~\vv$.}
\label{free} 
\end{center}
\end{figure}

On the other hand, given a complicated braid, it is intuitively
obvious that a certain number of discretization points are 
necessary to capture the topology correctly. 
If the period $d$ is too small $\Conf_d^n~{\rel}~\vv$ may
contain more than one path component with the same
topological braid class:
\begin{definition}
A relative braid class $[\uu~{\rel}~\vv]$ in $\Conf_d^n~{\rel}~\vv$
is called {\em free} if 
\begin{equation} 
	(\Conf_d^n~\rel~\vv) \cap \{\uu~\rel~\vv\} = [\uu~\rel~\vv] ;
\end{equation}
that is, if any other discretized braid in 
$\Conf_d^n~{\rel}~\vv$ which has the same topological braid class as 
$\uu~\rel~\vv$ is in the same discretized braid class $[\uu~\rel~\vv]$.
\end{definition}
A braid class $[\uu]$ is free if the above definition is satisfied 
with $\vv = \emptyset$. Not all discretized braid classes 
are free: see Fig.~\ref{free}. 
%

\begin{figure}[hbt]
\begin{center}
\psfragscanon
\psfrag{a}[][]{\Large (a)}
\psfrag{b}[][]{\Large (b)}
\psfrag{E}[][]{\Large $\E$}
\includegraphics[angle=0,width=5.1in]{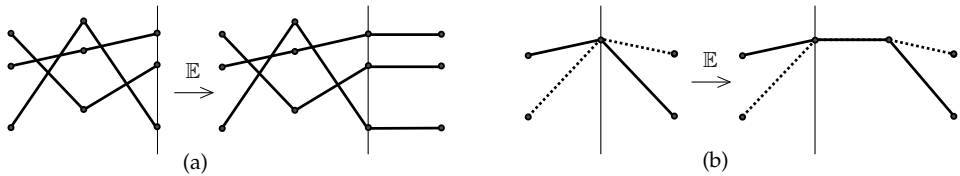}
\caption{(a) The action of $\E$ extends a braid by one period; 
occasionally, (b), $\E$ produces a singular braid. Vertical lines 
denote the $d^{th}$ discretization line.}
\label{fig_E} 
\end{center}
\end{figure}

Define the {\em extension map} $\E:\bar\Conf^n_d\to\bar\Conf^n_{d+1}$
via concatenation with the trivial braid of period one (as in 
Fig.~\ref{fig_E}(a)): 
\begin{equation}
	(\E\uu)^\alpha_i := \left\{
	\begin{array}{cl}
		u_i^{\alpha} &	i=0,\ldots,d 	\\
		u_d^{\alpha} &  i=d+1 .
	\end{array}\right.
\end{equation}
The reader may note (with a little effort) that the non-equivalent 
braids of Fig.~\ref{free} become equivalent under the image of $\E$. 
There are exceptional cases in which $\E\uu$ is a singular braid 
when $\uu$ is not: see Fig.~\ref{fig_E}(b). 
If the intersections at $i=d$ are generic then 
$\E\uu$ is a nonsingular braid. One can always find such a 
representative in $[\uu]$, again denoted by $\uu$.
Therefore the notation $[\E\uu]$ means that 
$\uu$ is chosen in $[\uu]$ with generic intersection at $i=d$. 
The same holds for relative classes $[\E\uu~\rel~\E\vv]$, i.e.
choose $\uu~\rel~\vv \in [\uu~\rel~\vv]$ such
that all intersections of $\uu \cup\vv$ at $i=d$ are generic.

Note that under the action of $\E$ boundedness of a braid class is not
necessarily preserved, i.e. $[\uu~\rel~\vv]$ may be bounded, and
$[\E\uu~\rel~\E\vv]$ unbounded.
For this reason we will prove a stabilization result for
{\em topological} bounded proper braid classes.

\subsection{A topological invariant}
Consider a period $d$ discretized relative braid pair $\uu~\rel~\vv$ 
which is not necessarily free. Collect all (a finite number) of the 
discretized braids $\uu(0),\ldots,\uu(m)$ such that the 
pairs $\uu(j)~\rel~\vv$ are all topologically isotopic to 
$\uu~\rel~\vv$ but not pairwise discretely isotopic.   
For the case of a free braid class, $m=1$. 

\begin{definition}
\label{def_wedge}
Given $\uu~\rel~\vv$ and $\uu(0),\ldots,\uu(m)$ as above, denote
by $\hh(\uu~\rel~\vv)$ the wedge of the homotopy indices of these
representatives,
\begin{equation}
\hh(\uu~\rel~\vv) := \bigvee_{j=0}^{m_d} 
	h(\uu(j)~\rel~\vv),
\end{equation}
where $\vee$ is the topological wedge which, in this context, 
identifies all the constituent exit sets to a single point. 
\end{definition}
This wedge product is well-defined by Theorem~\ref{isolating} by considering
the isolating neighborhood $N = \cup_j {\rm cl}[\uu(j)~\rel~\vv]$.
In general a union of isolating neighborhoods is not necessarily an 
isolating neighborhood again. However, since the word metric strictly 
decreases at $\Sigma$ the invariant set decomposes 
into the union of invariant sets of the individual components of $N$. 
Indeed, if an orbit intersects two components it must have passed 
through $\Sigma$: contradiction.

The principal topological result of this paper is that $\hh$ is an 
invariant of the {\em topological} bounded proper braid class 
$\bigl\{\uu ~\rel~\{\vv\}\bigr\}$.

\begin{theorem}
\label{topinvariant}
Given $\uu~\rel~\vv\in\Conf^n_d~\rel~\vv$ and 
$\tilde\uu~\rel~\tilde\vv\in\Conf^n_{\tilde{d}}~\rel~\tilde\vv$ which are
topologically isotopic as bounded proper braid pairs, then 
\begin{equation}
	\hh(\uu~\rel~\vv) = \hh(\tilde{\uu}~\rel~\tilde{\vv}).
\end{equation}
\end{theorem}

The key ingredients in this proof are that (1) the homotopy
index is invariant under $\E$ (Theorem~\ref{topinvariant}); 
and (2) discretized braids ``converge'' to topological braids 
under sufficiently many applications of $\E$ (Proposition~\ref{freestable}).

%
\begin{theorem}
\label{stabilization}
For $\uu~{\rel}~\vv$ any bounded proper discretized braid pair, the 
wedged homotopy index of Definition~\ref{def_wedge} is invariant
under the extension operator:
\begin{equation}
	\hh(\E\uu~{\rel}~\E\vv) = \hh(\uu~{\rel}~\vv). 
\end{equation}
\end{theorem}

{\it Proof.}
By the invariance of the index
with respect to the skeleton $\vv$, we may assume that $\vv$ is chosen 
to have all intersections generic ($v^\alpha_i\neq v^{\alpha'}_i$ for
all strands $\alpha\neq\alpha'$). Thus, from the proof of 
Lemma~\ref{explicit} in Appendix~\ref{app_1}, we may fix a recurrence 
relation $\RR$ having $\vv$ as fixed point(s) for which $\partial_1\RR_0=0$.

For $\epsilon>0$ consider the one-parameter family of augmented recurrence 
functions\footnote{Recall the indexing conventions: for a period
$d+1$ braid, $u_0^{\tau(\alpha)}=u_{d+1}^\alpha$, and $\RR_0:=\RR_{d+1}$.} 
$\RR^\epsilon=(\RR^\epsilon_i)_{i=0}^{d}$ on braids of period $d+1$:
\begin{equation}
\label{eq_singpert}
\begin{array}{rcl}
	\RR_i^\epsilon (u_{i-1}^\alpha,u_i^\alpha,u_{i+1}^\alpha) 
	&:=& \RR_i(u_{i-1}^\alpha,u_i^\alpha,u_{i+1}^\alpha),
	\quad i=0,..,d-1,\\
\epsilon\cdot\RR_{d}^\epsilon (u_{d-1}^\alpha,u_{d}^\alpha,u_{d+1}^{\alpha}) 
	&:=& u_{d+1}^{\alpha} - u_{d}^\alpha .
\end{array}
\end{equation}
Because of our choice of $\RR_0(r,s,t)=\RR_0(s,t)$ as being independent 
of the first variable, $\RR^\epsilon_0$ is decoupled from the extension of
the braid as $u_{d+1}^\alpha$ wraps around to $u_0^{\tau(\alpha)}$.
By construction the above system satisfies Axioms (A1)-(A2) for
all $\epsilon>0$ with, in particular, the strict monotonicity of (A1) 
holding only on one side. One therefore has a parabolic 
flow $\Psi^t_\epsilon$ on $\bar \Conf_{d+1}^n$ for all $\epsilon>0$. 
In the singular limit $\epsilon=0$, this forces $u_d^\alpha=u_{d+1}^\alpha$, 
and one obtains the flow $\Psi^t_0 = \E\circ\Psi^t$.

Since the skeleton $\vv$ has only generic intersections, 
$\E\vv$ is a nonsingular braid. From Equation (\ref{eq_singpert}), 
all stationary solutions of $\Psi^t$ are stationary solutions for 
$\Psi^t_\epsilon$, i.e., $\Psi^t_\epsilon(\E\vv) = \E\vv$, 
for all $\epsilon\ge 0$. 
Notice that this is not true in general for non-constant solutions.

\begin{figure}[hbt]
\begin{center}
\psfragscanon
\psfrag{M}[][]{\Large $\MM$}
\psfrag{S}[][]{\Large $\Sigma$}
\psfrag{0}[][]{\Large $N_0$}
\psfrag{D}[][]{\Large $D(r)$}
\psfrag{N}[][]{\Large $\NN_{d+1}$}
\psfrag{n}[][]{\Large $\NN_d$}
\psfrag{e}[l][l]{\Large $N_0(2\epsilon C)$}
\psfrag{E}[][]{\Large $\Sigma_{\E}$}
\psfrag{K}[][]{\Large $K$}
\includegraphics[angle=0,width=4.0in]{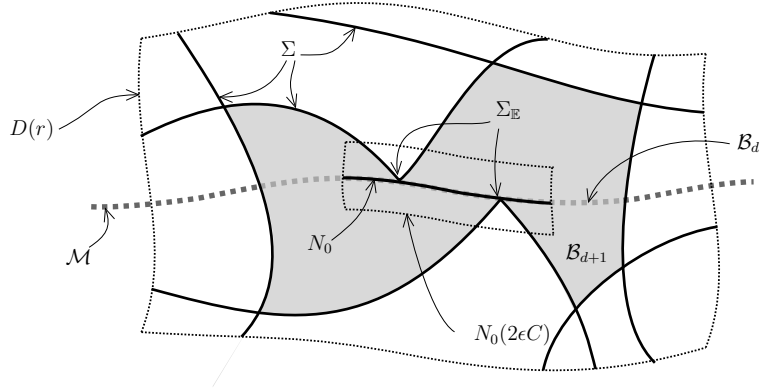}
\caption{The rescaled flow acts on $\NN_{d+1}$, the period $d+1$ braid 
classes. The submanifold $\MM$ is a critical manifold of fixed
points at $\epsilon=0$. Any appropriate isolating neighborhood
$N_0$ in $\NN_d$ thickens to an isolating neighborhood 
$N_0(2\epsilon C)$ which is not necessarily contained in $\NN_{d+1}$.}
\label{fig_stable} 
\end{center}
\end{figure}

Denote by $\NN_{d+1}\subset\Conf^n_{d+1}~\rel~\E\vv$ the subset of
relative braids which are topologically isotopic to $\E\uu~\rel~\E\vv$.
Likewise, denote by $\NN_{d}\subset\bar\Conf^n_{d+1}$ the image under 
$\E$ of the subset of braids in $\Conf^n_d~\rel~\vv$ which 
are topologically isotopic to $\uu~\rel~\vv$.
In other words,
\begin{equation}
	\NN_{d+1} := 
	\bigl\{\E\uu~\rel~\E\vv\bigr\}\cap\Conf^n_{d+1}~\rel~\E\vv
\;\;\;\; ; \;\;\;\;
	\NN_{d} := 
	\E\Bigl(\bigl\{\uu~\rel~\vv\bigr\}\cap\Conf^n_{d}~\rel~\vv\Bigr).
\end{equation}
As per the paragraph
preceding Definition~\ref{def_wedge}, there are a finite number
of connected components of each of these sets. Clearly, 
$\NN_d$ is a codimension-$n$ subset of ${\rm cl}(\NN_{d+1})$. 
Since not all braids in $\bigl\{\uu~\rel~\vv\bigr\}\cap\Conf^n_{d}~\rel~\vv$ 
have generic intersections, the set $\B_d$ may tangentially intersect
the boundary of $\NN_{d+1}$. We will denote this set of 
{\it $\E$-singular braids} by
$\Sigma_{\E} := \partial \NN_{d+1} \cap \NN_d$:
see Fig.~\ref{fig_stable}.

By performing an appropriate change of coordinates (cf. \cite{ConleyFife}),
we can recast the parabolic system $\RR^\epsilon$ as a singular 
perturbation problem. Let
$\x=(x_j)_{j=1}^{nd}$, with $x_{i+1+(\alpha-1)d} := u^{\alpha}_i$, 
and let $\y=(y_\alpha)_{\alpha=1}^n$, 
with $y_\alpha := (u^{\alpha}_{d+1}-u^{\alpha}_{d})$.
Upon rescaling time as $\tau:=t/\epsilon$, the vector field induced
by our choice of $\RR^\epsilon$ is of the form
\begin{equation}
\label{eq_ConleyFife}
\begin{array}{rcl}
\frac{\ds d\x}{\ds d\tau} &=& \epsilon X(\x,\y),\\
\frac{\ds d\y}{\ds d\tau} &=& -\y + \epsilon Y(\x),
\end{array}
\end{equation}
for some (unspecified) vector fields $X$ and $Y$ with the 
functional dependence indicated.
The product flow of this vector field (\ref{eq_ConleyFife}) in the 
new coordinates is denoted by $\Phi^\tau_\epsilon$ and is well-defined on 
$\bar\Conf^{n}_{d+1}$. In the case $\epsilon=0$, the set 
$\MM := \{\y = 0\}\subset\bar\Conf^n_{d+1}$ is a submanifold 
of fixed points containing $\NN_d$ 
for which the flow $\Phi^\tau_0$ is transversally nondegenerate 
(since here $\y'=-\y$). By construction 
${\rm cl}(\NN_d )= {\rm cl}(\NN_{d+1}) \cap \MM$,
as illustrated in Fig.~\ref{fig_stable} (in the simple case where all braid
classes are free and $\NN_{d+1}$ is thus connected).

The remainder of the analysis is a technique in  
singular perturbation theory following \cite{ConleyFife}: one 
relates the $\tau$-dynamics of Equation~(\ref{eq_ConleyFife}) to those
of the $t$-dynamics on $\MM$, whose orbits are of the form $(\x(t),0)$,
where $\x(t)$ satisfies the limiting equation $d\x/dt = X(\x,0)$.
The Conley index theory is well-suited to this situation.

For any compact set $D\subset\MM$ and $r\in\R$, let $D(r) := 
\{(\x,\y)~|~(\x,0) \in D,~~\Vert \y \Vert \le r\}$ denote the
``product'' radius $r$ neighborhood in $\bar\Conf^n_{d+1}$. Denote
by $C=C(D)$ the maximal value $C:=\max_{D} \Vert Y(\x) \Vert$.
Due to the specific form of (\ref{eq_ConleyFife}), we obtain the 
following uniform squeezing lemma.
\begin{lemma}
\label{estimate}
If $S$ is any invariant set of $\Phi^\tau_\epsilon$ contained in some $D(r)$, 
then in fact $S \subset D(\epsilon C)$.
Moreover, 
for all points $(\x,\y)$ with $\x \in D$ and $\Vert \y \Vert = 2\epsilon C$ it holds
that ${d\over d\tau} \Vert y\Vert <0$.
\end{lemma}
{\it Proof.}
Let $(\x,\y)(\tau)$ be an orbit in $S$ contained in some $D(r)$.
Take the inner product of the $\y$-equation with $\y$: 
\begin{eqnarray*}
\langle \frac{\ds d\y}{\ds d\tau},\y \rangle (\tau_0) 
	&=& -\Vert \y(\tau_0) \Vert^2 + \epsilon
		\langle Y(\x(\tau_0)),\y(\tau_0) \rangle,\\
	&\le& - \Vert \y \Vert ^2+\epsilon C \Vert \y \Vert.
\end{eqnarray*}
Hence $\frac{d}{d\tau} \Vert \y \Vert \leq - \Vert \y \Vert  + \epsilon C$,
and we conclude that if $\Vert \y (\tau_0) \Vert > \eps C$ for some $\tau_0 \in \R$, then
$\frac{d}{d\tau} \Vert \y \Vert < 0$. Consequently
$\Vert y (\tau) \Vert $ grows unbounded for $\tau < \tau_0$ and therefore 
$(\x,\y) \not\in S$, a contradiction.
Thus $\Vert y (\tau) \Vert \leq \epsilon C$ for all $\tau \in \R$.

For points $(\x,\y)$ with $\x \in D$ and $\Vert\y \Vert = 2\epsilon C$, 
the above inequality gives that 
$\frac{d}{d\tau} \Vert \y \Vert \leq - \Vert \y \Vert  + \epsilon C < 0$.
\fp
\vsp

By compactness of the proper braid class, it is clear that $\NN_{d+1}$,
and thus the maximal isolated invariant set of $\Phi^\tau_\epsilon$ given by
$S_\epsilon := {\Inv}(\NN_{d+1},\Phi^\tau_\epsilon)$\footnote{Since
$\NN_{d+1}$ is a proper braid class $S_\epsilon$ is contained in 
its interior.}, is strictly contained 
(and thus isolated) in $D(r)$ for some compact $D\subset\MM$ and some $r$ 
sufficiently large. Fix $C:=C(D)$ as above. 
Lemma \ref{estimate} now implies that as $\epsilon$ becomes 
small, $S_\epsilon$ is squeezed into $D(\epsilon C)$ --- a small 
neighborhood of a compact subset $D$ of the critical manifold $\MM$,
as in Fig.~\ref{fig_stable}.\footnote{
If one applies singular perturbation theory it is possible to construct
an invariant manifold $\MM_\epsilon \subset D(\epsilon C)$. 
The manifold $\MM_\epsilon$ lies strictly within $\NN_{d+1}$ and 
intersects $\MM$ at rest points of the $\Phi^t_0$.}

This proximity of $S_\epsilon$ to $\MM$ allows one to compare the 
dynamics of the $\epsilon=0$ and $\epsilon>0$ flows. 
Let $N_0\subset \NN_{d}\subset\MM$ be an isolating neighborhood 
(isolating block with corners) for the maximal $t$-dynamics 
invariant set $S_0 := {\Inv}(\NN_{d},\Psi^t_0)$ within the braid 
class $\NN_d$. Combining Lemma~\ref{estimate} above, Theorem 2.3C 
of \cite{ConleyFife}, and the existence theorems for isolating 
blocks \cite{WilsonYorke}, one concludes that if $(N_0,N_0^-)$ 
is an index pair for the limiting equations $d\x/dt = X(\x,0)$ then 
$N_0(2\epsilon C)$ is an isolating block for $\Phi^t_\epsilon$
for $0<\epsilon \le \epsilon^*(N_0)$ with $\epsilon^*$ sufficiently small.
A suitable index pair for the flow $\Phi^\tau_\epsilon$ of 
Equation~(\ref{eq_ConleyFife}) is thus given by
\begin{equation}
	\left( N_0(2\epsilon C), N_0^-(2\epsilon C)\right) .
\end{equation}
Clearly, then, the homotopy index of $S_0$ is equal to the 
homotopy index of $\Inv(N_0(2\epsilon C))$ for all $\epsilon$ 
sufficiently small. It remains to show that this captures the 
maximal invariant set $S_\epsilon$.

\begin{lemma}
\label{squeeze}
For all sufficiently small $\epsilon$, 
$\Inv(N_0(2\epsilon C),\Phi^\tau_\epsilon)=S_\epsilon$.
\end{lemma}
{\it Proof.} 
By the choice of $D$ it holds that
$S_\epsilon  \subset D(2\epsilon C)$.
We start by proving that $S_\epsilon \subset N_0(2\epsilon C)$
for $\epsilon$ sufficiently small.
%
%
%
Assume by contradiction that $ S_{{\epsilon_j}} 
\not \subset N_0(2{{\epsilon_j}} C)$
for some sequence $\epsilon_j\to 0$. Then, since $N_0(2\epsilon C)$ 
is an isolating neighborhood  for $\epsilon \le \epsilon^*$, 
there exist orbits $(\x_{\epsilon_j},\y_{\epsilon_j})$ in $S_{{\epsilon_j}}$
such that $(\x_{\epsilon_j},\y_{\epsilon_j})(\tau_j) 
\in D(2\epsilon_j C) - N_0(2\epsilon_j C)$, for some $\tau_j \in \R$.
Define $(\widetilde\x_{\epsilon_j},\widetilde\y_{\epsilon_j})(\tau)
= (\x_{\epsilon_j},\y_{\epsilon_j})(\tau - \tau_j)$, and
set $(\a_{\epsilon_j},\b_{\epsilon_j})(t) = 
(\widetilde\x_{\epsilon_j},\widetilde\y_{\epsilon_j})(\tau)$.
The sequence $(\a_{\epsilon_j},\b_{\epsilon_j})$ satisfies the equations
\begin{equation}
\frac{d}{dt} \a_{\epsilon_j}
= X(\a_{\epsilon_j},\b_{\epsilon_j}),~~~~\frac{d}{dt} \b_{\epsilon_j}=
-\frac{1}{\epsilon} \b_{\epsilon_j} + Y(\a_{\epsilon_j}).
\end{equation}
By assumption $\|\b_{\epsilon_j}(t)\| \le C \epsilon_j$, and 
$\|\a_{\epsilon_j}\|, \|d\a_{\epsilon_j}/dt\| \le C$, for all $t\in \R$ and 
all $\epsilon_j$. An Arzela-Ascoli argument then yields the existence 
of an orbit $(\a_*(t),0)\subset \NN_d$, with 
$(\a_*(0),0) \in {\rm cl}(\NN_d-N_0)$, satisfying the equation 
$\frac{d\a_*}{dt} = X(\a_*,0)$. By definition, 
$(\a_*,0) \in \Inv(\NN_d)=\Inv(N_0) \subset {\rm int}(N_0)$, a contradiction,
which proves that $S_\epsilon \subset N_0(2\epsilon C)$ for $\epsilon$ 
sufficiently small.

\begin{figure}[hbt]
\begin{center}
\psfragscanon
\psfrag{1}[lt][l]{\LARGE $u_i^\alpha-u_i^{\alpha'}$}
\psfrag{2}[lt][l]{\LARGE $u_{i+1}^\alpha-u_{i+1}^{\alpha'}$}
\includegraphics[angle=0,width=5.0in]{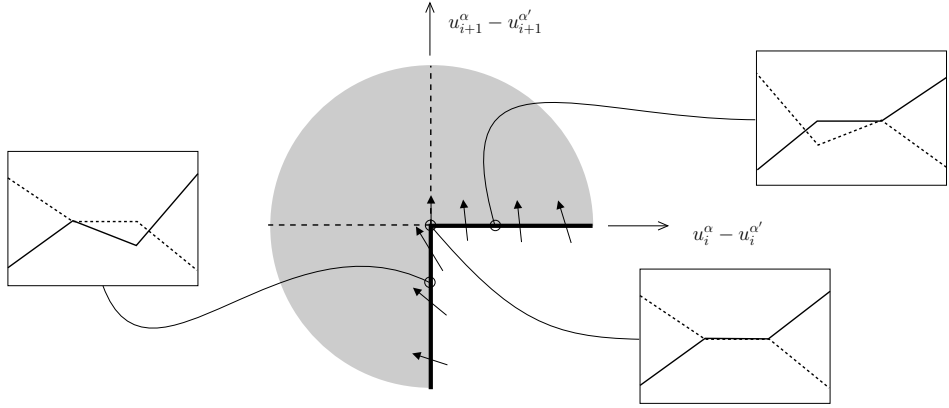}
\caption{The local picture of a generic singular tangency between
strands $\alpha$ (solid) and $\alpha'$ (dashed). The shaded region 
represents $\NN_{d+1}$.}
\label{fig_tangency} 
\end{center}
\end{figure}

The boundary of  $N_0(2\epsilon C)$ splits as $b_1 \cup b_2$, with
$$
b_1 = \{(\x,\y)~|~\Vert \y\Vert =2\epsilon C\}, ~~~~{\rm and}  ~~~~
b_2 = \{(\x,\y)~|~\x \in \partial N_0\}.
$$
Since the compact set $N_0 $ is contained in $\NN_d$, 
the boundary component $b_2$ is  contained
in $\NN_{d+1}$ provided that $\epsilon$ is sufficiently small.
If the set $\Sigma_{\E}$ is non-empty then the boundary component 
$b_1$  never lies entirely in $\NN_{d+1}$
regardless of $\epsilon$.
As $\epsilon \to 0$ the set  $N_0(2\epsilon C) - \bigl(\NN_{d+1} 
\cap N_0(2\epsilon C)\bigr)$
is contained is arbitrary small neighborhood of $\Sigma_{\E}$.
Independent of the parabolic flow in question, and thus of $\epsilon$, 
there exists a neighborhood $K\subset\Sigma_{d+1}^n~\rel~\vv$ of 
$\Sigma_{\E}$ on which the co-orientation of the
boundary is pointed inside the braid class $\NN_{d+1}$. In other words
for every parabolic system the points in $K$ enter $\NN_{d+1}$ 
under the flow, see Fig.~\ref{fig_tangency}.
%
%
By using coordinates $u^{\alpha}_i-u^{\alpha'}_i$ and 
$u^{\alpha}_{i+1}-u^{\alpha'}_{i+1}$ adapted to the singular strands,
it it easily seen (Fig.~\ref{fig_tangency}) that the braids are
simplified by moving into the set $\NN_{d+1}$.

We now show that $\Inv(N_0(2\epsilon C))\subset\NN_{d+1}\cap N_0(2\epsilon C)$.
If not, then there exist points  $(\x_{\epsilon_j},\y_{\epsilon_j}) \in 
\bigl[N_0(2\epsilon_j C)-\bigl(\NN_{d+1}\cap N_0(2\epsilon_j C)\bigr)\bigr] 
\cap\Inv(N_0(2\epsilon_j C))$ for some sequence $\epsilon_j \to 0$.
Consider the $\alpha$-limit sets $\alpha_{\epsilon_j} 
((\x_{\epsilon_j},\y_{\epsilon_j}))$.
Since $(\x_{\epsilon_j},\y_{\epsilon_j}) \in \Inv(N_0(2\epsilon_j C))$, and 
since $\Phi^\tau_{\epsilon_j}((\x_{\epsilon_j},\y_{\epsilon_j}))$ cannot
enter $\NN_{d+1} \cap N_0(2\epsilon_j C)$ in backward time due to 
the co-orientation of $K$, it follows that $\alpha_{\epsilon_j} 
((\x_{\epsilon_j},\y_{\epsilon_j}))$ is contained in 
$N_0(2\epsilon_j C) - \bigl(\NN_{d+1} \cap N_0(2\epsilon_j C)\bigr)$.

By a similar Arzela-Ascoli argument as before, 
this yields a  set $\alpha_0 \subset \Sigma_{\E}$
which is invariant for the flow $\Psi^t_0$.
However due to the form of the vector field the associated  flow $\Psi^t_0$ 
cannot contain an invariant set  in $\Sigma_{\E}$, which proves that 
$\Inv(N_0(2\epsilon C)) \subset \NN_{d+1} \cap N_0(2\epsilon C)$
for $\epsilon$ sufficiently small.

Finally, knowing that $S_\epsilon \subset \Inv(N_0(2\epsilon C))$, 
and that for sufficiently small $\epsilon$ it holds $\Inv(N_0(2\epsilon C) =
\Inv(\NN_{d+1} \cap N_0(2\epsilon C)) = S_\epsilon$, 
it follows that $S_\epsilon = \Inv(N_0(2\epsilon C))$, which proves the lemma.
\fp
\vsp

Theorem~\ref{stabilization} now follows. Since, by Theorem~\ref{isolating}, 
the homotopy index is independent of the parabolic flow
used to compute it, one may choose the parabolic flow $\Phi^\tau_\epsilon$ 
for $\epsilon>0$ sufficiently small. The homotopy index of $\Phi^\tau_\epsilon$
on the maximal invariant set $S_\epsilon$ yields the wedge of 
all the connected components: $\hh(\E\uu~\rel~\E\vv)$.
We have computed that this index is equal to the index of $\Psi^t$ 
on the original braid class: $\hh(\uu~\rel~\vv)$. 
\fp
\vsp

\begin{remark}
{\rm
The proof of Theorem~\ref{stabilization} implies that any component
of the period-$(d+1)$ braid class $\NN_{d+1}$ which does not 
intersect $\MM$ must necessarily have trivial index. 
}
\end{remark}

\begin{remark}
\label{remextra}
{\rm
The above procedure also yields a stabilization result for bounded proper 
classes which are not bounded as topological classes.
In this case one simply augments the skeleton $\vv$ by two constant 
strands as follows. Define the {\em augmented braid} 
$\vv^* := \vv \cup \vv^- \cup \vv^+$, where
\begin{equation}
v^-_{i} := \min_{\alpha,i}{v^\alpha_i}-1,\quad  v^+_{i} := 
	\max_{\alpha,i}{v^\alpha_i}+1 .
\end{equation}
Suppose $[\uu~\rel~\vv] \subset \Conf^n_{d_0} ~\rel~\vv$ is bounded for 
some period $d_0$. It now holds that 
$h(\uu~\rel~\vv) = h(\uu~\rel~\vv^*)$,
and $\bigl\{\uu~\rel~\{\vv^*\}\bigr\}$ is a bounded class.
It therefore follows from Theorem \ref{topinvariant} that
\begin{eqnarray}\label{extra}
\bigvee_{j=0}^{m_{d_0}} 
	h\bigl(\uu(j)~\rel~\vv\bigr)=
\hh(\uu~\rel~\vv^*),
\end{eqnarray}
where $\hh$ can be evaluated via any discrete representative of  
$\bigl\{\uu~\rel~\{\vv^*\}\bigr\}$ of any admissible period.
}
\end{remark}

\subsection{Eventually free classes}
At the end of this subsection, we complete the proof of 
Theorem~\ref{topinvariant}. The preliminary step is to show 
that discretized braid classes are eventually free under $\E$.

Given a braid $\uu\in\Conf^n_d$, consider the extension $\E\uu$ 
of period $d+1$.
Assume at first the simple case in which $d=1$, so that $\E\uu$ 
is a period-2 braid. Draw the braid diagram $\braid(\E\uu)$ 
as defined in \S\ref{braidspace} in the domain $[0,2]\times\R$.
Choose any 1-parameter family of curves 
$\gamma_s: t\mapsto (f_s(t),t) \in (0,2)\times\R$ 
such that $\gamma_0:t\mapsto(1,t)$ and so that $\gamma_s$ is 
transverse\footnote{At the anchor points, the transversality should
be topological as opposed to smooth.} to the braid diagram $\beta(\E\uu)$ 
for all $s$. Define the braid $\gamma_s\cdot\E\uu$ as follows:
\begin{equation}
	(\gamma_s\cdot\E\uu)_i^\alpha := \left\{
	\begin{array}{cl}
	(\E\uu)_i^\alpha & : i = 0,2 \\
	\gamma_s\cap(\E\uu)^\alpha & : i = 1
	\end{array}\right. .
\end{equation}
The point $\gamma_s\cap(\E\uu)^\alpha$ is well-defined 
since $\gamma_s$ is always transverse to the braid strands and $\gamma_0$
intersects each strand but once. 

\begin{lemma}
\label{sweep}
For any such family of curves $\gamma_s$, $[\gamma_s\cdot\E\uu]=[\E\uu]$. 
\end{lemma}
{\em Proof.}
It suffices to show that this path of braids does not 
intersect the singular braids $\Sigma$. Since $\uu$ 
is assumed to be a nonsingular braid, every crossing of two strands 
in the braid diagram of $\E\uu$ is a transversal crossing between 
$i=0$ and $i=1$. Thus, if for some $s$, 
$\gamma_s(t)\cap(\E\uu)^\alpha = \gamma_s(t)\cap(\E\uu)^{\alpha'}$
for distinct strands $\alpha$ and $\alpha'$, then 
\begin{equation}
	\left(\E\uu^\alpha_0 - \E\uu^{\alpha'}_0\right)
	\left(\E\uu^\alpha_1 - \E\uu^{\alpha'}_1\right)  < 0 .
\end{equation}
The braid $\gamma_{s}\cdot\E\uu$ has a crossing of the $\alpha$ and 
$\alpha'$ strands at $i=1$. Checking the transversality of this 
crossing yields
\begin{equation}
\begin{array}{l}
\left((\gamma_{s}\cdot\E\uu)_0^\alpha - 
(\gamma_{s}\cdot\E\uu)_0^{\alpha'}\right)
\left((\gamma_{s}\cdot\E\uu)_2^\alpha - 
(\gamma_{s}\cdot\E\uu)_2^{\alpha'}\right)
\\
= 
\left((\E\uu)_0^\alpha - (\E\uu)_0^{\alpha'}\right)
\left((\E\uu)_2^\alpha - (\E\uu)_2^{\alpha'}\right)
\\
= 
\left((\E\uu)_0^\alpha - (\E\uu)_0^{\alpha'}\right)
\left((\E\uu)_1^\alpha - (\E\uu)_1^{\alpha'}\right) < 0.
\end{array}
\end{equation}
Thus the crossing is transverse and the braid is never singular.
\fp
\vsp

Note that the proof of Lemma~\ref{sweep} does not require the 
braid $\E\uu$ to be a closed braid diagram since the isotopy fixes the 
endpoints: the proof is equally 
valid for any localized region of a braid in which one spatial segment 
has crossings and the next segment has flat strands. 

\begin{corollary}
The ``shifted'' extension operator which inserts a trivial period-1 braid
at the $i^{th}$ discretization point in a braid has the same action 
on components of $\Conf_d$ as does $\E$. 
\end{corollary}

\begin{figure}[htb]
\begin{center}
\psfragscanon
\psfrag{i}[][]{\Large $\sigma_i$}
\psfrag{j}[][]{\Large $\sigma_j$}
\psfrag{0}[][]{\Large $\sigma_i$}
\psfrag{1}[][]{\Large $\sigma_{i+1}$}
\includegraphics[angle=0,width=5.0in]{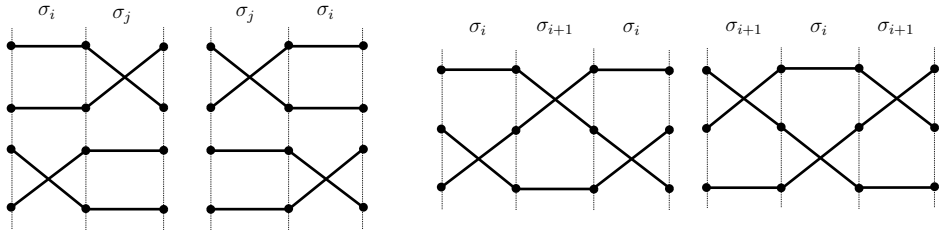}
\caption{Relations in the braid group via discrete isotopy.}
\label{fig_disciso}
\end{center}
\end{figure}

\begin{proposition}
\label{freestable}
The period-$d$ discretized braid class $[\uu]$ is free when
$d>\word{\uu}$. 
\end{proposition}
{\em Proof.}
We must show that any braid $\uu'\in\Conf^n_{d}$ which has the 
same topological type as $\uu$ is discretely isotopic to $\uu$. 
Place both $\uu$ and $\uu'$ in general position so as to record
the sequences of crossings using the generators of the $n$-strand 
positive braid semigroup, $\{\sigma_i\}$, as in \S\ref{braidspace}. 
Recall the braid group has relations $\sigma_i\sigma_j=\sigma_j\sigma_i$ 
for $\vert i-j\vert>1$ and $\sigma_i\sigma_{i+1}\sigma_i = 
\sigma_{i+1}\sigma_i\sigma_{i+1}$; closure requires making conjugacy 
classes equivalent. 

The conjugacy relation can be realized by a discrete isotopy as follows:
since $d>\word{\uu}$, $\uu$ must possess some discretization interval
on which there are no crossings. Lemma~\ref{sweep} then implies that
this interval without crossings commutes with all neighboring discretization
intervals via discrete isotopies. Performing $d$ consecutive exchanges
shifts the entire braid over by one discretization interval. This
generates the conjugacy relation.

To realize the remaining braid relations in a discrete isotopy,
assume first that $\uu$ and $\uu'$ are of the form that there is 
at most one crossing per discretization interval. It is then 
easy to see from Fig.~\ref{fig_disciso} that the braid relations 
can be executed via discrete isotopy.

In the case where $\uu$ (and/or $\uu'$) exhibits multiple crossings
on some discretization intervals, it must be the case 
that a corresponding number of other discretization intervals do not 
possess any crossings (since $d>\word{\uu}$).
Again, by inductively utilizing Lemma~\ref{sweep},
we may redistribute the intervals-without-crossing and ``comb'' out
the multiple crossings via discrete 
isotopies so as to have at most one crossing per 
discretization interval.
\fp

\vsp
{\em Proof of Theorem~\ref{topinvariant}:} Assume that 
$\bigl\{\uu~\rel~\{\vv\}\bigr\}=\bigl\{\uu'~\rel~\{\vv'\}\bigr\}$. This 
implies that there is a path of topological braid diagrams taking
the pair $(\uu,\vv)$ to $(\uu',\vv')$. This path may be 
chosen so as to follow a sequence of standard relations for
closed braids. From the proof of Proposition~\ref{freestable},
these relations may be performed by a discretized isotopy
to connect the pair $(\E^j\uu,\E^j\vv)$ to $(\E^k\uu',\E^k\vv')$ 
for $j$ and $k$ sufficiently large, and of the right relative
size to make the periods of both pairs equal. For this choice, then,  
$\bigl[\E^j\uu~\rel~[\E^j\vv]\bigr] = \bigl[\E^k\uu'~\rel~[\E^k\vv']\bigr]$, 
and their homotopy indices agree. An application of Theorem~\ref{stabilization} 
completes the proof.
\fp 
\vsp

We suspect that all braids in the image of $\E$ are free: a result
which, if true, would simplify index computations yet further.


%
\section{Duality}
\label{duality}

For purposes of computation of the index, we will often pass
to the homological level.
In this setting, there is a natural duality made possible by the 
fact that the index pair used to compute the index of a braid class
can be chosen to be a manifold pair.

\begin{definition}
The {\em duality operator} on discretized braids is the map 
$\Dual:\bar\Conf_{2p}^n\to\bar\Conf_{2p}^n$ given by 
\begin{equation}
	(\Dual\uu)^\alpha_i  :=  (-1)^iu^\alpha_i .
\end{equation}
\end{definition}

\begin{figure}[hbt]
\begin{center}
\psfragscanon
\psfrag{0}[][]{$i=0$}
\psfrag{1}[][]{$i=1$}
\psfrag{2}[][]{$i=2$}
\psfrag{3}[][]{$i=3$}
\psfrag{4}[][]{$i=4$}
\psfrag{D}[][]{$\Dual$}
\includegraphics[angle=0,width=5.1in]{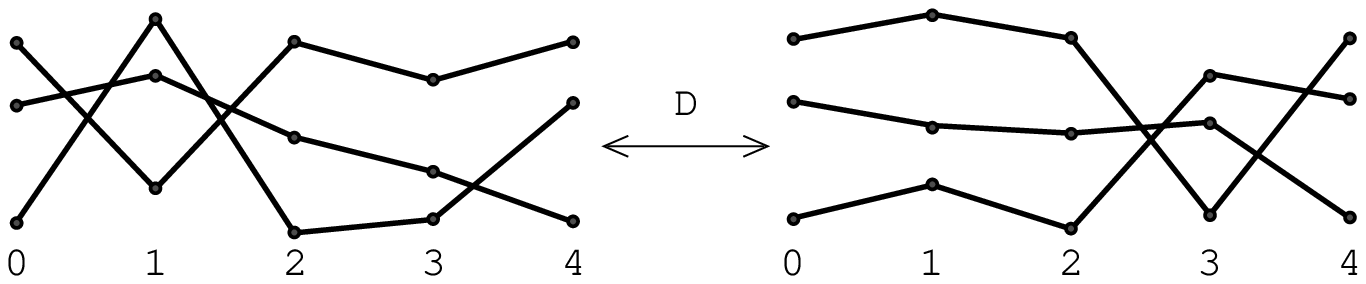}
\caption{The topological action of $\Dual$.}
\label{fig_duals} 
\end{center}
\end{figure}

Clearly $\Dual$ induces a map on relative braid diagrams by 
defining $\Dual(\uu~\rel~\vv)$ to be $\Dual\uu~\rel~\Dual\vv$. 
The topological action of $\Dual$ is to insert a half-twist 
at each spatial segment of the braid. This has the effect of linking 
unlinked strands, and, since $\Dual$ is an involution, 
linked strands are unlinked by $\Dual$: see Fig.~\ref{fig_duals}.

 For the duality statements to follow, we assume that 
all braids considered have even periods and that all of the 
braid classes and their duals are proper, so that the homotopy index
is well-defined.

\begin{lemma}
\label{lem_duality}
The duality map $\Dual$ respects braid classes: if 
$[\uu]=[\uu']$ then $[\Dual(\uu)]=[\Dual(\uu')]$. Bounded
braid classes are taken to bounded braid classes by $\Dual$.
\end{lemma}
\pf
It suffices to show that the map $\Dual$ is a homeomorphism 
on the pair $(\bar\Conf_{2p}^n,\Sigma)$. 
This is true on $\bar\Conf^n_{2p}$ since $\Dual$ is a 
smooth involution ($\Dual^{-1}=\Dual$). If $\uu\in\Sigma$ with
$u_i^\alpha=u_i^{\alpha'}$ and
\begin{equation}
\label{eq_crossing}
(u_{i-1}^{\alpha}-u_{i-1}^{\alpha'})(u_{i+1}^{\alpha}-u_{i+1}^{\alpha'})
\geq 0 ,	
\end{equation}
then applying the operator $\Dual$ yields points 
$\Dual u_i^\alpha=\Dual u_i^{\alpha'}$ with each term in the 
above inequality multiplied by $-1$ (if $i$ is even) 
or by $+1$ (if $i$ is odd): in either case, the quantity 
is still non-negative and thus $\Dual\uu\in\Sigma$. Boundedness
is clearly preserved.
\fp

\begin{theorem}
\label{thm_duality}

\begin{enumerate}
\item[(a)]	The effect of $\Dual$ on the index pair is to 
	reverse the direction of the parabolic flow.
\item[(b)]	For $[\uu~\rel~\vv]\subset\Conf^n_{2p}~\rel~\vv$ of period 
	$2p$ with $n$ free strands,	
	\begin{equation}\label{dual}
		CH_*(h(\Dual(\uu~\rel~\vv));\R) \cong 
		CH_{2np-*}(h(\uu~\rel~\vv);\R) .
	\end{equation}
\item[(c)]	For $[\uu~\rel~\vv]\subset\Conf^n_{2p}~\rel~\vv$ 
         of period $2p$ with $n$ free strands,	
	\begin{equation}\label{dual2}
		CH_*(\hh(\Dual(\uu~\rel~\vv));\R) \cong 
		CH_{2np-*}(\hh(\uu~\rel~\vv);\R).
	\end{equation}
\end{enumerate}
\end{theorem}
{\it Proof:}
For (a), let $(N,N^-)$ denote an index pair associated to a proper relative
braid class $[\uu~\rel~\vv]$. Dualizing sends $N$ to a homeomorphic
space $\Dual(N)$. The following local argument shows that the exit 
set of the dual braid class is in fact the complement (in the boundary)
of the exit set of the domain braid: specifically, 
\[
	(\Dual(N))^- = {\rm cl}\left\{\partial(\Dual(N))-\Dual(N^-)\right\} .
\]

Let $\ww\in[\uu~\rel~\vv]\cap\Sigma$. At any singular anchor point
of $\ww$, \ie, where $w_i^\alpha=w_i^{\alpha'}$ and the 
transversality condition is not satisfied, then it follows from 
Axiom (A2) that
\begin{equation}
\sign\left\{\frac{d}{dt}(w_i^{\alpha} - w_i^{\alpha'})\right\}
= \sign\left\{w_{i-1}^{\alpha} - w_{i-1}^{\alpha'}\right\} .
\end{equation}
(Depending on the form of (A2) employed, one might 
use $w_{i+1}^{\alpha}-w_{i+1}^{\alpha'}$ on the
right hand side without loss.) Since the subscripts on the left side 
have the opposite parity of the subscripts on the right side,
taking the dual braid (which multiplies the anchor points by $(-1)^i$
and $(-1)^{i-1}$ respectively) alters the sign of the terms. 
Thus, the operator $\Dual$ reverses the direction of the  
parabolic flow.

 From this, we may compute the Conley index of the dual braid by 
reversing the time-orientation of the flow. Since one can 
choose the index pair used to compute the index to be an oriented manifold
pair (specifically, an isolating block: see, \eg, \cite{WilsonYorke}), 
one may then apply a Poincar\'e-Lefschetz duality argument as in 
\cite{McCord} and use the fact that the dimension is $2np$ to obtain the
duality formula for homology. This yields (b).

The final claim (c) follows from (b) by showing that 
$\Dual$ is bijective on {\em topological} braid
classes within $\bar\Conf^n_{2p}$. 
Assume that $[\uu~\rel~\vv]$ and $[\uu'~\rel~\vv]$ are 
distinct braid classes in $\Conf^n_{2p}$ of the same topological 
type. Since $\Dual$ is a homeomorphism on $\Conf^n_{2p}$, the 
dual classes $[\Dual\uu~\rel~\Dual\vv]$ and $[\Dual\uu'~\rel~\Dual\vv]$
are distinct. Claim (c) follows upon showing that these
duals are still topologically the same braid class.

Proposition~\ref{freestable} implies that $[(\E^{2k}\uu)~\rel~(\E^{2k}\vv)]
= [(\E^{2k}\uu')~\rel~(\E^{2k}\vv)]$ for $k$ sufficiently large  
since $\{\uu~\rel~\vv\}=\{\uu'~\rel~\vv\}$. By Lemma~\ref{lem_duality},
\[
\Dual\left[(\E^{2k}\uu)~\rel~(\E^{2k}\vv)\right] =  
\Dual\left[(\E^{2k}\uu')~\rel~(\E^{2k}\vv)\right] ,
\]
which, by Lemma~\ref{braidtype} means that these braids are topologically
the same. The topological action of dualizing the $2k$-stabilizations
of $\uu~\rel~\vv$ and $\uu'~\rel~\vv$ is to add $k$ full twists.
Since the full twist is in the center of the braid group (this 
element commutes with all other elements of the braid group 
\cite{Birman}), one can factor the dual braids within the topological
braid group and mod out by $k$ full twists, yielding that 
$\left\{\Dual\uu~\rel~\Dual\vv\right\}
=	\left\{\Dual\uu'~\rel~\Dual\vv\right\}$.
\fp
\vsp

We use this homological duality to complete a crucial computation
in the proof of the forcing theorems (e.g., Theorem~\ref{H}) at the 
end of this paper. The following small corollary uses duality to 
give the first step towards answering the question of just what the
homotopy index measures topologically about a braid class. Recall 
the definition of an augmented braid from Remark~\ref{remextra}.

\begin{corollary}
Consider the dual of any augmented
proper relative braid. Adding a full twist to this dual braid shifts 
the homology of the index up by two dimensions.
\end{corollary}
\pf
Assume that $\Dual[\uu~\rel~\vv^*]$ is the dual of an augmented braid
in period $2p$ (the augmentation is required to keep the braid class bounded
upon adding a full twist). The prior augmentation implies that the outer two
strands of $\Dual\vv$ ``maximally link'' the remainder of the 
relative braid. The effect of adding a full twist to this 
braid can be realized by instead stabilizing $[\uu~\rel~\vv^*]$
twice and then dualizing. The homological duality implies that
for each connected component of the topological class,
\begin{equation}
\label{shift}
\begin{array}{rcl}
	CH_*(h(\Dual\E^2(\uu~\rel~\vv^*))) 
&\cong&		CH_{2np+2-*}(h(\E^2(\uu~\rel~\vv^*))) \\
&\cong&		CH_{2np+2-*}(h(\uu~\rel~\vv^*)) \\
&\cong&		CH_{*-2}(h(\Dual(\uu~\rel~\vv^*))),
\end{array} 
\end{equation}
which gives the desired result for the
index $\hh$ via Theorem \ref{thm_duality}.
\fp

\begin{remark}\label{shift2}
{\em
The homotopy version of \rmref{shift} can be achieved by following a similar
procedure as in \S\ref{stable}. One obtains a double-suspension of the 
homotopy index, as opposed to a shift in homology. 
}
\end{remark}

\begin{remark} {\em
Given a braid class $[\uu]$ of odd period $p=2d+1$, the image under
$\Dual$ is {\em not} necessarily a discretized braid at all: without
some symmetry condition, the braid will not ``close up'' at the ends.
To circumvent this, define the dual of $\uu$ to be the braid
$\Dual(\uu^2)$ --- the dual of the period $2p$ extension of 
$\uu$. The analogue of Theorem~\ref{thm_duality} above 
is that 
\begin{equation}
	CH_*(\hh(\Sym(\Dual(\uu~\rel\vv)));\R)
	\cong CH_{np-*}(\hh(\uu~\rel~\vv);\R) ,
\end{equation}
where $\Sym$ denotes the subset of the braid class which consists
of symmetric braids: $u^\alpha_i=u^\alpha_{2p-i}$ for all $i$.
}\end{remark}


%
%
\section{Morse theory}
\label{VI}
\label{Morseth}\label{morseth}

It is clear that the Morse-theoretic content of the homotopy index
on braids holds implications for the dynamics of parabolic flows and thus
zeros of parabolic recurrence relations. With this in mind, we 
restrict ourselves to bounded proper braid classes.

Recall that the {\em characteristic polynomial} of an index pair
$(N,N^-)$ is the polynomial 
\begin{equation}\label{betti}
	CP_t(N) := \sum_{k\geq 0} \beta_k t^{k}; \quad\quad
	\beta_k(N) := \dim CH_k(N;\real) = \dim H_k(N,N^-;\real).
\end{equation}
The {\em Morse relations} in the setting of the Conley index
(see \cite{ConleyZehn2}) state that, if $N$ has a Morse decomposition 
into distinct isolating subsets $\{N_a\}_{a=1}^C$, then 
\begin{equation}
\label{Morserelations}
	\sum_{a=1}^C CP_t(N_a) = CP_t(N) + (1+t)Q_t ,
\end{equation}
for some polynomial $Q_t$ with {\em nonnegative} integer coefficients. 

\subsection{The exact, nondegenerate case}
For parabolic recurrence relations which satisfy (A3) (gradient type)
it holds that if $h(\uu~\rel~\vv)\neq 0$, then $\RR$ has at 
least one fixed point in $[\uu~\rel~\vv]$. Indeed, one has:
\begin{lemma}
For an exact nondegenerate parabolic flow on a bounded proper 
relative braid class, the sum of the Betti numbers $\beta_k$ of $h$, as
defined in \rmref{betti},  is a lower bound on the 
number of fixed points of the flow on that braid class. 
\end{lemma}
\pf 
The details of this standard Morse theory argument are provided for the 
sake of completeness. Choose $\Psi^t$ a nondegenerate gradient parabolic 
flow on $[\uu~\rel~\vv]$ (in particular, $\Psi^t$ fixes $\vv$ for all time). 
Enumerate the [finite number of] fixed points $\{\uu_a\}_{a=1}^C$ of $\Psi^t$ 
on this [bounded] braid class. By nondegeneracy, the fixed point set
may be taken to be a Morse decomposition  of $\Inv(N)$. The characteristic
polynomial of each fixed point is merely $t^{\mu^*(\uu_a)}$, where 
$\mu^*(\uu_a)$ is the Morse co-index of $\uu_a$. Substituting $t=1$ into 
Equation~(\ref{Morserelations}) yields the lower bound
\begin{equation}\label{morseI}
\#{\mbox{Fix}}([\uu~\rel~\vv],\Psi^t) \ge \sum_k\beta_k(h) .
\end{equation}
\fp

On the level of the topological braid invariant $\hh$, one needs
to sum over all the path components as follows. 
As in Theorem \ref{topinvariant}, choose period-$d$ representatives 
$\uu(j)$ ($j$ from 0 to $m$) for each path component of the topological
class $\bigl\{ \uu~\rel~\{\vv\}\bigr\}$. If we consider fixed 
points in the union $\cup_{j=0}^{m}[\uu(j)~\rel~\vv]$, we obtain
the following Morse inequalities from \rmref{morseI} and 
Theorem~\ref{topinvariant}:
\begin{equation}\label{morseII}
\#{\mbox{Fix}}(\cup_{j=0}^{m}[\uu(j)~\rel~\vv],\Psi^t) \ge \sum_k\beta_k(\hh) ,
\end{equation}
where $\beta_k(\hh)$ is the $k^{th}$ Betti number of $\hh(\uu~\rel~\vv^*)$.
Thus, again, the sum of the Betti numbers is a lower bound, with the
proviso that some components may not contain any critical points.

If the topological class $\bigl\{ \uu~\rel~\{\vv\}\bigr\}$ is bounded the inequality
\rmref{morseII}  holds with the invariant $\hh(\uu~\rel~\vv)$.

\subsection{The exact, degenerate case}
Here a coarse lower bound still exists.

\begin{lemma}\label{mono}
For an arbitrary exact parabolic flow on a bounded relative braid 
class, the number of fixed points is bounded below by the number 
of distinct nonzero monomials in the characteristic polynomial $CP_t(h)$. 
\end{lemma}
\pf
Assuming that $\#{\mbox{Fix}}$ is finite, all critical points are 
isolated and form a Morse decomposition of $\Inv(N)$. The specific 
nature of parabolic recurrence relations reveals that the dimension of  
the null space of the linearized matrix at an isolated critical point 
is at most 2, see e.g. \cite{vanMoer}. Using this fact Dancer proves 
\cite{Dancer}, via the degenerate version of the Morse lemma due to 
Gromoll and Meyer, that $CH_k(\uu_a) \not = 0$ for {\it at most} one 
index $k = k_0$. Equation~(\ref{Morserelations}) implies that, 
\begin{equation}
\sum_{a=1}^{C_d} CP_t(\uu_a) \ge CP_t(h) 
\end{equation}
on the level of polynomials. As the result of Dancer 
implies that for each $a$, $CP_t(\uu_a) = A t^k$, for some $A\ge 0$, 
it follows that the number of critical points needs to be {\it at least} 
the number of non-trivial monomials in $CP_t(h)$.
\fp


As before, if we instead use the topological invariant $\hh$ for 
$\bigl\{ \uu~\rel~\{\vv^*\}\bigr\}$ we obtain that the number of 
monomials in $CP_t(\hh)$ is a lower bound for the total sum of 
fixed points over the topologically equivalent path-components. 

More elaborate estimates in some cases can be obtained via the extension of the Conley index
due to Floer \cite{Floer}.

\subsection{The non-exact case}
If we consider parabolic recurrence relations that are not necessarily 
exact, the homotopy index may still provide information about solutions 
of $\RR=0$. This is more delicate because of the possibility of periodic 
solutions for the flow $u_i' = \RR_i(u_{i-1},u_i,u_{i+1})$.
For example, if $CP_t(h)~{\rm mod}~(1+t) =0$, the index does not 
provide information about additional solutions for $\RR=0$, as a simple 
counterexample shows. However, if 
$
CP_t(h) ~{\rm mod}~(1+t) \not = 0,
$
then there exists at least one solution
of $\RR=0$ with the specified relative braid class. Specifically,

\begin{lemma}
An arbitrary parabolic flow on a bounded relative braid 
class is forced to have a fixed point if $\chi(h):=CP_{-1}(h)$ is 
nonzero. If the flow is nondegenerate, then the number of 
fixed points is bounded below by the quantity
\begin{equation}
	\left.\left( CP_t(h)~{\rm mod}_{\Z^+[t]}~(1+t)\right)\right\vert_{t=1}
\end{equation}
\end{lemma}
\pf
Set $N = {\rm cl}([\uu~\rel~\vv])$. As the vector field $\RR$ has no 
zeros at $\partial N$, the Brouwer degree, ${\rm deg}(\RR,N,0)$,
may be computed via a small perturbation $\widetilde \RR$ and is given
by\footnote{We choose to define the degree via $-d\widetilde\RR$ in order to 
simplify the formulae.}
\[
{\rm deg}(\RR,N,0) := \sum_{\uu \in N, \widetilde \RR(\uu)=0} 
	{\rm sign}~\det(-d\widetilde\RR(\uu)).
\]
For a generic perturbation $\widetilde\RR$ the associated
parabolic flow $\widetilde \Psi^t$ is a Morse-Smale flow \cite{FusOliva}.
The (finite) collection of rest points $\{\uu_a\}$ 
and periodic orbits $\{\gamma_b\}$ of $\widetilde \Psi^t$
then yields a Morse decomposition of ${\rm Inv}(N)$, and the Morse 
inequalities are
\[	\sum_a CP_t(\uu_a) + \sum_b CP_t(\gamma_b) = CP_t(h) + (1+t)Q_t. \]
The indices of the fixed points are given by 
$CP_t(\uu_a) = t^{\mu^*(\uu_a)}$, where $\mu^*$
is the number of eigenvalues of $d\widetilde \RR(\uu_a)$ 
with positive real part, and the indices of periodic orbits are 
given by $CP_t(\gamma_b) = (1+t)t^{\mu^*(\gamma_b)}$.
Upon substitution of $t=-1$ we obtain
\begin{eqnarray*}
{\rm deg}(\RR,N,0) &=& {\rm deg}(\widetilde \RR,N,0) 
= \sum_a (-1)^{\mu^*(\uu_a)}\\
&=&\sum_a CP_{-1}(\uu_a) = CP_{-1}(h) = \chi(h).
\end{eqnarray*}
Thus, if the Euler characteristic of $h$ is non-trivial, 
then $\RR$ has at least one zero in $N$.

In the generic case the Morse relations
give even more information. One has $CP_t(h) = p_1(t) + (1+t)p_2(t)$,
with $p_1,p_2 \in \Z_+[t]$, and $CP_t(h) ~{\rm mod}_{\Z^+[t]}~(1+t) = p_1(t)$.
It then follows that $\sum_{a} CP_t(\uu_a)
\ge CP_t(h) ~{\rm mod}_{\Z^+[t]}~(1+t)$, proving the stated
lower bound.
\fp

The obvious extension of these results to the full index
$\hh$ is left to the reader.



\section{Second order Lagrangian systems}\label{second}

In this final third of the paper, we apply the developed machinery to 
the problem of forcing closed characteristics in {\em second} order 
Lagrangian systems of twist type. The vast literature on fourth order 
differential equations coming from second order Lagrangians includes 
many physical models in nonlinear elasticity, nonlinear optics,
physics of solids, Ginzburg-Landau equations, etc. (see \S\ref{prelude}).
In this context we mention the work of \cite{ABK,Mizel,PT1,PT2}.

\subsection{Twist systems}\label{VIII1}
We recall from \S\ref{I} that closed characteristics at an energy level 
$E$ are concatenations of monotone laps between minima and maxima 
$(u_i)_{i \in \Z}$, which are periodic sequences with even period $2p$. The 
extrema are restricted to the set $\U$, whose connected components are denoted
by $I_E$: interval components (see \S\ref{prelude2} for the precise 
definition). The problem of finding closed characteristics can, in most cases, 
be formulated as a finite dimensional variational problem on the extrema 
$(u_i)$. The following {\em twist hypothesis}, introduced in \cite{VV1}, 
is key:
\begin{itemize}
\item[{\bf (T):}] {\it $\inf \{ J_E[u] =
	\int_0^\tau\bigl( L(u,u_x,u_{xx})+E\bigr) dx \,|\,
	u \in X_\tau(u_1,u_2),\, \tau \in \R^+ \}$  has a minimizer
	$u(x;u_1,u_2)$ for all $(u_1,u_2) \in \{I_E \times I_E~|~ u_1 \not =
	u_2\}$, and $u$ and $\tau$ are $C^1$-smooth functions   of
	$(u_1,u_2)$.}
\end{itemize}
Here $X_\tau=X_\tau(u_1,u_2)=\{u \in C^2([0,\tau])~|~ u(0)
=u_1,~u(\tau)=u_2,~u_x(0)=u_x(\tau)=0~{\rm and}~ u_x|_{(0,\tau)} >0\}.$

Hypothesis (T) is a weaker version of the
hypothesis that assumes that the monotone
laps between extrema are unique (see, \eg, \cite{Kwap1, Kwap2, VV1}).
Hypothesis (T) is valid for large classes of Lagrangians $L$.
For example, if $L(u,v,w) = \frac{1}{2} w^2 + K(u,v)$, the following 
two inequalities ensure the validity of (T):
\begin{enumerate}
\item[{ (a)}]
$\frac{\partial K}{\partial v} v - K(u,v) -E \leq 0$, and 
\item[{ (b)}] 
$\frac{\partial^2 K}{\partial {v}^2} |v|^2 - \frac{5}{2} 
\bigl\{\frac{\partial K}{\partial v} v - K(u,v) -E \bigr\} \geq 0$  
for all $u \in I_E$ and $v \in \R$.  
\end{enumerate}
Many physical models, such as the Swift-Hohenberg equation \rmref{SH}, 
meet these requirements, although these conditions are not always met.
In those cases numerical calculations still predict the validity of (T), which
leaves the impression that the results obtained for twist systems
carry over to many more systems for which Hypothesis (T)  
is hard to check.\footnote{Another method to implement the ideas 
used in this paper is to set up a curve-shortening flow for second 
order Lagrangian systems in the $(u,u')$ plane.} 
For these reasons twist systems play a important 
role in understanding second order Lagrangian systems.
For a direct application of this see \cite{KV}. 

The existence of minimizing laps is valid under very mild hypotheses on
$K$ (see \cite{KV}). In that case 
(b) above is enough to guarantee the validity of (T).
An example of a Lagrangian that satisfies (T), but not (a) is given
by the Erickson beam-model \cite{K1,PT1,TrusZan} $L(u,u_x,u_{xx}) =
\frac{\alpha}{2} |u_{xx}|^2 + \frac{1}{4} (|u_x|^2-1)^2 + \frac{\beta}{2} u^2$.

\subsection{Discretization of the variational principle}\label{VIII2}
We commence by repeating the underlying variational principle for
obtaining closed characteristics as described in \cite{VV1}.
In the present context a {\it broken geodesic} is a $C^2$-concatenation
of monotone laps (alternating between increasing and decreasing laps)
given by Hypothesis (T). A closed characteristic $u$ at energy level $E$
is a ($C^2$-smooth) function $u:~[0,\tau] \to \R$, $0<\tau<\infty$,
which is stationary for the action $J_E[u]$ with respect to variations
$\delta u \in C^2_{\rm per}(0,\tau)$, and $\delta \tau \in \R^+$, and
as such is a `smooth broken geodesic'. 

The following result, a translation of results implicit in \cite{VV1},
is the motivation and basis for the applications of the machinery in 
the first two-thirds of this paper.
\begin{theorem}
\label{discretization}
Extremal points $\{u_i\}$ for bounded solutions of second order Lagrangian 
twist systems are solutions of an exact parabolic recurrence relation 
with the constraints that (i) $(-1)^{i} u_i<(-1)^{i}u_{i+1}$; 
and (ii) the recurrence relation blows up along any sequence satisfying 
$u_i=u_{i+1}$. 
\end{theorem}
{\em Proof:}
For simplicity, we restrict to the case of a nonsingular energy
level $E$: for singular energy levels, a slightly more involved 
argument is required. 
Denote by $I$ the interior of $I_E$, and by $\Delta(I) = \Delta :=
\{(u_1,u_2) \in I\times I~|~u_1=u_2 \}$ the diagonal. Then define 
the {\em generating function}
\begin{equation}
	S:(I\times I)-\Delta \to \R \quad ;\quad
	S(u_1,u_2) := 	\int_0^\tau\bigl( L(u,u_x,u_{xx})+E\bigr) dx;
\end{equation}
the action of the minimizing lap from $u_1$ to $u_2$.
That $S$ is a well-defined function is the content of Hypothesis (T). 
The {\em action functional} associated to $S$ for a period $2p$ system
is the function 
\[	W_{2p}(\uu) := \sum_{i=0}^{2p-1} S(u_i,u_{i+1})	.	\]
Several properties of $S$ follow from \cite{VV1}:
\begin{enumerate}
\item[(a)] [{\it smoothness}]
$S \in C^2(I \times I \backslash \Delta)$.

\item[(b)] [{\it monotonicity}]
$\partial_1 \partial_2 S (u_1,u_2) >0$ for all $u_1\not=u_2\in I$.

\item[(c)] [{\it diagonal singularity}]
${\displaystyle \lim_{u_1\nearrow u_2} -\partial_1 S(u_1,u_2) =
\lim_{u_2\searrow u_1} \partial_2 S(u_1,u_2) =  }$

${\displaystyle \lim_{u_1\searrow u_2} \partial_1 S(u_1,u_2) =
 \lim_{u_2\nearrow u_1} -\partial_2 S(u_1,u_2) = +\infty}$.
\end{enumerate}
In general the function $\partial_1S(u_1,u_2)$
is strictly increasing in $u_2$ for all $u_1 \le u_2 \in  I_E$,
and similarly $\partial_2S(u_1,u_2)$ is strictly increasing in $u_1$.
The function $S$ also has the additional property that $S|_\Delta \equiv 0$.

Critical points of $W_{2p}$ satisfy the exact recurrence relation
\begin{eqnarray}\label{rec}
	\RR_i(u_{i-1},u_i,u_{i+1}):=\partial_2 S(u_{i-1},u_i) +
	\partial_1 S(u_i,u_{i+1})=0,
\end{eqnarray}
where $\RR_i(r,s,t)$ is both well-defined and $C^1$ on the domains
\begin{eqnarray*}
	\Omega_i &=& \{
	(r,s,t) \in I^3~|~  (-1)^{i+1} (s-r)>0,~ (-1)^{i+1} (s-t)>0\},
\end{eqnarray*}
by Property (a). The recurrence function $\RR$ is periodic with $d=2$, 
as are the domains $\Omega$.\footnote{We could also work
with  sequences $\uu$ that satisfy $(-1)^{i+1}(u_{i+1}-u_i)>0$.}
Property (b) implies that Axiom (A1) is satisfied.
Indeed, $\partial_1 \RR_i = \partial_1 \partial_2 S(u_{i-1},u_i) >0$, and 
$\partial_3 \RR_i =  \partial_1 \partial_2 S(u_i,u_{i+1}) >0$.

Property (c) provides information about the behavior of $\RR$ at
the diagonal boundaries of $\Omega_i$, namely, 
\begin{equation}
\begin{array}{rcl}
\label{diagsing}
 	\lim_{s \searrow r}\RR_i(r,s,t) =
 	\lim_{s \searrow t}\RR_i(r,s,t) = +\infty \\
 	\lim_{s \nearrow r}\RR_i(r,s,t) =
 	\lim_{s \nearrow t}\RR_i(r,s,t) = -\infty
\end{array}
\end{equation}
\fp
\vsp

The parabolic recurrence relations generated by second order Lagrangians
are defined on the constrained polygonal domains $\Omega_i$. 
\begin{definition}
A parabolic recurrence relation
is said to be of up-down type if  \rmref{diagsing} is satisfied.
\end{definition}

In the next subsection we demonstrate that the up-down recurrence relations 
can be embedded  into the standard theory as developed in 
\S\ref{braidspace}-\S\ref{VI}. 

\subsection{Up-down restriction}\label{IX}

The variational set-up for second order Lagrangians introduces
a few complications into the scheme of parabolic recurrence 
relations as discussed in \S\ref{braidspace}-\S\ref{VI}.
The problem of boundary conditions will be considered in the following
section. Here, we retool the machinery to deal with the fact that 
maxima and minima are forced to alternate. Such braids we call
{\em up-down} braids.\footnote{The more natural term {\em alternating}
has an entirely different meaning in knot theory.}

\subsubsection{The space $\EE$}\label{IX2}
%
\begin{definition}\label{up-down}
The spaces of general/nonsingular/singular up-down braid diagrams 
are defined respectively as:
\begin{eqnarray*}
\bar\EE_{2p}^n &:=& \bar\DD_{2p}^n \cap \left\{ \uu \, : \,
(-1)^i \bigl( u_{i+1}^\alpha - u_i^\alpha \bigr)> 0 \quad \forall i, \alpha
\right\} ,
\\
\EE_{2p}^n &:=& \DD_{2p}^n \cap \left\{ \uu \, : \,
(-1)^i \bigl( u_{i+1}^\alpha - u_i^\alpha \bigr)> 0 \quad \forall i, \alpha
\right\} ,
\\
\Sigma_{\EE} &:=& \bar\EE_{2p}^n - \EE_{2p}^n.
\end{eqnarray*}
Path components of $\EE_{2p}^n$ comprise the up-down braid types $[\uu]_\EE$,
and path components in $\EE_{2p}^n~\rel~\vv$ comprise the relative up-down
braid types $[\uu~\rel~\vv]_\EE$.
\end{definition}

The set $\bar\EE_{2p}^n $ has a boundary in $\bar\DD_{2p}^n$ 
\begin{equation}
\partial \bar\EE_{2p}^n = \partial\Bigl(\bar\DD_{2p}^n \cap \left\{ \uu \, : \,
(-1)^i \bigl( u_{i+1}^\alpha - u_i^\alpha \bigr)\geq 0 \quad \forall i, \alpha
\right\}\Bigr)
\end{equation}
Such braids, called {\it horizontal singularities}, are not included
in the definition of $\bar\EE_{2p}^n$ since the recurrence relation
\rmref{rec} 
does {\em not} induce a well-defined flow on the boundary 
$\partial \bar\EE_{2p}^n$. 

\begin{lemma}
\label{lem_Boundary}
For any parabolic flow of up-down type on $\bar\EE_{2p}^n$, the 
flow blows up in a neighborhood of $\partial\bar\EE_{2p}^n$ in such a manner 
that the vector field points into $\bar\EE_{2p}^n$. All of the conclusions
of Theorem~\ref{isolating} hold upon considering 
the $\epsilon$-closure of braid classes $[\uu~\rel~\vv]_\EE$ in
$\bar\EE_{2p}^n$, denoted 
$$
{\rm cl}_{\bar\EE,\epsilon} [\uu~\rel~\vv]_\EE:= \Bigl\{\uu~\rel~\vv\in
{\rm cl}_{\bar\EE} [\uu~\rel~\vv]_\EE~:~(-1)^i \bigl( u_{i+1}^\alpha - u_i^\alpha \bigr)
\geq \epsilon \quad \forall i, \alpha\Bigr\},
$$ 
for all $\epsilon>0$ sufficiently small.
\end{lemma}
{\em Proof:}
The proof that any parabolic flow $\Psi^t$ of up-down type 
acts here so as to strictly decrease the word metric at singular
braids is the same 
proof as used in Proposition~\ref{word}. The only difficulty 
arises in what happens at the boundary of $\bar\EE_{2p}^n$: we must 
show that $\Psi^t$ respects the up-down restriction in forward time.

Define the function 
$$
\epsilon (\uu) = \min_{i,\alpha} |u_i^\alpha- u_{i+1}^\alpha|.
$$
Clearly, if $\epsilon (\uu)=0$, then $\uu \in \partial \bar\EE_{2p}^n$.
Let $\uu \in \bar\EE_{2p}^n$, and consider the evolution $\Psi^t(\uu)$, $t>0$.
We compute ${d\over dt} \epsilon( \Psi^t(\uu))$ as $\epsilon(\Psi^t(\uu))$ becomes small.
Using \rmref{rec} it follows that 
\begin{eqnarray*}
{d\over dt} (u_i^\alpha -  u_{i+1}^\alpha )&=& \RR_i(u_{i-1}^\alpha,u_{i}^\alpha,u_{i+1}^\alpha)
- \RR_{i+1}(u_{i}^\alpha,u_{i+1}^\alpha,u_{i+2}^\alpha) \to
\infty,\\
&~&~~~~\hbox{as}~~~u_i\searrow u_{i+1}, ~~~(i~~{\rm odd}),\\
{d\over dt} (u_{i+1}^\alpha -  u_{i}^\alpha ) &=& \RR_{i+1}(u_{i}^\alpha,u_{i+1}^\alpha,u_{i+2}^\alpha) 
- \RR_i(u_{i-1}^\alpha,u_{i}^\alpha,u_{i+1}^\alpha) \to
\infty,\\
&~&~~~~\hbox{as}~~~u_i\nearrow u_{i+1}, ~~~(i~~{\rm even}).
\end{eqnarray*}
These inequalities show that ${d\over dt} \epsilon( \Psi^t(\uu))>0$ as
soon as  $\epsilon(\Psi^t(\uu))$ becomes  too small.
Due to the boundedness of $[\uu~\rel~\vv]_\EE$ 
and the infinite repulsion at $\partial \bar \EE_{2p}^n$, 
we can choose a uniform $\epsilon(\uu~\rel~\vv)>0$
so that  
${d\over dt} \epsilon( \Psi^t(\uu))>0$ for $\epsilon(\Psi^t(\uu)) \le \epsilon(\uu~\rel~\vv)$, and
thus
${\rm cl}_{\bar\EE,\epsilon}[\uu~\rel~\vv]_\EE$ is an isolating neighborhood for all
$0<\epsilon\le \epsilon(\uu~\rel~\vv)$.
\fp
\subsubsection{Universality for up-down braids}\label{IX3}
We now show that the topological information contained in
up-down braid classes can be continued to the canonical case
described in \S\ref{braidspace}.
As always, we restrict attention to proper, bounded braid classes, 
proper being defined as in Definition~\ref{proper}, and bounded 
meaning that the set $[\uu~ \rel ~\vv]_\EE$ is bounded in $\bar \DD_{2p}^n$.
Note that an up-down braid class $[\uu~\rel~\vv]_\EE$ can sometimes 
be bounded while $[\uu~\rel~\vv]$ is not.
To bounded proper up-down braids we assign a homotopy index. From 
Lemma \ref{lem_Boundary} it follows that for $\epsilon$ sufficiently small
the set $N_{\EE,\epsilon}:={\rm cl}_{\bar\EE,\epsilon}[\uu~\rel~\vv]_\EE$ 
is an isolating neighborhood in $\bar\EE_{2p}^n$ whose Conley index,
$$
	h(\uu~\rel~\vv,\EE) := h(N_{\EE,\epsilon}) ,
$$
is well-defined with respect to any parabolic flow $\Psi^t$ generated 
by a parabolic recurrence relation of up-down type, and is independent of $\epsilon$.
As before, non-triviality of 
$h(N_{\EE,\epsilon})$ implies existence of a non-trivial
invariant set inside $N_{\EE,\epsilon}$ (see \S\ref{updownmorse}).

The obvious question is what relationship holds between the homotopy 
index $h(\uu~\rel~\vv,\EE)$ and that of a braid class without the 
up-down restriction. To answer this, augment the skeleton $\vv$ as follows:
define $\vv^* = \vv \cup \vv^- \cup \vv^+$, where
$$
v^-_{i} := \min_{\alpha,i}{v^\alpha_i} -1 +(-1)^{i+1},\quad  
v^+_{i} := \max_{\alpha,i}{v^\alpha_i} +1 +(-1)^{i+1}.
$$
The topological braid class $\{\uu~\rel~\vv^*\}$ is bounded and proper. 
Indeed, boundedness follows from adding the strands $\vv^\pm$ which 
bound $\uu$, since 
$\min_{\alpha,i} v_i^\alpha \le u_i^\alpha \le \max_{\alpha,i} v_i^\alpha$.
Properness is satisfied since $\{\uu~\rel~\vv\}$ is proper.
\begin{theorem}\label{equiv}
For any bounded proper up-down braid class $[\uu~\rel~\vv]_\EE$ in 
$\EE_{2p}^n~\rel~\vv$, 
$$
	h(\uu~\rel~\vv,\EE) = h(\uu~ \rel~\vv^*).
$$
\end{theorem}

{\it Proof.} From Lemma~\ref{explicit}  in Appendix A we obtain  
a parabolic recurrence relation $\RR^0$ (not necessarily up-down type) 
for which $\vv^*$ is a solution.
We denote the associated parabolic flow by $\Psi^t_0$.
Define two functions $k_1$ and $k_2$ in $C^1(\R)$, with
$k_1'\geq 0\geq k_2'$, and
$k_1(\tau) =0$ for $\tau\le -2\delta$, $k_1(-\delta)\ge K$, and
$k_2(\tau) =0$ for $\tau\ge 2\delta$, $k_2(\delta)\ge K$, for
some $\delta>0$ and $K>0$ to be specified later.
Introduce a new recurrence function
$\RR_i^1 (r,s,t) = \RR_i^0(r,s,t) + k_2 (s-r) + k_1(t-s)$ for $i$ odd, and
$\RR_i^1 (r,s,t) = \RR_i^0(r,s,t) - k_1 (s-r) - k_2(t-s)$ for $i$ even.
The associated parabolic flow will be denoted by $\Psi^t_1$, and
$\Psi^t_1(\vv^*) = \vv^*$ by construction by choosing $\delta$ sufficiently small.
Indeed, if we choose $\delta < \epsilon (\vv)$, the augmented skeleton is a
fixed point for $\Psi^t_1$.

Since the braid class $[\uu~\rel~\vv^*]$ is bounded and proper,
$N_1 = {\rm cl} [\uu~\rel~\vv^*]$ is an isolating neighborhood
with invariant set {\sc Inv}$(N_1)$. If we choose $K$ large enough, and
$\delta$ sufficiently small, then the invariant set {\sc Inv}$(N_1)$ lies
entirely in ${\rm cl}_{\bar\EE,\epsilon} [\uu~\rel~\vv^*]_\EE = 
{\rm cl}_{\bar\EE,\epsilon} [\uu~\rel~\vv]_\EE = N_{\EE,\epsilon}$.
Indeed, for large $K$ we have that for each $i$, $\RR_i^1 (r,s,t)$
has a fixed sign on the complement of  $N_{\EE,\epsilon}$.
Therefore, $ h(\uu~\rel~\vv^*) = h(N_1) = h(N_{\EE,\epsilon})$.
Now restrict the flow 
$\Psi^t_1$ to $N_{\EE,\epsilon} \subset \bar\EE_{2p}^n~\rel~\vv$.
We may now construct a homotopy between $\Psi^t_1$ and $\Psi^t$, via
$(1-\lambda)\RR + \lambda \RR^1$ (see Appendix A), where
$\RR$ and the associated flow $\Psi^t$ are defined by \rmref{rec}.
The braid $\vv^*$ is stationary along the homotopy and therefore
$$
	h(N_1) =  h(N_{\EE,\epsilon},\Psi^t_1) = h(N_{\EE,\epsilon},\Psi^t),
$$
which proves the theorem.
\fp
\vskip.2cm
We point out that similar results can be proved for other domains $\Omega_i$ 
with various boundary conditions. The key observation is that the up-down 
constraint is really just an addition to the braid skeleton.
\subsubsection{Morse theory}\label{updownmorse}
For bounded proper up-down braid classes  $[\uu~\rel~\vv]_\EE$ 
the Morse theory of \S\ref{VI} applies. Combining this with 
Lemma~\ref{lem_Boundary} and Theorem~\ref{equiv}, the topological information 
is given by the invariant $\hh$ of the topological braid type 
$\bigl\{\uu~\rel~\{\vv^*\}\bigr\}$.
\begin{corollary}\label{udmorse}
On bounded proper up-down braid classes, the total number of fixed points of an 
exact parabolic up-down recurrence relation is bounded below by the 
number of monomials in the critical polynomial $CP_t(\hh)$ of
the homotopy index. 
\end{corollary}

{\it Proof.} 
Since all critical point are contained in $N_{\EE,\epsilon}$ the
corollary follows from the Lemmas \ref{mono}, \ref{lem_Boundary} and
Theorem \ref{equiv}.
\fp
\vsp

%
%

\section{Multiplicity of closed characteristics}\label{X}\label{multi}

We now have assembled the tools necessary to prove Theorem~\ref{H}, the 
general forcing theorem for closed characteristics in terms of braids, and
Theorems \ref{H3} and \ref{H4}, the application to singular and
near-singular energy levels. 
Given one or more closed characteristics, we keep track of the braiding
of the associated strands, including at will any period-two shifts. 
Fixing these strands as a skeleton, we add hypothetical
free strands and compute the homotopy index. If nonzero, this index then
forces the existence of the free strand as an existing solution, which, 
when added to the skeleton, allows one to iterate the argument 
with the goal of
producing an infinite family of forced closed characteristics. 

The following lemma (whose proof is straightforward and thus omitted) 
will be used repeatedly for proving existence of closed characteristics. 
\begin{lemma}
\label{coverings}
Assume that $\RR$ is a parabolic recurrence relation on $\Conf_d^n$ with
$\uu$ a solution. Then, for each integer $N>1$, there exists a lifted 
parabolic recurrence relation on $\Conf_{Nd}^n$ for which every lift of 
$\uu$ is a solution. Furthermore, any solution to the lifted dynamics
on $\Conf_{Nd}^n$ projects to some period-$d$ solution.
\footnote{This does not imply a $d$-periodic solution, but merely 
a braid diagram $\uu$ of period $d$.} 
\end{lemma}

The primary difficulties in the proof of the forcing theorems 
are (i) computing the index (we will use all features of the 
machinery developed thus far, including stabilization and duality); 
and (ii) asymptotics/boundary conditions related to the three types 
of closed interval components $I_E$: a compact
interval, the entire real line, and the semi-infinite ray. 

All of the forcing theorems are couched in a little braid-theoretic 
language:
\begin{definition}
\label{intersection}
The {\em intersection number} of two strands $\uu^\alpha$, $\uu^{\alpha'}$
of a braid $\uu$ is the number of crossings in the braid diagram, denoted
\[
	\intnum(\uu^\alpha,\uu^{\alpha'}) := 
	\# {\mbox{ of crossings of strands }} 
\]
The {\em trivial braid} on $n$ strands is any braid (topological or 
discrete) whose braid diagram has no crossings whatsoever,
\ie, $\iota(\uu^\alpha,\uu^{\alpha'})=0$, for all $\alpha,\alpha'$. 
The {\em full-twist braid} on $n$ strands, is the braid of $n$ connected
components, each of which has exactly two crossings with every other 
strand, i.e., $\iota(\uu^\alpha,\uu^{\alpha'})=2$ for all 
$\alpha \not = \alpha'$. 
\end{definition}

Among discrete braids of period two, the trivial braid and the full twist
are duals in the sense of \S\ref{duality}.

\subsection{Compact interval components}\label{X1}
Let $E$ be a regular energy level for which the set $\U$ contains a compact
interval component $I_E$. 

\begin{theorem}\label{A}
Suppose that a twist system with compact $I_E$ possesses one or more closed 
characteristics which, as a discrete braid diagram, 
form a nontrivial braid. Then there exists an infinity of non-simple, 
geometrically distinct closed characteristics in $I_E$.
\end{theorem}

In preparation for the proof of Theorem \ref{A} we state a technical lemma,
whose [short] proof may be found in \cite{VV1}.
\begin{lemma}\label{tech}
Let $I_E = [u_-,u_+]$, then there exists a $\delta_0>0$ such that 
\begin{enumerate}
\item $\RR_1(u_-+\delta,u_-,u_-+\delta)>0$,
$\RR_1(u_+,u_+-\delta,u_+)<0$, and
\item  $\RR_2(u_-,u_-+\delta,u_-)>0$, $\RR_2(u_+-\delta,u_+,u_+-\delta)<0$,
\end{enumerate}
for any $0<\delta\le \delta_0$.
\end{lemma}

{\it Proof of Theorem \ref{A}.}
Via Theorem~\ref{discretization}, finding closed characteristics is
equivalent to solving the recurrence relation given by \rmref{rec}.
Define the domains
\[
\Omega_i^\delta = \left\{ \begin{array}{c} 
\{(u_{i-1},u_i,u_{i+1}) \in I_E^3 ~|~
u_-+\delta <u_{i\pm 1} <u_i-{\delta/2}<u_+-\delta\},~~ i~{\rm odd,}\\
 \{(u_{i-1},u_i,u_{i+1}) \in I_E^3 ~|~
u_-+\delta < u_i+{\delta/2} < u_{i\pm 1} < u_+-\delta\},~~ i~{\rm even},
\end{array}\right.
\]
For any integer $p\ge 1$ denote by $\Omega_{2p}$ the 
set of $2p$-periodic sequences $\{u_i\}$ for
which  $(u_{i-1},u_i,u_{i+1}) \in \Omega_i^\delta$.
By Lemma~\ref{tech}, choosing $0<\delta<\delta_0$ small enough forces 
the vector field $\RR=(\RR_i)$ to be everywhere transverse to 
$\partial \Omega_{2p}$, making $\Omega_{2p}$ positively invariant
for the induced parabolic flow $\Psi^t$.

By Lemma~\ref{coverings}, one can lift the assumed solution(s) to a 
pair of period $2p$ single-stranded solutions to \rmref{rec}, 
$\vv^1$ and $\vv^2$, satisfying $\intnum(\vv^1,\vv^2)\neq 0$, for
some $p\ge 1$. Define the cones
\begin{eqnarray*}
\label{eq_cone}
	C_-&=& \{ \uu\in \Omega_{2p}~|~ 
	u_i\le v_i^\alpha, ~\alpha=1,2\},~~{\rm and}\\
	C_+ &=& \{ \uu \in \Omega_{2p}~|~ 
	u_i\ge v_i^\alpha, ~\alpha=1,2\}.
\end{eqnarray*}
The combination of the facts $\intnum(\vv^1,\vv^2)=r>0$, Axiom (A1), 
and the behavior of $\RR$ on 
$\partial\Omega_{2p}$ implies that on the boundaries of the cones
$C_-$ and $C_+$ the vector field $\RR$ is everywhere transverse and
pointing inward. Therefore, $C_-$ and $C_+$ are also positively invariant 
with respect to the parabolic flow $\Psi^t$. Consequently,
$W_{2p}$ has global maxima $\vv^-$ and $\vv^+$ on ${\rm int}(C_-)$ and ${\rm
int}(C_+)$ respectively. The maxima $\vv^-$ and $\vv^+$ have the property that
$v^-_i < v^\alpha_i < v^+_i$, $\alpha=1,2$.
As a braid diagram, $\vv=\{\vv^1,\vv^2,\vv^-,\vv^+\}$ is a stationary
skeleton for the induced parabolic flow $\Psi^t$.

Having found the solutions $\vv^-$ and $\vv^+$ we now choose a
compact interval $I \varsubsetneq I_E$, such that the skeletal strands
are all contained in $I$. In this way we obtain a proper parabolic flow 
(circumventing boundary singularities) which can be extended to a parabolic 
flow on $\bar\EE_{2p}^1 ~\rel~\vv$.
Let $[\uu~\rel~\vv]_\EE$ be the relative braid class with a period $2p$
free strand $\uu=\{u_i\}$ which links the strands $\vv^1$ and $\vv^2$ 
with intersection number $2q$ while satisfying $v_i^- < u_i < v_i^+$: see 
Fig.~\ref{StrI} below.

\begin{figure}[hbt]
\begin{center}
\psfragscanon
\psfrag{u}[][]{\large $\uu$}
\psfrag{1}[][]{\large $\vv^1$}
\psfrag{2}[][]{\large $\vv^2$}
\psfrag{+}[][]{\large $\vv^+$}
\psfrag{-}[][]{\large $\vv^-$}
\includegraphics[angle=0,height=3cm,width=9cm]{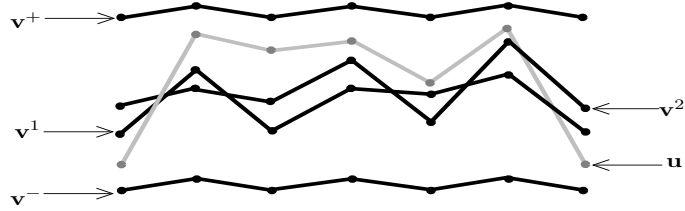}
\caption{A representative braid class for the compact case: $q=1, ~r=4,~2p=6$.}
\label{StrI}
\end{center}
\end{figure}

As an up-down braid class, $[\uu~\rel~\vv]_\EE$ is a bounded proper braid class
provided $0<2q<r\le 2p$, and the Morse
theory  discussed in \S\ref{Morseth} and \S \ref{updownmorse}
then requires the evaluation of the invariant $\hh$ of
the topological class $\bigl\{ \uu~\rel~\{\vv^*\}\bigr\}$.
In this case, since
$h(\uu~\rel~\vv,\EE)=h(\uu~ \rel~\vv^*)=h(\uu~ \rel~\vv)$,  
augmentation is not needed, and $\hh(\uu~\rel~\vv^*) = \hh(\uu~\rel~\vv)$.
The nontriviality of the homotopy index $\hh$ is given
by the following lemma, whose proof we delay until \S\ref{compu}.
\begin{lemma}\label{comp1}
The Conley homology of $\hh(\uu~\rel~\vv)$ is given by:
\begin{equation}
\label{CHbraid1}
CH_k(\hh) = \left\{\begin{array}{ccl}
	 \R &:& k= 2q-1,~2q \\
	  0 &:& {\mbox{ else.}}
	\end{array}\right.
\end{equation}
In particular $CP_t(\hh) = t^{2q-1}(1+t)$.
\end{lemma}
 From the Morse theory of Corollary \ref{udmorse} we derive 
that for each $q$ satisfying $0<2q<r\le 2p$ there exist
at least two distinct period-$2p$ solutions of \rmref{rec}, 
which generically are of index $2q$ and $2q-1$. In this manner, 
the number of solutions depends on $r$ and $p$. To construct 
infinitely many, we consider $m$-fold coverings of the skeleton
$\vv$, i.e., one periodically extends $\vv$ to a skeleton 
contained in $\EE_{2pm}^4$, $m\ge 1$.
Now $q$ must satisfy $0<2q<rm\le 2pm$. By choosing triples $(q,p,m)$
such that $(q,pm)$ are relative prime, we obtain 
the same Conley homology as above, and therefore an infinity of pairs
of geometrically distinct solutions of \rmref{rec}, which, via 
Lemma~\ref{coverings} and Theorem~\ref{discretization} yield 
an infinity of closed characteristics.
\fp
\vskip.2cm
Note that if we set $q_m=q$  and $p_m=pm$, then
the admissible ratios $\frac{q_m}{p_m}$ for finding closed characteristics 
are determined by the relation
\begin{eqnarray}\label{ratio}
	0< \frac{q_m}{p_m} < \frac{r}{2p}.
\end{eqnarray}
Thus if $\vv^1$ and $\vv^2$ are maximally linked, i.e. $r=2p$, 
then closed characteristics exist for all ratios in $\Q \cap (0,1)$.

\subsection{Non-compact interval components: $I_E = \R$}\label{X2}
On non-compact interval components, closed characteristics need not exist.
An easy example of such a system is given by the quadratic Lagrangian
$L= \frac{1}{2} |u_{xx}|^2 + \frac{\alpha}{2} |u_x|^2 + \frac{1}{2} |u|^2 $,
with $\alpha>-2$. Clearly $I_E=\R$ for all $E>0$, and the Lagrangian system
has no closed characteristics for those energy levels.
For $\alpha<-2$ the existence of closed characteristics strongly
depends on the eigenvalues of the linearization around $0$. 
To treat non-compact interval components,
some prior knowledge about asymptotic behavior of the system is needed.
We adopt an asymptotic condition shared by most
physical Lagrangians: {\it dissipativity}.

\begin{definition}
\label{def_dissipative}
A second order Lagrangian system is {\em dissipative} on an interval 
component $I_E=\R$ if there exist pairs $u_1^*< u_2^*$, with $-u_1^*$ and $u_2^*$ arbitrarily large,  such that
\begin{eqnarray*}
	&~&-\partial_1 S(u_1^*, u_2^*) >0, \quad \partial_2 S(u_1^*, u_2^*) >0,
		\quad {\rm and} \quad \\
       &~&~~ \partial_1 S(u_2^*, u_1^*) >0, \quad -\partial_2 S(u_2^*, u_1^*) >0.
\end{eqnarray*}
\end{definition}
Dissipative Lagrangians admit a strong forcing theorem:
\begin{theorem}\label{D}
Suppose that a dissipative twist system with $I_E=\R$ possesses one or more
closed characteristic(s) which, as discrete braid diagram
in the period-two projection, forms a link which is not a full-twist 
(Definition \ref{intersection}). Then there exists an infinity of 
non-simple, geometrically distinct closed characteristics in $I_E$. 
\end{theorem} 
{\it Proof.}
After taking the $p$-fold covering of the period-two projection for some $p\ge 1$, 
the hypotheses imply the existence two sequences $\vv^1$ and $\vv^2$ 
that form a braid diagram in $\EE_{2p}^2$ whose intersection
number is not maximal, i.e. $0\le \intnum(\vv^1,\vv^2) = r <2p$.
Following Definition~\ref{def_dissipative}, choose 
$I= [u_1^*,u_2^*]$,
with $u_1^*<u_2^*$ such that $u_1^* < v^1_i,v^2_i < u_2^*$ for all $i$,
 and let $\Omega_i^\delta$ and $\Omega_{2p}$ be as 
in the proof of Theorem \ref{A},
with $u_1^*$ and $u_2^*$ playing the role of $u_-+\delta$ and $u_+ -\delta$ respectively
for some $\delta>0$ small.
Furthermore define the set 
$$
C := \{ \uu \in \Omega_{2p}~|~\intnum(\uu,\vv^1) = \intnum(\uu,\vv^2)=2p\}.
$$
Since $0\le \intnum(\vv^1,\vv^2)<2p$, the vector field $\RR$ given by 
\rmref{rec} is transverse to $\partial C$. Moreover, the set $C$ is 
contractible, compact, and $\RR$ is pointing outward at the boundary 
$\partial C$. The set $C$ is therefore negatively invariant for the 
induced parabolic flow $\Psi^t$. Consequently, there exists a global 
minimum $\vv^3$ in the interior of $C$. Define the skeleton $\vv$ to be 
$\vv:=\{\vv^1,\vv^2,\vv^3\}$.

\begin{figure}[hbt]
\begin{center}
\psfragscanon
\psfrag{u}[][]{\large $\uu$}
\psfrag{1}[][]{\large $\vv^1$}
\psfrag{2}[][]{\large $\vv^2$}
\psfrag{-}[][]{\large $\vv^3$}
\includegraphics[angle=0,height=3cm,width=9cm]{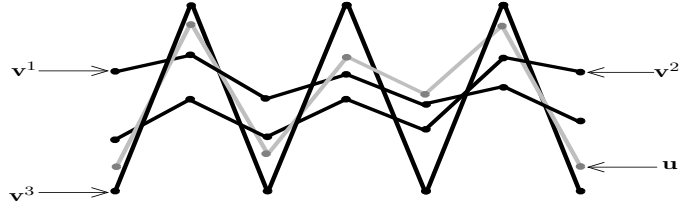}
\caption{A representative braid class for the case $I_E=\R$: $q=2,~r=1,~2p=6$.}
\label{StrII}
\end{center}
\end{figure}

Consider the up-down relative braid class $[\uu~\rel~\vv]_\EE$
described as follows: choose $\uu$ to be a $2p$-periodic strand with 
$(-1)^i u_i \ge (-1)^i v_i^3$, such that $\uu$ has intersection number $2q$ with
each of the strands $\vv^1\cup \vv^2$, $0\le r< 2q<2p$, as in 
Fig.~ \ref{StrII}. For $p\ge 2$, $[\uu~\rel~\vv]_\EE$ is a bounded 
proper up-down braid class. As before, in order to apply the Morse theory 
of Corollary \ref{udmorse}, it suffices to compute 
the homology index of the topological braid class 
$\bigl\{ \uu~\rel~\{\vv^*\}\bigr\}$:
%
\begin{lemma}\label{comp2}
The Conley homology of $\hh(\uu~\rel~\vv^*)$ is given by:
\begin{equation}
\label{CHbraid2}
CH_k(\hh) = \left\{\begin{array}{ccl}
	 \R &:& k= 2q,~2q+1 \\
	  0 &:& {\mbox{ else}}
	\end{array}\right.
\end{equation}
In particular $CP_t(\hh) = t^{2q}(1+t)$.
\end{lemma}
By the same covering/projection argument as in the proof of Theorem~\ref{A}, 
infinitely many solutions are constructed within the admissible ratios
\begin{eqnarray}\label{ratio2}
	 \frac{r}{2p}< \frac{q_m}{p_m} <1.
\end{eqnarray}
\fp

Theorem~\ref{D}  also implies that the existence of a single {\it non-simple}
closed characteristic yields  infinitely many other closed characteristics.
In the case of two unlinked closed characteristics all possible ratios
in $\Q \cap (0,1)$ can be realized.
 
\subsection{Half spaces $I_E \simeq \R^\pm$}\label{X3}

The case $I_E = [\bar u,\infty)$ (or $I_E=(-\infty,\bar u]$) shares much
with both the compact case and the the case $I_E=\R$. Since these 
$I_E$ are non-compact we again impose a dissipativity condition.

\begin{definition}
\label{def_dissipative2}
A second order Lagrangian system is {\em dissipative} on an interval component 
$I_E = [\bar u,\infty)$ 
if there exist arbitrarily large  points $u^*>\bar{u}$  such that
\begin{eqnarray*}
&~&\partial_1 S(\bar u, u^*) >0, \quad \partial_2 S(\bar u, u^*) >0, \quad {\rm and} \quad \\
&~&\partial_1 S(u^*,\bar u) >0, \quad \partial_2 S(u^*,\bar u) >0.
\end{eqnarray*}
\end{definition}
For dissipative Lagrangians we obtain the same general result as Theorem \ref{A}.
\begin{theorem}\label{F}
Suppose that a dissipative twist system with $I_E\simeq\R^\pm$ possesses  
one or more closed characteristics which, as a discrete braid 
diagram, form a nontrivial braid. Then there exists an infinity of non-simple, 
geometrically distinct closed characteristics in $I_E$. 
\end{theorem} 
{\it Proof.}
We will give an outline of the proof since 
the arguments are more-or-less the same as 
in the proofs of Theorems \ref{A} and \ref{D}.
Assume without loss of generality that $I_E = [\bar u,\infty)$.
By assumption there exist two sequences $\vv^1$ and $\vv^2$ which
form a nontrivial braid in $\EE_{2p}^2$, and thus
$0<r=\intnum(\vv^1,\vv^2)\le 2p$. 
Defining the cone $C_-$ as in the proof of Theorem \ref{A} yields a global maximum $\vv^-$
which contributes to the skeleton $\widetilde \vv = \{\vv^1,\vv^2,\vv^-\}$.
Consider the braid class $[\uu~\rel~\widetilde\vv]_\EE$ defined by adding
the strand $\uu$ such that $u_i>v_i^-$ and $\uu$ links with the strands 
$\vv^1$ and $\vv^2$ with intersection number $2q$, $0<2q<r$.

Notice, in contrast to our previous examples, that $[\uu~\rel~\widetilde\vv]_\EE$
is not bounded. In order to incorporate the dissipative boundary condition 
that $u_i \to u^*$ is attracting, we add one additional strand $\vv^+$.
Set $\vv^+_i = \bar u$ for $i$ even, and $\vv^+_i = u^*$, for $i$ odd.
As in the proof of Theorem \ref{D} choose $u^*$ large enough such that
$v^1_i,v_i^2<u^*$.
Let $\RR^\dagger$ be a parabolic recurrence relation such that 
$\RR^\dagger(\vv^+)=0$. Using $\RR^\dagger$ one can construct yet another 
recurrence relation $\RR^{\dagger\dagger}$ which coincides with $\RR$ on  
$[\uu~\rel~\widetilde\vv]_\EE$ and which has $\vv^+$ as a fixed point (use cut-off functions). 
By definition the skeleton $\vv = \{\vv^1,\vv^2,\vv^-,\vv^+\}$ 
is stationary with respect to the recurrence relation $\RR^{\dagger\dagger}=0$.

Now let $ [\uu~\rel~\vv]_\EE$ be as before, with the additional 
requirement that $(-1)^{i+1} u_i < (-1)^{i+1} v_i^+$.
This defines a bounded proper up-down braid class.
The homology index of the topological class $\bigl\{\uu~\rel~\{\vv^*\}\bigr\}$
is given by the following lemma (see \S\ref{compu}).
\begin{lemma}\label{comp3}
The Conley homology of $\hh(\uu~\rel~\vv^*)$ is given by
\begin{equation}
\label{CHbraid3}
CH_k(\hh) = \left\{\begin{array}{ccl}
	 \R &:& k= 2q-1,~2q \\
	  0 &:& {\mbox{ else}}
	\end{array}\right.
\end{equation}
In particular $CP_t(\hh) = t^{2q-1}(1+t)$.
\end{lemma} 
For the remainder of the proof we refer to that of Theorem \ref{A}.
\fp

\begin{figure}[hbt]
\begin{center}
\psfragscanon
\psfrag{u}[][]{\large $\uu$}
\psfrag{1}[][]{\large $\vv^1$}
\psfrag{2}[][]{\large $\vv^2$}
\psfrag{+}[][]{\large $\vv^+$}
\psfrag{-}[][]{\large $\vv^-$}
\includegraphics[angle=0,height=3cm,width=9cm]{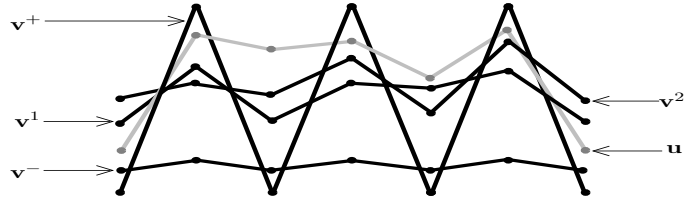}
\caption{A representative braid class for the case $I_E=\R^\pm$: 
	$q=1,~r=4,~2p=6$.}
\label{StrIII}
\end{center}
\end{figure}

\subsection{A general multiplicity result and singular energy levels}\label{X4}

\subsubsection{Proof of Theorem \ref{H}}
Lagrangians for which the above mentioned dissipativity conditions are
satisfied for all (non-compact) interval components at energy $E$,
are called {\it dissipative} at $E$.\footnote{One class of 
Lagrangians that is dissipative on all its regular energy levels 
is described by 
$$
\lim_{\lambda \to \infty} 
\lambda^{-s}L(\lambda u,\lambda^\frac{2+s}{4} v,\lambda^\frac{s}{2} w) = 
c_1  |w|^2 + c_2 |u|^s, {\rm for~some}~s>2,~{\rm and}~c_1, c_2 >0,
$$
pointwise in  $(u,v,w)$.}
For such Lagrangians the results for the three different types of interval 
components are summarized in Theorem \ref{H} in \S\ref{prelude}.
The fact that the presence of a non-simple closed characteristic, when 
represented as a braid, yields
a non-trivial, non-maximally linked braid diagram,  allows us
to apply all three Theorems \ref{A}, \ref{D}, and \ref{F}, proving 
Theorem \ref{H}.

\subsubsection{Singular energy levels}
The forcing theorems in \S \ref{X1} - \S\ref{X3} are applicable for all 
regular energy levels provided the correct configuration of 
closed characteristics can be found a priori.
In this section we will discuss the role of
singular energy levels; they may create
configurations which force the existence of (infinitely) many periodic orbits. 
The equilibrium points in these singular energy levels act as seeds for
the infinite family of closed characteristics.

For singular energy levels the set $\U$ is the union of several interval
components, for which at least one interval component contains an equilibrium 
point. If $\partial^2_u L(u_*,0,0)>0$ at an equilibrium point $u_*$, 
then such a point is called {\it non-degenerate} and is contained in the 
interior of an interval component.
For applying our results of the previous section
the nature of the equilibrium points may play a role.

\subsubsection{Case I: $I_E=\R$}
We examine the case of a singular energy level $E=0$ such that $I_{E}=\R$ and
$I_{E}$ contains at least two equilibrium points.
One observes that if the equilibria can be regarded as 
periodic orbits then Theorem~\ref{D} would apply: a regularization argument
makes this rigorous.  Let $E=0$ be the energy level in which $\U$
is the concatenation of three interval components $(-\infty,a] \cup [a,b]
\cup [b,\infty)$, i.e., the equilibria are $a$ and $b$.  
We remark that the nature of the equilibrium points is irrelevant; there is
a global reason for the existence of two unlinked periodic orbit in the
energy levels $E \in (0,c_0)$ for some small $c_0>0$, see \cite{VV1}. In 
these regular
energy levels we can apply Theorem~\ref{D}, and a limit procedure ensures
that the periodic solutions persist to the degenerate energy level $E=0$,
proving Theorem \ref{H3}.

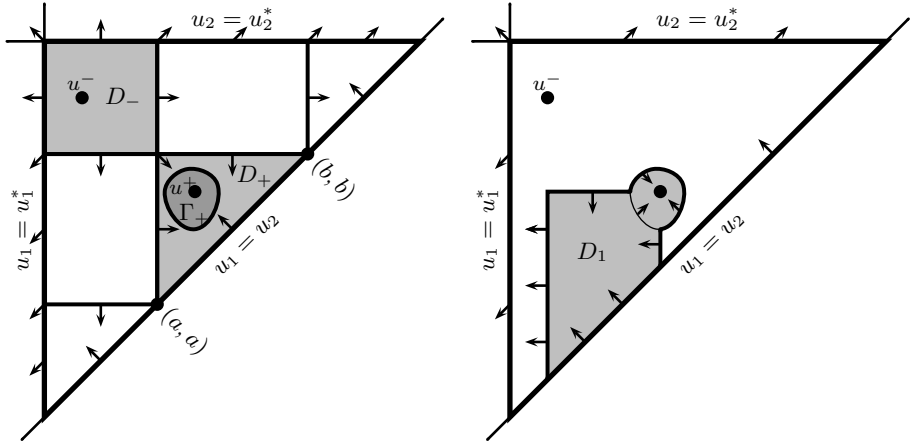
\begin{figure}
\centerline{\psset{xunit=1cm,yunit=1cm}
\begin{pspicture}(-3.1,-3)(3,3.1)
 
  
  \psset{linewidth=1.5pt,fillcolor=grey1,fillstyle=solid}
  \pspolygon(-2.5,1)(-2.5,2.5)(-1,2.5)(-1,1)
  \pspolygon(-1,-1)(-1,1)(1,1)


  \psset{linewidth=1.5pt,fillcolor=grey2}
  \psccurve(-0.9,0.5)(-0.6,0.8)(-0.2,0.6)(-0.5,0)  
 
  \psset{fillstyle=none}
  \pspolygon(-1,1)(-1,2.5)(1,2.5)(1,1)
  \pspolygon(-1,1)(-2.5,1)(-2.5,-1)(-1,-1)
  \psset{linewidth=2pt}
  \pspolygon(-2.5,-2.5)(-2.5,2.5)(2.5,2.5)

  \rput(-0.5,0.2){\small $\Gamma_+$}
  \rput(0.3,0.7){\small $D_+$}
  \rput(-1.45,1.75){\small $D_-$}

  \psset{linewidth=1pt}
  \psline{c-c}(-2.8,-2.8)(2.8,2.8)
  \psline{c-c}(-3,2.5)(2.5,2.5)
  \psline{c-c}(-2.5,3)(-2.5,-2.5)
  \rput(0,2.8){\small$u_2=u_2^*$}
  \rput{90}(-2.8,0){\small$u_1=u_1^*$}
  \rput{45}(0.2,-0.2){\small$u_1=u_2$}
  \psset{dotsize=5pt 0}
  \psdots(-1,-1)(1,1)
  \psdots(-0.5,0.5)(-2,1.75)
  \rput[rt](-0.45,0.7){\footnotesize $u^+$}
  \rput[b](-2,1.85){\footnotesize $u^-$}
  \rput[l]{-45}(1.1,0.9){\small $(b,b)$}
  \rput[l]{-45}(-0.9,-1.1){\small $(a,a)$}

  \psline{->}(-2.5,2.5)(-2.7,2.7)
  \psline{->}(-1,1)(-0.8,0.8)
  \psline{->}(0,0)(-0.2,0.2)

  \psline{->}(-2.5,1.75)(-2.8,1.75)
  \psline{->}(-2.5,1)(-2.7,0.8)
  \psline{->}(-2.5,0)(-2.7,-0.2)
  \psline{->}(-2.5,-1)(-2.7,-1.2)
  \psline{->}(-2.5,-1.75)(-2.7,-1.95)

  \psline{->}(-1.75,2.5)(-1.75,2.8)
  \psline{->}(-1,2.5)(-0.8,2.7)
  \psline{->}(0,2.5)(0.2,2.7)
  \psline{->}(1,2.5)(1.2,2.7)
  \psline{->}(1.75,2.5)(1.95,2.7)

  \psline{->}(-1.75,1)(-1.75,0.7)
  \psline{->}(-1.75,-1)(-1.75,-1.3)
  \psline{->}(-1.75,-1.75)(-1.95,-1.55)
  \psline{->}(-1,0)(-0.7,0)

  \psline{->}(-1,1.75)(-0.7,1.75)
  \psline{->}(1,1.75)(1.3,1.75)
  \psline{->}(1.75,1.75)(1.55,1.95)
  \psline{->}(0,1)(0,0.7)

\end{pspicture}
 \psset{xunit=1cm,yunit=1cm}
\begin{pspicture}(-3.1,-3)(3,3.1)
 
  
  \psset{linewidth=1.5pt,fillcolor=grey1,fillstyle=solid}


  \psccurve(-0.9,0.5)(-0.6,0.8)(-0.2,0.6)(-0.5,0)  
  \pspolygon*[linecolor=grey1,linewidth=0](-0.5,0.5)(-2,0.5)(-2,-2)(-0.5,-0.5)

  \psset{fillstyle=none}
  \psccurve[linewidth=0.5pt](-0.9,0.5)(-0.6,0.8)(-0.2,0.6)(-0.5,0)  
  \psline{c-c}(-0.9,0.5)(-2,0.5)(-2,-2)
  \psline{c-c}(-0.5,-0.5)(-0.5,0)  
  \psset{linewidth=2pt}
  \pspolygon(-2.5,-2.5)(-2.5,2.5)(2.5,2.5)

  \rput(-1.4,-0.3){\small $D_1$}

  \psset{linewidth=1pt}
  \psline{c-c}(-2.8,-2.8)(2.8,2.8)
  \psline{c-c}(-3,2.5)(2.5,2.5)
  \psline{c-c}(-2.5,3)(-2.5,-2.5)
  \rput(0,2.8){\small$u_2=u_2^*$}
  \rput{90}(-2.8,0){\small$u_1=u_1^*$}
  \rput{45}(0.2,-0.2){\small$u_1=u_2$}
  \psset{dotsize=5pt 0}
  \psdots(-0.5,0.5)(-2,1.75)
  \rput[b](-2,1.85){\footnotesize $u^-$}

  \psline{->}(-2.5,2.5)(-2.7,2.7)
  \psline{->}(-0.75,0.75)(-0.6,0.6)
  \psline{->}(-0.25,0.25)(-0.4,0.4)
  \psline{->}(-0.8,0.2)(-0.65,0.35)

  \psline{->}(-2.5,1)(-2.7,0.8)
  \psline{->}(-2.5,-1)(-2.7,-1.2)

  \psline{->}(-1,2.5)(-0.8,2.7)
  \psline{->}(1,2.5)(1.2,2.7)

  \psline{->}(-1.5,-1.5)(-1.7,-1.3)
  \psline{->}(-1,-1)(-1.2,-0.8)
  \psline{->}(-0.5,-0.2)(-0.8,-0.2)
  \psline{->}(-2,-1.5)(-2.3,-1.5)
  \psline{->}(-2,-0.75)(-2.3,-0.75)
  \psline{->}(-2,0)(-2.3,0)
  \psline{->}(-2,0)(-2.3,0)
  \psline{->}(-1.4,0.5)(-1.4,0.2)

  \psline{->}(1,1)(0.8,1.2)

\end{pspicture}
 }
\caption{The gradient of $W_2$ for the case with two equilibria and
dissipative boundary conditions. On the left, for $E=0$, 
the regions $\D_\pm$ with the
maxima and minima $\uu^\pm$ are depicted, as well as the 
superlevel set $\Gamma^+$. On the right, for $E \in (0,c_0)$,
the region $D_1$, containing an
index 1 point, is indicated.}
\label{f:twosimple}
\end{figure}

Recall from~\cite{VV1} that  two equilibrium points
imply the existence of maximum $\uu^+$ and minimum $\uu^-$, both simple
closed characteristics, see Fig.~\ref{f:twosimple}.
Define the regions
$D_+=\{ (u_1,u_2) \,|\, u_2-u_1 > 0,~u_1 \ge a,~u_2\le b\}$,
and $D_- = \{ (u_1,u_2) \,|\, u_1^* \le u_1 \le a,~~b\le u_2\le u_2^*\}$,
where $(u_1^*,u_2^*)$ is the point where the dissipativity condition of
Definition~\ref{def_dissipative} is satisfied.
Then $\uu^+ \in D_+$ and $\uu^- \in D_-$.

Since $W_2$ is a $C^2$-function on $\mbox{int}(D_+)$ it follows from Sard's
theorem that there exists a regular value $e^+$ such that $0\le\max_{\partial
D_+} W_2 < e^+ <\max_{D_+} W_2$.  
Consider the connected component
of the super-level set $\{ W_2 \geq e^+ \}$ which contains $\uu^+$.
The outer
boundary of this component is a  smooth circle and $\nabla W_2$ points
inwards on this boundary circle.
Let $\Gamma^+$ be the {\it interior} of the outer boundary circle 
in question.
%
%
By continuity it follows
that there exists a positive constant $c_0$ such that 
$\Gamma^+$ remains an isolating neighborhood for $E \in (0,c_0)$.
In the following let $E \in (0,c_0)$ be arbitrary.

Define 
$D_1 = \{(u_1,u_2) \,|\, u_2-u_1 \ge \epsilon,~
u_1^-\le u_1 \le u_1^+,~u_2\le u_2^+\} \cup \Gamma^+$.
It follows from the properties of $S$ (see \S\ref{second}) that $D_1$
is again an isolating neighborhood, see Fig.~\ref{f:twosimple}.
It holds that $CP_t(D_1) = 0$, and $\{ D_1\backslash \Gamma^+, \Gamma^+\}$
forms a Morse decomposition. The Morse relations \rmref{Morserelations} 
yield
$$
CP_t(\Gamma^+) + CP_t(D_1\backslash \Gamma^+) =
1+CP_t(D_1\backslash \Gamma^+) = (1+t)Q_t,
$$
where $Q_t$ is a nonnegative polynomial.  This implies that $D_1\backslash
\Gamma^+$ contains an index 1 solution $\uu^1$.  We can now define $D_2=
\{(u_1,u_2) \,|\, u_2-u_1 \ge \epsilon,~ u_1^+\le u_1,~u_2^+ \le u_2\le
u_2^-\} \cup \Gamma^+$.  In exactly the same way we find an index 1
solution $\uu^2 \in D_2$.  Notice, that by construction $\intnum(\uu^1,\uu^2) =
0$.  Theorem~\ref{D} now yields an infinity of closed characteristics for all
$0<E<c_0$.  As described in \S\ref{X2} these periodic solutions can be
characterized by $p$ and $q$, where $(p,q)$ is any pair of integers such
that $q<p$ and $p$ and $q$ are relative prime (or $p=q=1$).  Here $2p$ is
the period of the solution $\uu_{p,q}$ and
$2q=\intnum(\uu_{p,q},\uu^1)=\intnum(\uu_{p,q},\uu^2)$.

In the limit $E \to 0$ the solutions $\uu^1$ and $\uu^2$ may collapse onto
the two equilibrium points (if they are centers).  
Nevertheless, the infinite family of solutions
still exists in the limit $E=0$, because the extrema of the associated
closed characteristics may only  coalesce in pairs at the equilibrium points.  
This follows from the uniqueness of the initial value problem
of the Hamiltonian system.  Hence in the limit $E \to 0$ the type $(p,q)$ of
the periodic solution is conserved when we count extrema \emph{with}
multiplicity and intersections \emph{without} multiplicity.

Note that when the equilibria are saddle-foci then $\uu^1$ and $\uu^2$
stay away from $\pm 1$ in the limit $E \to 0$. Extrema
may still coalesce at the equilibrium points as $E \to 0$, but
intersections are counted with respect to $\uu^1$ and $\uu^2$.
Finally, in the regular energy levels $E \in (0,c_0)$, Theorem~\ref{D} 
provides at least two solutions of each type (except $p=q=1$); in
the limit $E=0$ one cannot exclude the possibility that  two solutions of 
the same type coincide.
\fp


\begin{remark}\label{improper}
{\rm
Theorem \ref{H3}, proved in this subsection, is immediately applicable 
to the Swift-Hohenberg model \rmref{SH} as described in \S\ref{prelude}. 
Notice that if the parameter $\alpha$ satisifies $\alpha>1$, then 
Theorem \ref{H3} yields the existence of infinitely many closed 
chararacteristics at energy $E= -\frac{(\alpha-1)^2}{4}$, and nearby levels.
However from the physical point of view it is also of interest to consider 
the case $\alpha\le 1$. In that case there exists only one singular 
energy level and one equilibrium point. This case can be treated with 
our theory, but the nature of the equilibrium point comes into play. 
If an equilibrium point is a saddle or saddle-focus, 
it is possible that no additional periodic orbits exist 
(see \cite{vdB}). However, if the equilibrium point is a center an 
initial non-simple closed characteristic can be found by analyzing an 
improper braid class, which by, Theorem \ref{H},  then yields infinitely 
many closed characteristics. The techniques involved are very similar to 
those used in the  present and subsequent sections. We do not present
the details here as this falls outside of the scope of this paper.
}
\end{remark}

\subsubsection{Cases II and III: $I_E=[a,b]$ or $I_E=\R^\pm$}
The remaining cases are dealt with in Theorem \ref{H4}. We will restrict 
the proof here to the case that $I_E$ contains an equilibrium point 
that is a saddle-focus --- the center case can be treated as in 
\cite{Angenent1}.\footnote{Indeed, for energy levels $E+c$, $c$ 
sufficiently small,  a small simple closed charecteristic exists
due to the center nature of the equilibrium point at $E$; spectrum 
$\{\pm ai,\pm bi\}$, $a<b$. This small simple closed characteristic will 
have a non-trivial rotation number close to $\frac{a}{b}$. The fact 
that the rotation number is non-zero allows one to use the arguments 
in \cite{Angenent1} to construct a non-simple closed characteristic. 
As a matter of fact a linked braid diagram is created this way.}
It also follows for the previous that there is no real difference 
between $I_E$ being compact or a half-line. For simplicity we consider 
the case that $I_E$ is compact.
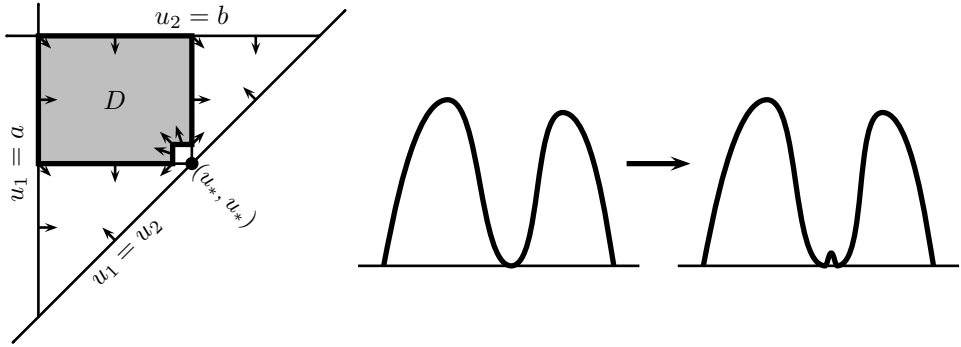
\begin{figure}[h]
\centerline{\psset{xunit=1.7cm,yunit=1.7cm}
\begin{pspicture}(-1.6,-1.5)(6.1,1.35)
 
  
  \psset{linewidth=2pt,fillcolor=grey1,fillstyle=solid}
  \pspolygon(-1.2,0)(-1.2,1)(0,1)(0,0.15)(-0.15,0.15)(-0.15,0)
  \rput(-0.6,0.5){$D$}

  \psset{linewidth=1pt,fillstyle=none}
  \psline{c-c}(0,0)(-0.15,0)
  \psline{c-c}(0,0)(0,0.15)
  \psline{c-c}(-1.4,-1.4)(1.2,1.2)
  \psline{c-c}(-1.45,1)(1,1)
  \psline{c-c}(-1.2,1.25)(-1.2,-1.2)
  \rput(0,1.15){$u_2=b$}
  \rput{90}(-1.35,0){$u_1=a$}
  \rput{45}(-0.5,-0.7){$u_1=u_2$}
  \psdots[dotsize=5pt 0](0,0)
  \rput[l]{-45}(0.03,-0.03){\small $(u_*,u_*)$}

  \psline{->}(-0.6,0)(-0.6,-0.15)
  \psline{->}(-0.6,1)(-0.6,0.85)
  \psline{->}(0.5,1)(0.5,0.85)
  \psline{->}(0,0.5)(0.15,0.5)
  \psline{->}(-1.2,0.5)(-1.05,0.5)
  \psline{->}(-1.2,-0.5)(-1.05,-0.5)
  \psline{->}(0.5,0.5)(0.4,0.6)
  \psline{->}(-0.6,-0.6)(-0.7,-0.5)

  \psline{->}(-1.2,1)(-1.1,0.9)
  \psline{->}(-1.2,0)(-1.1,-0.1)
  \psline{->}(0,1)(0.1,0.9)

  \psline{->}(-0.15,0.15)(-0.25,0.25)
  \psline{->}(-0.15,0.075)(-0.3,0.125)
  \psline{->}(-0.075,0.15)(-0.125,0.3)
  \psline{->}(0,0.15)(0.1,0.25)
  \psline{->}(-0.15,0)(-0.25,-0.1)
  
  \psline(1.3,-0.8)(3.5,-0.8)
  \psline(3.8,-0.8)(6,-0.8)
  \psset{linewidth=2pt,curvature=0.7 0.3 0}
  \psecurve(1.3,-10)(1.5,-0.8)(2,0.5)(2.5,-0.8)(2.9,0.4)(3.3,-0.8)(3.5,-10)
  \psecurve(3.8,-10)(4,-0.8)(4.5,0.5)(4.95,-0.8)(5,-0.7)(5.05,-0.8)(5.4,0.4)(5.8,-0.8)(6,-10)
  \psline{->}(3.4,0)(3.9,0)

\end{pspicture}
 }
\caption{[left] The gradient of $W_2$ for the case of one saddle-focus 
equilibrium and compact boundary conditions. Clearly a saddle point is 
found in $D$. [right] The perturbation of one equilibrium to three equilibria.}
\label{f:onesadfoc}
\end{figure}

Let us first make some preliminary observations. When $u_*$ is a saddle-focus,
then in $E=0$ there exists a solution $\uu^1$ such that $u^1_1 < u_* < u^1_2$.
This follows from the fact that there is a point $(u^*_1,u^*_2)$, 
$u^*_1<u^*_2$, close to $(u_*,u_*)$ at which the vector $\nabla W_2$ 
points to the north-west (see Fig.~\ref{f:onesadfoc} and~\cite{VV1}).
This solution $\uu^1$ is a saddle point, its rotation number being unknown.
The impression is that $(u_*,u_*)$ is a minimum (with $\tau =0$), and
if $u_*$ were a periodic solution, then one would have a linked pair
$(u_*,u_*)$ and $\uu^1$ to which one could apply Theorem~\ref{A}.
Since $u_*$ is a saddle-focus it does not perturb to a
periodic solution for $E>0$. Hence we need to use a different
regularization which conveys the information 
that $u_*$ acts as a minimum. The form of the perturbation that
 we have in mind is depicted in Fig.~\ref{f:onesadfoc}, 
where we have drawn the ``potential'' $L(u,0,0)$.

This idea can be formalized as follows.
Choose a function $T \in C^\infty_0[0,\infty)$ such that $0 \leq T(s) \leq 1$,
$T(s)=1$ for $x \leq \frac{1}{2}$, $T(s)$ strictly decreases on
$(\frac{1}{2},1)$ and $T(s)=0$ for $x \geq 1$. 
Add a perturbation 
$$
  \Phi_\eps(u)= \int_{u_*}^{u} - 2 C_0 \, (s-u_*) \,
  T \Big(\frac{|s-u_*|}{\epsilon} \Big)    \, ds
$$
to the Lagrangian, i.e.\
$\tilde{L} = L + \Phi_\eps(u)$, where $C_0=\partial^2_u L(u_*,0,0)$.  
The new Euler-Lagrange equation
 near $u_*$ becomes
$$
  \partial^2_{u_{xx}} L u_{xxxx} + 
  \big[2 \partial^2_{u_{xx}u} L - 
  \partial^2_{u_x} L  \big] u_{xx} + 
  \partial^2_{u} L 
  \textstyle \left[ 1-2T \big( \frac{|u-u_*|}{\epsilon}\big) 
  \right] (u-u_*) = O(U^2),
$$
where all partial derivatives of $L$ are evaluated at $(u_*,0,0)$, and where
$U$ is the vector $(u-u_*,u_x,u_{xx},u_{xxx})$ in phase space.
Hence for all small $\epsilon$ there are now two additional equilibria near
$u_*$, denoted by $\widehat{u} \in (u_*-\epsilon,u_*-\epsilon/2)$ and 
$\widetilde{u}\in(u_*+\epsilon/2,u_*+\epsilon)$.
Since $(u_*-\widehat{u})-(\widetilde{u}-u_*)= O(\epsilon^2)$, the difference
between $\widetilde{E}(\widehat{u})$ and $\widetilde{E}(\widetilde{u})$ 
is $O(\epsilon^2)$. To level this difference we add another small 
perturbation to $\widetilde{L}$ of the form 
$\Psi(u)= \int_{u_*}^{u} C_\epsilon T(\frac{|s-u_*|}{2\epsilon}) ds$,
i.e. $\widehat{L}(u) = \widetilde{L}(u) + \Psi(u)$, where $C_\epsilon$ 
is chosen so that $\hat{E}(\widehat{u})=\hat{E}(\widetilde{u})$
(of course $\widehat{u}$ and $\widetilde{u}$ shift slightly),
and $C_\epsilon = O(\epsilon^2)$.
Using the same analysis as before  we conclude
that a neighborhood of $u_*$ in the energy level $E(\widehat{u})$
looks just like Fig.~\ref{f:twosimple}.
In $B=\{ (u_1,u_2) \,|\, u^*_1<u_1<\widehat{u} ,~\widetilde{u}<u_2<u^*_2 \}$ 
we find a minimum. Choose an regular energy level $E_\eps$ slightly 
larger than $\widehat{E}(\widehat{u})=\widehat{E}(\widetilde{u})$ 
(with $E_\epsilon = O(\epsilon))$, such that the minimum in $B$ persists.
Taking this minimum and the original $\uu^1$ --- which persists since we 
have only used small perturbations, preserving  $D$ 
(see Fig.~\ref{f:onesadfoc}) as an isolating neighborhood --- we apply 
Theorem~\ref{A}.

Finally, we take the limit $\eps \to 0$. The solutions now converge 
to solutions of the original equation in the degenerate energy level.
It follows that in the energy level $E=0$ a solution of type $(p,q)$ exists, 
where the number of extrema has to be counted with multiplicity since extrema
can coalesce in pairs at $u_*$.

%
%

\section{Computation of the homotopy index}\label{XI}\label{compu}

Theorems~\ref{A}, \ref{D}, and \ref{F} hang on the homology computations of 
the homotopy invariant for certain canonical braid classes (Lemmas 
\ref{comp1}, \ref{comp2}, and \ref{comp3}). Our strategy (as in, 
e.g., \cite{Angenent1}) is to choose a sufficiently simple system 
(an integrable Hamiltonian system) which exhibits the braids 
in question and to compute the homotopy index via knowing the structure
of an unstable manifold. By the topological invariance of the homotopy index, 
any computable case suffices to give the index for any period $d$. 

\begin{figure}[hbt]
\begin{center}
\includegraphics[angle=0,height=4cm,width=9cm]{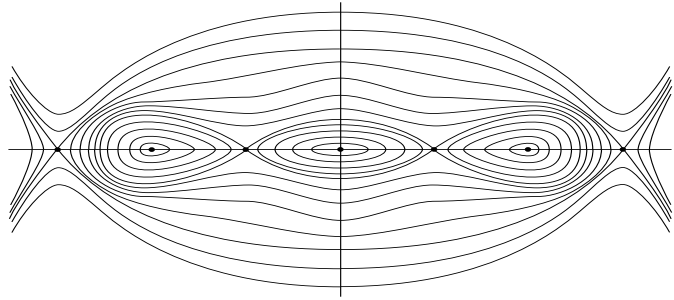}
\caption{The integrable model in the $(u,u_x)$ plane; there are centers
at $0, \pm 2$ and saddles at $\pm 1, \pm 3$.}
\label{intI}
\end{center}
\end{figure}

Consider the first-order Lagrangian system given by the Lagrangian
$L_\lambda (u,u_x) = \frac{1}{2} |u_x|^2 + \lambda F(u)$, where we choose
$F(u)$ to be an even four-well potential, with $F''(u) \ge -1$, and $F''(0) =-1$.
The Lagrangian system $(L_\lambda,dx)$ defines an integrable
Hamiltonian system on $\R^2$, with phase portrait given in Fig. \ref{intI}.

Linearization about bounded solutions $u(x)$ of the above Lagrangian system
yields the quadratic form 
$$
Q[\phi] = \int_0^1 |\phi_x|^2 dt + \lambda \int_0^1 F''(u(x)) \phi^2 dx \ge
\int_0^1 ( \pi^2 -\lambda)  \phi^2 dx,\quad \phi \in H_0^1(0,1),
$$
which is strictly positive for all $0<\lambda<\pi^2$.
For such choices of $\lambda$ the time-1 map defined via the  
induced Hamiltonian flow $\psi^x$, i.e., 
$(u,p_u)=(u,u_x) \mapsto \psi^1(u,p_u)$, 
is an area preserving monotone twist map.
The generating function of the twist map is given by the minimization problem
\begin{eqnarray*}
S_\lambda (u_1,u_2) = \inf_{q\in X(u_1,u_2)} \int_0^1 L_\lambda (u,u_x)dx,
\end{eqnarray*}
where $X(u_1,u_2) = \{u \in H^1(0,1)~|~u(0)=u_1,~u(1) = u_2\}$.\footnote{
The strict positivity of the quadratic form $Q$ via the choice of $\lambda$ 
yields a smooth family of hyperbolic minimizers.}
The function $S_\lambda $ is a smooth function on $\R^2$, with
$\partial_1 \partial_2 S_\lambda  >0$.
The recurrence function $\RR_\lambda (u_{i-1},u_i,u_{i+1}) =
\partial_2 S_\lambda (u_{i-1},u_i) + \partial_1 S_\lambda (u_i,u_{i+1})$ 
satisfies Axioms (A1)-(A3), and thus defines an
exact (autonomous) parabolic recurrence relation on ${\bf X} = \R^\Z$.
We choose the potential $F$ such that the bounded solutions within the 
heteroclinic loop between $u=-1$ and $u=+1$ have the property that the 
period $T_\lambda $ is an increasing function of the amplitude $A$, 
and $T_\lambda(A) \to \frac{2\pi}{\sqrt{\lambda}}$, as $A\to 0$.

This single integrable system is enough to compute the homotopy index
of the three families of braid classes in 
Lemmas  \ref{comp1}, \ref{comp2}, and \ref{comp3} in \S\ref{X}.

We begin by identifying the following periodic solutions.
Set $\vv^{1,\pm} = \{v_i^{1,\pm}\}$, $v_i^{1,\pm} = \pm 3$, and 
$\vv^{2,\pm} = \{v_i^{2,\pm}\}$, $v_i^{2,\pm} = \pm 1$.
Let $\widehat u(t)$ be a solution of $(L_\lambda,dx)$ with 
$\widehat u_x(0) = 0$ (minimum), $ |\widehat u(x)| <1$, and 
$T_1(A(\widehat u)) = 2\tau_0 >2\pi$, $\tau_0 \in \N$.
For arbitrary $\lambda \le 1$ this implies that
$$
	T_\lambda(A(\widehat u)) = \frac{T_1(A(\widehat u))}{\sqrt{\lambda}} 
	= \frac{2\tau_0}{\sqrt{\lambda}},
$$
where we choose $\lambda$ so that $\frac{1}{\sqrt{\lambda}} \in \N$. 
For $r\ge 1$  set $d:=\frac{\tau_0 r}{\sqrt{\lambda}}$ and  
define $\vv^3 := \{v_i^3\}$, with $v_i^3 = \widehat u (i)$, and 
$\vv^4 = \{v_i^4\}$, with $v_i^4 = \widehat u (i+ \tau_0/\sqrt{\lambda})$, 
$i=0,...,d$. Clearly, $\intnum(\vv^3,\vv^4)= r$, for all 
$\frac{1}{\sqrt{\lambda}} \in \N$.

Next choose $\widetilde u(x)$, a solution of $(L_\lambda,dx)$ with 
$\widetilde u_x(0) = 0$ (minimum), which oscillates around both
equilibria $-2$ and $+2$, and in between the equilibria $-3$ and $+3$, and 
$T_1(A(\widetilde u)) = 2\tau_1 > 2\pi$, $\tau_1 \in \N$. As before
$$
	T_\lambda(A(\widetilde u)) = \frac{T_1(A(\widetilde u))}{\sqrt{\lambda}} 
	= \frac{2\tau_1}{\sqrt{\lambda}}.
$$
Let $2p\ge r$ and choose $\tau_0, \tau_1 \ge 4$ such that
$$
\frac{\tau_0}{\tau_1} = \frac{2p}{r} \ge 1 \quad\quad\quad (\tau_0 \ge \tau_1).
$$
Set $\vv^5 = \{v_i^5\}$, $v^5_i = \widetilde u (i)$, and $\vv^6 = \{v_i^6\}$, 
with $v_i^6 = \widetilde u (i+ \tau_1/\sqrt{\lambda})$, $i=0,...,d$.
For $x \in [0,d]$ the solutions $\widehat u$ and $\widetilde u$ have exactly 
$2p$ intersections. Therefore, if we choose $\lambda$ sufficiently small, 
i.e. $\frac{1}{\sqrt{\lambda}} \in \N$ is large, then it also holds that 
$\intnum(\vv^{3,4}, \vv^{5,6}) = 2p$.

Finally we choose the unique periodic solution $u(x)$, with 
$|u(x)|<1$, $u_x(0)=0$ (minimum), and $T_1(A(u)) = 2\tau_2>2\pi$, 
$\tau_2 \in \N$. Let $0<2q<r\le 2p$, and choose $\tau_2$, and 
consequently the amplitude $A$, so that
$$
	\frac{\tau_0}{\tau_2} = \frac{2q}{r} < 1 \quad\quad\quad 
	(\tau_0 < \tau_2, \quad A(\widehat u)) < A(u)).
$$
The solution $u$ is part of a hyperbolic circle of solutions 
$u_s(x)$, $s \in \R/\Z$. Define $(\uu(s))_{s \in \R/\Z}$, with 
$\uu(s) = \{ u_i(s)\}$, where $u_i(s) = u_s(i+2\tau_2 s/\sqrt{\lambda})$.
As before, since the intersection number of $\widehat u$ and $u_s$ is equal 
to $2q$, it holds that $\intnum(\uu(s),\vv^{3,4}) =2q$, for $\lambda$ 
sufficiently small. Moreover, $\intnum(\uu(s),\vv^{5,6}) = 2p$. From 
this point on $\lambda$ is fixed.
We now consider three different skeleta $\vv$.

I: $\vv= \{\vv^{2,-},\vv^{2,+},\vv^3,\vv^4\}$.  
The relative braid class $[\uu~\rel~\vv]_{\rm I}$
is defined as follows: $v_i^{2,-} \le u_i \le v_i^{2,+}$, 
and $\uu$ links with the strands $\vv^3$ and
$\vv^4$ with intersection number $2q$, $0\le 2q<r$.
The topological class $\bigl\{\uu~\rel~\{\vv\}\bigr\}$
is precisely that of Lemma~\ref{comp1} [Fig.~\ref{StrI}]
and as such is bounded and proper.

II: $\vv= \{\vv^{2,-},\vv^{1,+},\vv^3,\vv^4,\vv^5\}$.
The relative braid class $[\uu~\rel~\vv]_{\rm II}$
is defined as follows: $v_i^{2,-} \le u_i \le v_i^{1,+}$, 
$\uu$ links with the strands $\vv^3$ and 
$\vv^4$ with intersection number $2q$, $0\le 2q<r$, 
and $\uu$ links with $\vv^5$ with intersection number $2p$. 
The topological class $\bigl\{\uu~\rel~\{\vv\}\bigr\}$ 
is precisely that of Lemma~\ref{comp3} [Fig.~\ref{StrIII}]
and as such is bounded and proper.

\begin{figure}[hbt]
\begin{center}
\psfragscanon
\psfrag{D}[][]{\Large $\Dual$}
\includegraphics[angle=0,height=1.5in,width=5.1in]{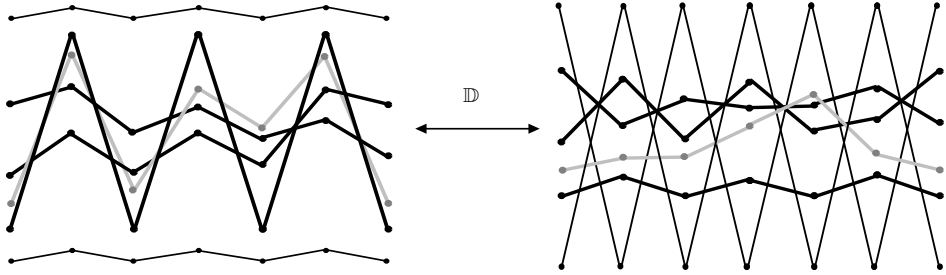}
\caption{The augmentation of the braid from Lemma~\ref{comp2} [left]
is the dual of the type III braid [right].}
\label{flipit}
\end{center}
\end{figure}

III: $\vv= \{\vv^{2,-},\vv^3,\vv^4,\vv^5,\vv^6\}$.
The relative braid class $[\uu~\rel~\vv]_{\rm III}$
is defined as follows: $v_i^{2,-} \le u_i$, 
$\uu$ links with the strands $\vv^3$ and 
$\vv^4$ with intersection number $2q$, 
and $\uu$ links with $\vv^5$ and $\vv^6$ with intersection number $2p$.
The topological class $\bigl\{\uu~\rel~\{\vv\}\bigr\}$ 
is {\it not} bounded [Fig.~\ref{flipit}[right]]. 
The augmentation of this braid class is bounded.

\vsp
\noindent
{\bf Cases I and II:} Since the topological classes are bounded and proper, 
the invariant $\hh$ is independent of period of the chosen representative, 
and can be easily computed from the integrable model. 
The closure of the collection of topologically equivalent braid classes 
is an isolating neighborhood for the parabolic flow $\Psi^t$ induced by the
recurrence relation $\RR_\lambda=0$ (defined via $(L_\lambda,dx)$).
The invariant set is given by the normally hyperbolic circle 
$\{\uu(s)\}_{s\in \R/\Z}$. For this reason the index $\hh$ can be 
computed via the connected component that contains
the critical circle; we denote this neighborhood by $N$. 
The Conley index of $N$ can be determined 
via computing $W^u(\{\uu(s)\})$, the unstable manifold associated to this 
circle. This computation is precisely 
that appearing in the calculations of \cite[pp. 372]{Angenent1}:
$W^u(\{\uu(s)\})$ is orientable and of dimension $2q$, and thus
\begin{equation}
\hh(\uu~\rel~\vv)=h(N)
\simeq \left(S^1 \times S^{2q-1}\right) / 
	\left(S^1 \times \{{\rm pt}\}\right) 
\simeq S^{2q-1}\vee S^{2q}.
\end{equation}
The Conley homology is given by 
$CH_k(\hh) = \R$ for $k=2q-1,2q$, and $CH_k(\hh) = 0$ elsewhere.
This completes the proofs of the Lemmas \ref{comp1} and \ref{comp3}.
\fp

\vsp
\noindent
{\bf Case III:} It holds that 
$$
\bigl\{\uu~\rel~\vv\bigr\}\cap\left(\Conf^1_{2p}~\rel~\vv\right)
\neq \emptyset.
$$
The discrete class for period $2p$ is bounded, but for periods $d>2p$ this 
is not the case. However, by augmenting the braid, we obtain from 
\rmref{extra} that
$$
	\hh(\uu~\rel~\vv) = \hh(\uu~\rel~\vv^*) ,
$$
where $\vv^* = \vv \cup \{\vv^{1,-},\vv^{1,+}\}$.
Since the topological class $\bigl\{\uu~\rel~\{\vv^*\}\bigr\}$ is 
bounded and proper, we may use the previous calculations to conclude that
$$
\hh(\uu~\rel~\vv^*)
\simeq (S^1 \times S^{2q-1})/(S^1 \times \{{\rm pt}\})
\simeq S^{2q-1}\vee S^{2q}.
$$
Our motivation for this computation is to complete the proof of 
Lemma~\ref{comp2}. Let $[\uu'~\rel~\vv']$ denote the period $2p$ 
braid class described by Fig. \ref{StrII}, with intersection numbers
denoted by $2q'$ and $2r'$, and 
let $[\uu~\rel~\vv]$ denote a type-III braid of period $2p$. Then,
it is straightforward to see (as illustrated in Fig.~\ref{flipit}) that, 
for $q'=p-q$ and $r'=2p-r$,
\begin{equation}
	\bigl[\uu'~\rel~[\left.\vv'\right.^*]\bigr] 
	= \Dual\Bigl(\bigr[\uu~\rel~[\vv]\bigr]\Bigr).
\end{equation}
Lemma \ref{comp2} gives the index for the augmented class 
$\bigl\{\uu'~\rel~\{\left.\vv'\right.^*\}\bigr\}$,
which is bounded and proper as a topological class.
The above considerations allow us to compute the homology of 
$\hh(\uu'~\rel~\left.\vv'\right.^*)$ via Theorem~\ref{thm_duality}:
\begin{equation}
\begin{array}{ccl}
	  	CH_*\left(\hh(\uu'~\rel~\left.\vv'\right.^*)\right)
	&\cong& CH_*\left(\hh(\Dual\uu~\rel~\Dual\vv)\right)	\\
	&\cong& CH_{2p-*}\left(\hh(\uu~\rel~\vv)\right) \\
        &\cong& CH_{2p-*}\left(\hh(\uu~\rel~\vv^*)\right) \\
	&\cong& \left\{	
		\begin{array}{ccl}
		\R &:& 2p-* = 2q-1,2q \\
		0 &:& {\mbox{ else}}
		\end{array}\right. \\
	&\cong&	\left\{
		\begin{array}{ccl}
		\R &:& *=2q',2q'+1 \\
		0 &:& {\mbox{ else}}
		\end{array} \right.
\end{array}. 
\end{equation}
The intersection numbers $2q'$ and $2r'$ are exactly 
those of Lemma \ref{comp2}, completing the proof.
\fp

%
%

\section{Postlude}\label{postlude}
\subsection{Extensions and questions: dynamics}
There are several ways in which the basic machinery introduced in 
this paper can be generalized to other dynamical systems.

\noindent {\em 1. Scalar uniformly parabolic PDEs.}
Theorem~\ref{topinvariant} suggests strongly that the 
homotopy index for discretized braids extends to and agrees 
with an analogous index for parabolic dynamics on spaces of 
smooth curves via uniformly parabolic PDE's. This is true \cite{AGV}.

\noindent {\em 2. Periodicity in the range.}
Although we consider the anchor points to be in $\R$, one may just
as well constrain the anchor points to lie in $S^1$ and work in the 
universal cover. Such additional structure is used in the 
theory of annulus twist maps 
\cite{Angenent1,Aubry-LeD,LeCalvez,LeCalvez2,Mather}.
All of our results immediately carry over to this setting. We note that
compact-type boundary conditions necessarily follow. 
\vskip.1cm

\noindent {\em 3. Aperiodic dynamics.}
One can also extend the theory to 
include braid diagrams with ``infinite length'' strands. To be more 
precise, consider braid diagrams on the infinite 1-d lattice, and omit 
the spatial periodicity.
In this context several compactness issue come into play. To name a few:
(a) the parabolic flows generated by aperiodic recurrence relations no 
longer live on a finite dimensional space but on the infinite dimensional 
space $\ell^\infty(\R)$. See \cite{Angenent1} for a case similar to this;
(b) a priori, the Conley index should be replaced with an infinite 
dimensional analogue such as that developed by Rybakowski \cite{Ryb}. 
However, if one considers braid diagrams with finite word metric, 
the stabilization theory of \S\ref{stable} allows one to define 
the necessary invariants via the finite dimensional theory in this paper.
This is not unlike the procedure one can follow in the treatment of 
parabolic PDE's \cite{AGV}.
\vskip.1cm

\noindent {\em 4. Fixed boundary conditions.}
Our decision to use {\em closed} braid diagrams is motivated by 
applications in Lagrangian systems; however, one can also fix the end 
points of the braid diagrams. In this setting one can define a braid 
invariant in the same spirit as is done for closed braids.
The proof of stabilization is not sensitive to the type of boundary 
conditions used. Such an extension of the theory to include fixed 
endpoints is useful in applications to parabolic PDE's \cite{AGV}.
\vskip.1cm

\noindent {\em 5. Traveling waves and period orbits.} 
The stationary solutions we find in this
paper are but the beginnings of a dynamical skeleton for the systems 
considered. The next logical step would be to classify connecting
orbits between stationary solutions: several authors (\eg, \cite{Mallet-Paret})
have considered these problems analytically in the context of 
traveling wave phenomena in monotone lattice dynamics. There is 
a precedent of using Conley indices to prove existence theorems 
for connecting orbits (\eg, \cite{Mischaik}): we
anticipate that such applications are possible in our setting. 
One could as well allow the skeletal strands to be 
part of a periodic motion (in the case of non-exact recurrence 
relations). In this setting one could look for both fixed points and 
periodic solutions of a given braid class.
\vskip.1cm

\noindent {\em 6. Long-range coupling.}
Assume that the recurrence relations $\RR_i$ are functions of the form 
$\RR_i(u_{i-n},\ldots,u_{i+n})$ for some $n$. Even if a strong 
monotonicity condition holds, $\partial_j\RR_i>0$ for all $j\neq i$ 
the proof of Proposition~\ref{word} still encounters a difficulty: 
two strands with a simple (codimension-one) tangency can have enough
local crossings to negate the parabolic systems' separation. 
Such monotone systems do exhibit an ordering principle \cite{Angenent2, Hir}
(initially nonintersecting strands will never intersect), but additional
braiding phenomena is not automatically present.
\vskip.1cm

\noindent {\em 7. Higher-dimensional lattice dynamics.}
In parabolic PDE's of spatial dimension greater than one, the straightforward
generalization of the lap number (number of connected components of an
intersection of graphs) does not obey a monotonicity property
(due to the fact that for graphs of $\R^n$ with $n>1$, critical 
points of non-zero index and co-index can pass through each other).
Finding a suitable form of dynamics which retains some isolation 
remains an important and challenging problem.
\vskip.1cm

\noindent {\em 8. Arbitrary second-order Lagrangians.}
Our principal dynamical goal is to prove existence theorems for
periodic orbits with a minimal amount of assumptions, particularly
``genericity'' assumptions (which are, in practice, rarely verifiable). 
To this end, we have been successful for second-order Lagrangians
modulo the twist assumption. Although this assumption is provably
satisfied in numerous contexts, we believe that it is not,
strictly speaking, necessary. Its principal utility is in the 
reduction of the problem to a finite-dimensional recurrence 
relation. We believe that the forcing results proved 
in \S\ref{second} are valid for {\em all} second-order Lagrangian
systems. See for instance \cite{KV} for a result on that behalf.
We propose that a version of the curve-shortening 
techniques in the spirit of Angenent \cite{Angenent3} should yield a homotopy
index for smooth curves to which our forcing theorems apply.

\subsection{Extensions and questions: topology}
The homotopy index is, as a topological invariant of braid pairs, 
utterly useless. Nevertheless, there is topological meaning 
intrinsic to this index, the precise topological interpretation
of which is as yet unclear. One observes that the index 
captures some Morse-theoretic data about braid classes.
Any topological interpretation is certainly related to linking data of the 
free strands with the skeleton, as evidenced by the examples in this
paper. Though the total amount of linking should provide some upper 
bound to the dimension of the homotopy index, linking numbers alone are 
insufficient to characterize the homotopy index. 

We close with several related questions about the 
homotopy index itself. 
\vskip.1cm

\noindent {\em 1. Realization.} It is clear that given any polynomial in 
$t$, there exists a braid pair whose homological Poincar\'e 
polynomial agrees with this. [Idea: take Example 3 of \S\ref{conley}
and stack disjoint copies of the skeleton vertically, using as many
free copies and strands as necessary to obtain the desired homology.]
Can a realization theorem be proved for the homotopy index itself?
As a first step to this, consider replacing the real coefficients 
in the homological index with integral coefficients. Does torsion
ever occur? We believe not, with the possible exception of a $\Z_2$ torsion. 
\vskip.1cm

\noindent {\em 2. Product formulas and the braid group.}
Perhaps the most pressing problem for the homotopy index is to 
determine a product formula for the concatenation of two braids
with compatible skeleta. This would eliminate the need for 
computing the index via continuation to an integrable model
system as in \S\ref{compu}. 
However, since we work on spaces of {\em closed} braids, a product
formula is not well-defined. The group
structure on the braid group $B_n$ does not extend naturally to 
a group structure on conjugacy classes: where one ``cuts open'' the
braid to effect a gluing can change the resulting braid class
dramatically. The one instance in which a product operation is 
natural is a power of a closed braid. Here, splitting the closed
braid to an open braid and concatenating several copies then 
reclosing yields equivalent closed braids independent of the 
representative of the conjugacy class chosen.
Such a product/power formula, in conjunction with numerical methods
of index computation effective in moderately low dimensions,
would allow one to compute many invariants.
\vskip.1cm

\noindent {\em 3. Improper and unbounded classes}
In certain applications one  also needs to deal with improper
braid classes $[\uu~\rel~\vv]$. To such classes one can also assign 
an index. The interpretation of the index as a Morse theory will not
only depend on the topological data, but also on the behavior of the flow
$\Psi^t$ at $\Sigma^-$. The simplest case is when  
$\Sigma^- \cap \partial N$ consists of finitely many points.
This for example happens when $\uu$ consists of only one strand.
The homotopy index is 
then defined by the intrinsic definition in \rmref{intrinsic}.
The interpretation of the index and the associated Morse theory
depends on the linearization $D\Psi^t|_{\Sigma^- \cap \partial N}$. 
The definition of the index in the case of more complicated sets 
$\Sigma^- \cap \partial N$
and the Morse theoretic interpretation will be subject of future study.
Similar considerations hold for unbounded classes.
\vskip.1cm

\noindent {\em 4. General braids.}
The types of braids considered in this paper are positive braids. 
Naturally, one wishes to extend the ideas to all braids; however, 
several complications arise. First, passing to discretized braids
is invalid --- knowing the anchor points is insufficient data for
reconstructing the braid. Second, compactness is troublesome --- one 
cannot merely bound braid classes via augmentation.
We can model general braids dynamically using recurrence relations with 
nearest neighbor coupling allowing ``positive'', or  ``negative'' interaction.
This idea appears in the work of LeCalvez \cite{LeCalvez,LeCalvez2} 
and can be translated to our setting via a change of variables 
--- coordinate flips ---  of which our duality operator $\Dual$ is 
a particular example. However the compactness and discretization issues remain.
LeCalvez works in the setting of annulus maps, where one can circumvent 
these problems: the general setting is more problematic.
\vskip.1cm

\noindent{\em 5. Hamiltonian vs. Lagrangian.}
One approach to extending to arbitrary braids would be to switch from 
a Lagrangian setting to a Hamiltonian setting. Consider an $S^1$ family
of Hamiltonians $H_t$ on a symplectic surface $(M^2,\omega)$ which has
a ``skeleton'' of periodic orbits braided in $M\times S^1$. 
Adding ``free'' braid strands, one could define a relative Floer 
index for the system which should detect whether the free strands are
forced to exist as periodic orbits.


\appendix

%
%
\section{Construction of parabolic flows}
\label{app_1}
In this appendix, we construct particular parabolic flows
on braid diagrams in order to carry out the continuation 
arguments for the well-definedness of the Conley index
for proper braid diagrams. The constructions are explicit and 
are generated by recurrence relations 
$\RR = (\RR_i)_{i \in \Z}$, with $\RR_{i+d} = \RR_i$, which
are of the form:
\begin{eqnarray}\label{genform}
\RR_i(r,s,t) = a_i(r,s) + b_i(s,t) + c_i(s),\quad i=1,...,d,
\end{eqnarray}
with $a_i, b_i, c_i \in C^1(\R)$, and 
$\frac{\partial a_i}{\partial r}(r,s)>0$,
$\frac{\partial b_i}{\partial t}(s,t)\geq 0$ for all $(r,s,t) \in \R^3$.
By definition, such recurrence relations are parabolic.
\begin{lemma}
\label{explicit}
For any $\vv\in \DD_d^m$ there exists a parabolic flow $\Psi^t$
under which $\vv$ is stationary.
\end{lemma}
{\it Proof.}
In order to have $\Psi^t(\vv)=\vv$, the sequences
$\{\vv_\alpha\}$ need to satisfy 
$\RR_i(v^{\alpha}_{i-1},v^\alpha_{i},v^\alpha_{i+1})=0$, for
some parabolic recurrence relation. We will construct $\RR$ 
by specifying the appropriate functions $\{a_i, b_i, c_i\}$ 
as above. In the construction to follow, the reader should think 
of the anchor points $\{v_i^\alpha\}$ of the fixed braid $\vv$ as constants.

For each $i$ such that the values $\{v_i^\alpha\}_\alpha$ are distinct, 
one may choose $a_i(r,s)=r$, $b_i(s,t)=t$, and $c_i(s)$ to be any
$C^1$ function which interpolates the defined values
\[
	c_i(v^\alpha_i) := -(v^\alpha_{i-1}+v^\alpha_{i+1}) .
\]
This generates the desired parabolic flow. 

In the case where there are several strands $\alpha_1,\alpha_2,\ldots
\alpha_n$ for which $v_i^{\alpha_j}=v^*$ are all equal, the former 
construction is invalid: $c_i$ is not well-defined. According 
to Definition~\ref{PL}, we have for each $\alpha_j\neq\alpha_k$
	$(v_{i-1}^{\alpha_j}-v_{i-1}^{\alpha_k})
	(v_{i+1}^{\alpha_j}-v_{i+1}^{\alpha_k}) < 0$.
This implies that if we order the $\{\alpha_j\}_j$ so that 
$v_{i-1}^{\alpha_1}<v_{i-1}^{\alpha_2}<\cdots<v_{i-1}^{\alpha_n}$, 
then the corresponding sequence $\{v_{i+1}^{\alpha_j}\}_j$ satisfies
$v_{i+1}^{\alpha_n}<v_{i+1}^{\alpha_{n-1}}<\cdots<v_{i+1}^{\alpha_1}$. 
 From Lemma~\ref{stupidlemma} below, there exist increasing functions 
$f$ and $g$ such that 
\[	f(v_{i-1}^{\alpha_j})-f(v_{i-1}^{\alpha_k}) 
= 	g(v_{i+1}^{\alpha_k})-g(v_{i+1}^{\alpha_j}) \,\,\,\, \forall j,k  .
\]
Define $a_i$ and $b_i$ in the following manner: set
$a_i(r,v^*) := f(r)$ and $b_i(v^*,t) := g(t)$. Thus it follows that
there exists a well defined value 
\[	c_i(v^*) := -\left(
		f(v_{i-1}^{\alpha_j}) + g(v_{i+1}^{\alpha_j}) \right)	
\]
which is independent of $j$. For any other strands $\alpha'$, repeat
the procedure, defining the slices $a_i(r,v_i^{\alpha'})$, 
$b_i(v_i^{\alpha'},t)$ and the points $c_i(v_i^{\alpha'})$, 
choosing new functions $f$ and $g$ if necessary. To extend these
functions to global functions $a_i(r,s)$ and $b_i(s,t)$, 
simply perform a $C^1$ homotopy in $s$ without changing the 
monotonicity in the $r$ and $t$ variables: \eg, 
on the interval $[v_i^\alpha,v_i^{\alpha'}]$, choose a 
monotonic function $\xi(s)$ for which $\xi(v_i^\alpha)=
\xi'(v_i^{\alpha})=\xi'(v_i^{\alpha'})=0$ and 
$\xi(v_i^{\alpha'})=1$. Then set 
\[
	a_i(r,s) := (1-\xi(s))a_i(r,v_i^\alpha) + 
			\xi(s)a_i(r,v_i^{\alpha'}) . 
\] 
Such a procedure, performed on the appropriate $s$-intervals, 
yields a smooth $r$-monotonic interpolation. Repeat with $b_i(s,t)$. 
Finally, choose any function $c_i(s)$ which smoothly interpolates the
preassigned values.
These choices of $a_i$, $b_i$ and $c_i$ give the desired recurrence relation,
and consequently the parabolic flow $\Psi^t$.
\fp

\begin{lemma}
\label{stupidlemma}
Given two sequences of increasing real numbers $x_1<x_2<\cdots<x_n$ 
and $y_1<y_2<\cdots<y_n$, there exist a pair of strictly increasing
functions $f$ and $g$ such that 
\begin{equation}
\label{stupid}
	f(x_j)-f(x_k) = g(y_j)-g(y_k) \,\,\,\, \forall j,k
\end{equation}
\end{lemma}
{\it Proof.} Induct on $n$, noting the triviality of the case $n=1$.
Given increasing sequences $(x_i)_1^{N+1}$ and $(y_i)_1^{N+1}$, 
choose functions $f$ and $g$ which satisfy \rmref{stupid} for 
$j,k\leq N$: this is a restriction on $f$ and $g$ only for values
in $[x_1,x_N]$ since outside of this domain the functions can be 
arbitrary as long as they are increasing. Thus, modify $f$ and $g$
outside this interval to satisfy
\[
	f(x_{N+1}) = f(x_N)+C \,\, ; \,\, 
	g(y_{N+1}) = g(y_N)+C ,
\]
for some fixed constant $C>0$. These functions satisfy \rmref{stupid}
for all $j$ and $k$. 
\fp

\begin{lemma}\label{homotopy}
For any pair of equivalent braids $[\uu(0)]=[\uu(1)]$, 
there exists a path $\uu(\lambda)$ in $\Conf^n_d$ and a continuous 
family of parabolic flows $\Psi^t_\lambda$, such that 
$\Psi^t_\lambda(\uu(\lambda)) = \uu(\lambda)$,
for all $\lambda \in [0,1]$.
\end{lemma}

{\it Proof.}
Given $\uu$ any point in $\Conf^n_d$, consider any parabolic 
recurrence relation $\RR_\uu$ which fixes $\uu$ and which is
{\em strictly} monotonic in $r$ and $t$. From the proof of
Lemma~\ref{explicit}, $\RR_\uu$ exists. 
For every braid $\uu'$ sufficiently close to $\uu$, there 
exists $\phi$ a near-identity diffeomorphism of $\Conf^n_d$
which maps $\uu$ to $\uu'$. The recurrence relation 
$\RR_\uu\circ\phi^{-1}$ fixes $\uu'$ and is still parabolic  
since $\phi$ cannot destroy monotonicity. 
Choosing a short smooth path $\phi^t$ of
such diffeomorphisms to $\Id$ proves the lemma on small neighborhoods
in $\Conf^n_d$, which can be pieced together to yield arbitrary
paths. 
\fp

\end{sloppypar}
\end{document}